\documentclass[12pt]{article}
\usepackage{amsmath,amssymb}
\usepackage{graphicx}
\def\date{\hfill \number\day.  \number\month. \number\year }
\def\ledot{\le\!\!\!\raise2pt\hbox{$\scriptscriptstyle\bullet$}\;}
\def\ldot{\!<\!\!\!\raise1.5pt\hbox{$\scriptscriptstyle\bullet$}}
\def\and{\hskip2pt{\scriptstyle \land}\hskip2pt}

\def\smcap{\kern 2pt {\scriptstyle \cap}\kern 2pt}
\def\smcup{\kern 2pt {\scriptstyle \cup}\kern 2pt}
\def\ssm{{\smallsetminus}}

\def\qed{\hglue 0pt plus 1filll $\scriptstyle\square$}  

\def\cB{{\cal B}}
\def\cC{{\cal C}}
\def\cD{{\cal D}}
\def\cE{{\cal E}}
\def\cF{{\cal F}}
\def\cG{{\cal G}}
\def\cH{{\cal H}}

\def\cK{{\cal K}}

\def\cM{{\cal M}}

\def\cO{{\cal O}}
\def\cP{{\cal P}}

\def\cS{{\cal S}}

\def\1{1\kern-2.5pt {\rm l}}    

\def\Cs#1#2{{\rm Cs}_{\hskip1pt #1}{\hskip.5pt#2}} 

\def\Aut{\mathop{{\rm Aut}}} 
\def\GL#1#2{{\rm GL}_{#1}{#2}}

\def\SL#1#2{{\rm SL}_{#1}{#2}}  
\def\PSL#1#2{{\rm PSL}_{#1}{#2}} 

\def\Opr#1#2{{\rm O}^\prime_{#1}{#2}}  
\def\POpr#1#2{{\rm PO}^\prime_{#1}{#2}}
\def\SO#1#2{{\rm SO}_{#1}{#2}}

\def\Spin#1#2{{\rm Spin}_{#1}{#2}}   
\def\U#1#2{{\rm U}_{#1}{#2}}
\def\SU#1#2{{\rm SU}_{#1}{#2}} 
\def\PU#1#2{{\rm PU}_{#1}{#2}}
\def\PSU#1#2{{\rm PSU}_{#1}{#2}}  
\def\Sp#1#2{{\rm Sp}_{#1}{#2}}
\def\PSp#1#2{{\rm PSp}_{#1}{#2}}

\def\Gtwo{{\rm G}_2}

\def\SaU#1{{\rm U}^\alpha_{#1}\HH}

\def\lie{\mathop{\strut\frl}}

\font\Bbb=msbm10 at12 pt
\font\scriptBbb=msbm7  at12 pt
\font\scriptscriptBbb=msbm5 at12 pt

\newfam\Bbbfam
\textfont\Bbbfam=\Bbb
\scriptfont\Bbbfam=\scriptBbb
\scriptscriptfont\Bbbfam=\scriptscriptBbb

\def\CC{{\fam=\Bbbfam C}}

\def\HH{{\fam=\Bbbfam H}}

\def\OO{{\fam=\Bbbfam O}}

\def\RR{{\fam=\Bbbfam R}}
\def\sS{{\fam=\Bbbfam S}}
\def\TT{{\fam=\Bbbfam T}}

\def\ZZ{{\fam=\Bbbfam Z}}

\font\sansserif=cmss10 at 12 pt 
\font\scriptsansserif=cmss10 at 8.4 pt
\font\scriptscriptsansserif=cmss10 at 6 pt
\newfam\ssfam
\textfont\ssfam=\sansserif
\scriptfont\ssfam=\scriptsansserif
\scriptscriptfont\ssfam=\scriptscriptsansserif

\let\usuDelta=\Delta  
\let\usuGamma=\Gamma  
\let\usuLambda=\Lambda  
\let\usuOmega=\Omega  
\let\usuPhi=\Phi  
\let\usuPi=\Pi  
\let\usuPsi=\Psi  
\let\usuSigma=\Sigma  
\let\usuTheta=\Theta  
\let\usuUpsilon=\Upsilon  
\let\usuXi=\Xi  

\def\Alpha{{\fam=\ssfam A}}
\def\Beta{{\fam=\ssfam B}}
\def\Chi{{\fam=\ssfam \usuChi}}
\def\Delta{{\fam=\ssfam \usuDelta}}
\def\Epsilon{{\fam=\ssfam E}}
\def\Eta{{\fam=\ssfam H}}
\def\Gamma{{\fam=\ssfam \usuGamma}}
\def\Chi{{\fam=\ssfam X}}

\def\Kappa{{\fam=\ssfam K}}
\def\Lambda{{\fam=\ssfam \usuLambda}}
\def\Mu{{\fam=\ssfam M}}
\def\Nu{{\fam=\ssfam N}}
\def\Omega{{\fam=\ssfam \usuOmega}}
\def\Phi{{\fam=\ssfam \usuPhi}}
\def\Pi{{\fam=\ssfam \usuPi}}
\def\Psi{{\fam=\ssfam \usuPsi}}

\def\Rho{{\fam=\ssfam P}}
\def\Sigma{{\fam=\ssfam \usuSigma}}
\def\Tau{{\fam=\ssfam T}}
\def\Theta{{\fam=\ssfam \usuTheta}}
\def\Upsilon{{\fam=\ssfam \usuUpsilon}}
\def\Ypsilon{{\fam=\ssfam \usuUpsilon}}
\def\Xi{{\fam=\ssfam \usuXi}}
\def\Zeta{{\fam=\ssfam Z}}

\let\epsilon=\varepsilon
\let\theta=\vartheta
\let\phi=\varphi
\let\rho=\varrho

\font\teneufm=eufm10 at 12 pt
\font\seveneufm=eufm7 at 12 pt
\font\fiveeufm=eufm5 at 12 pt
\newfam\eufmfam
\textfont\eufmfam=\teneufm
\scriptfont\eufmfam=\seveneufm
\scriptscriptfont\eufmfam=\fiveeufm
\def\frak#1{{\fam\eufmfam\relax#1}}

\def\frB{{\frak B}}

\def\frL{{\frak L}}

\def\frS{{\frak S}}

\def\fre{{\frak e}}

\def\frl{{\frak l}}

\def\frx{{\frak x}}

\def\frz{{\frak z}}

\input cyracc.def

\font\tencyss=wncyss10 scaled1200
\def\cyss{\tencyss\cyracc}

\def\sm{\hskip-1pt$\raise1pt\hbox{$\scriptstyle\setminus$}$\hskip-1pt}

\begin{document}

\abovedisplayskip=3pt
\belowdisplayskip=3pt

\overfullrule=1pt
\font\bf=cmbx10 scaled 1240
\font\bbf=cmbx10 scaled 1750
\font\Bf=cmbx10 scaled 1500
\font\sbf=cmbx12 scaled 970
\font\ss=cmss12
\def\bullett{\raise1pt\hbox{$\scriptscriptstyle\bullet$}} 
\def\ssm{\smallsetminus}
\def\rk{\rm rk\,}
\let\hat=\widehat
\let\tilde=\widetilde
\let\bold=\bf
\let\§=\ss{}
\def\bib{\bibitem{} }
\font\teneufm=eufm10 scaled 1200 
\font\seveneufm=eufm10 
\newfam\eufmfam
\textfont\eufmfam=\teneufm
\scriptfont\eufmfam=\seveneufm

\def\frak#1{{\fam\eufmfam\relax#1}}
\def\frL{{\frak L}}
\def\fre{{\frak e}}

\def\3ast{$(\lower1pt\hbox{$\scriptstyle**$})\hskip-13pt\raise3.2pt\hbox{$\scriptstyle*$}\hskip6pt$}
\def\ldot{\!<\!\!\!\raise1.3pt\hbox{$\scriptscriptstyle\bullet$}}
\def\Rtimes{\times\hskip-3pt\raise0pt\hbox{$\vrule height 6pt width .8pt$}}
\def\smotimes{\begin{scriptsize}$\otimes$\end{scriptsize}}
\def\smoplus{\begin{scriptsize}$\oplus$\end{scriptsize}} 
\def\tsmotimes{\begin{scriptsize}$\tilde\otimes$\end{scriptsize}}
\def\tsmoplus{\begin{scriptsize}$\tilde\oplus$\end{scriptsize}} 
\def\GGamma{\Gamma\hskip-4.6pt\Gamma}
\def\Yu{{\cyss Yu}}

\centerline{\bbf  Compact 16-dimensional planes. An  update} 

\par\bigskip
\hspace{160pt} {by Helmut R. \sc Salzmann} \hfill 
\par\bigskip
\begin{abstract}
This paper is an addition to the book  \cite{cp} on {\it Compact projective planes\/}. Such planes, if 
connected and finite-dimensional, have a point space of topological dimension 2, 4, 8, or 16, 
the classical example in the last case being the projective closure of the affine plane over the 
octonion algebra. The final result in the book (which was published 20 years ago) is a complete
description of all planes admitting an automorphism group of dimension at least 40. Newer results on 8-dimensional planes have been collected in \cite{sz1}. Here, we present a classification of 
16-dimensional planes with a group of dimension ${\hskip-3pt\ge}35$, provided the group does not 
 fix exactly one flag, and we prove several further theorems, among them  criteria for a 
 connected group of automorphisms to be a Lie group.  My sincere thanks are due to 
 Hermann H\"ahl for many fruitful discussions.
\end{abstract}
\par\bigskip
{\Bf 1. Introduction}
\par\medskip
The last section {\it ``Principles of classification''\/} of the treatise {\it Compact projective planes\/} 
\cite{cp} consists of a detailed programme how to determine all sufficiently homogeneous compact 
$16$-dimensional topological planes. Here, we shall survey and amplify  the results that have been achieved in the meantime. For $2$- and $4$-dimensional planes see 
\cite{cp} \S\S\ 3 and 7; newer developments in the $8$-dimensional case are summarized in 
\cite{sz1}. It is an open problem if the dimension of a compact projective plane is necessarily finite;  
the point space of any finite dimensional compact plane has necessarily dimension $0$ or 
$2^m$ with  $m{\,\in\,}\{1,2,3,4\}$,\; see \cite{cp} 41.10 and 54.11. 
\par\smallskip
The classical $16$-dimensional plane $\cO$ is the projective closure of the affine plane over 
the (bi-associative) division algebra $\OO$ of the (real) octonions; for the properties of $\OO$ see 
2.5 below or \cite{cp} \S\,1,  cf. also \cite{eb} and 
\cite{es}. The  automorphism group  
$\Aut\cO$ is isomorphic to the $78$-dimensional exceptional simple Lie group ${\rm E}_6(-26)$, it is transitive on the set of quadrangles of $\cO$ (cf. \cite{cp} 17.\,6,\,7,\,10). In fact, $\cO$ is the only compact $16$-dimensional plane  admitting a point-transitive group $\Delta$ of automorphisms (L\"owen \cite{lw1}). 
Such a group $\Delta$ necessarily contains the compact flag-transitive elliptic motion group of 
$\cO$, see \cite{cp} 63.8.
\par\smallskip
Let $\cP{\,=\,}(P,\frL)$ be a compact projective plane with a point space $P$ of (covering) 
dimension  $\dim P{\,=\,}16$. Taken with the compact-open topology, the automorphism group
$\Sigma{\,=\,}\Aut\cP$ of all continuous collineations is a locally compact transformation group 
of $P$ of finite dimension  $\dim\Sigma$ with a countable basis, see \cite{cp} 44.3 and 83.26, and  
$\dim\Sigma{\,\ge\,}n$ if, and only if, $\Sigma$ contains a euclidean $n$-ball (\cite{cp} 93.\,5 and 6). In all known cases $\Sigma$ is even a Lie group, see also 2.3 below. 
\par\smallskip
Homogeneity of $\cP$ can best be expressed by the size of $\Sigma$ measured by its dimension.
The  most homogeneous planes $\cP$ are known explicitly: if  $\dim\Sigma{\,\ge\,}40$, then
$\cP$ is the projective closure of an affine plane coordinatized by a {\it mutation\/} 
$\OO_{(t)}{\,=\,}(\OO,+,\circ)$ of the octonions, where  
$c{\hskip1pt\circ}z{\,=\,}t{\cdot}cz{\,+\,}(1{-}t){\cdot}zc$ with $t{\,>\,}\frac{1}{2}$; either 
$t{\,=\,}1$, $\,\OO_{(t)}{\,=\,}\OO$, and $\cP{\,\cong\,}\cO$, or $\dim\Sigma{\,=\,}40$ and $\Sigma$ fixes the line at infinity and some point on this line  (\cite{cp} 87.7).
\par\smallskip
The ultimate goal is to describe all pairs $(\cP,\Delta)$, where $\Delta$ is a {\it connected\/} closed
subgroup of $\Aut\cP$ and $\dim\Delta{\,\ge\,}b$ for a suitable bound $b$ in the range
$27{\,\le\,}b{\,<\,}40$. This bound varies with the structure of $\Delta$ and the configuration 
$\cF_\Delta$ of the fixed elements (points and lines) of $\Delta$. In some cases, only $\Delta$ can be determined, but there seems to be no way to find an explicit description of all the corresponding planes.
\par\bigskip
{\Bf 2. Background}
\par\medskip
In the following, $\cP$ will always denote  a compact $16$-dimensional topological projective plane with point space $P$ and line space $\frL$, if not stated otherwise. Note that line pencils  $\frL_p{\,=\,}\{L{\,\in\,}\frL \mid p{\,\in\,}L\}$ are homeomorphic to lines.  Hence the dual of 
$\cP$ is also a compact $16$-dimensional plane.
\par\medskip
{\bf Notation} is standard and agrees with that in the book \cite{cp}. A {\it flag\/} is an incident
point-line pair. A topological plane with a $2$-dimensional point space will be called {\it flat\/}.
For a locally compact group $\Gamma$ and a closed subgroup $\Delta$ the coset space 
$\{\Delta\gamma\mid\gamma{\,\in\,}\Gamma\}$ will be denoted by $\Gamma/\Delta$, 
its dimension $\dim\Gamma{-}\dim\Delta$ by $\Gamma{:\hskip2pt}\Delta$. As in \cite{cp} 94.1, 
the fact that the topological groups $\Gamma$ and $\Delta$ are {\it locally\/} isomorphic\ will sometimes be symbolized by $\Gamma{\,\circeq\,}\Delta$. 
As customary, $\Cs\Delta\Ypsilon$ or just~$\Cs{}\Ypsilon$ is the centralizer of $\Ypsilon$ 
in $\Delta$; the center $\Cs{}\Delta$ is usually denoted by $\Zeta$. Distinguish between
the connected component $\Delta^{\hskip-1.5pt1}$ and the commutator subgroup $\Delta'$ 
of $\Delta$. A {\it Levi complement\/} of the radical $\sqrt\Delta$ is a maximal semi-simple 
subgroup of $\Delta$.  If $M^\Gamma{\,=\,}M$, then $\Gamma|_M$ is the group induced by 
$\Gamma$ on $M$. The configuration of all fixed points  and fixed lines of a subset 
$\Xi{\,\subseteq\,}\Delta$ will be denoted by $\cF_\Xi$. We write $\Delta_{[A]}$ for 
the subgroup of all {\it axial} collineations in $\Delta$ with axis $A$\; (i.e., of collineations fixing $A$ pointwise)  and, dually, $\Delta_{[c]}$ for the subgroup of collineations with center $c$. 
Let $\Delta_{[c,A]}{\,=\,}\Delta_{[c]}{\smcap}\Delta_{[A]}$ and put 
$\Delta_{[C,A]}{\,=\,}\bigcup_{c\in C}\Delta_{[c,A]}$. 
A subgroup $\Gamma{\,\le\,}\Delta$  is called {\it straight\/}, if each point orbit $x^\Gamma$
is contained in some line; by a result of Baer \cite{ba}, the group $\Gamma$ is then contained in 
$\Delta_{[c,A]}$  for some center~$c$ and axis $A$,  or  $\cF_{\Gamma}$ is a Baer subplane (see 2.2). An element  
$\gamma$ is straight, if it generates a straight cyclic group $\langle \gamma\rangle$, 
and $\gamma$  is said to be {\it planar\/}, if  $\cF_{\hskip-1pt\gamma}$ is a subplane.
\par\medskip
{\bf 2.1 Topology.}   {\it Each line $L$ of $\cP$  is homotopy equivalent $(\simeq)$ to an 
$8$-sphere $\sS_8$\/}. \\
This  fundamental theorem is due to L\"owen \cite{Lw}, see also \cite{cp} 54.11. If $L$ is a manifold, then $L$ is homeomorphic ($\approx)$ to $\sS_8\,$ (\cite{cp} 52.3). 
In this case, the spaces $P$ and $\frL$ are homeomorphic (Kramer \cite{kr2}).
So far, no example with 
$L{\,\not\approx\,}\sS_8$ has been found. Furthermore, lines have the following properties: \\
{\it Each line $\,pq\,$ is a locally homogeneous and locally contractible homology manifold, and it has domain invariance\/}, see \cite{cp} 51.12, 54.10, and 51.21. Moreoer, 
$pq\sm\{p,q\}{\,\simeq\,}\sS_7$ by \cite{cp} 52.5 combined with 54.11. \\
The Lefschetz fixed point theorem implies that  each homeomorphism $\phi{\,:\,}P{\,\to\,}P$ has a fixed point. By duality, {\it each automorphism of $\cP$ fixes a point and a line\/}, 
see \cite{cp} 55.\,19,\,45. 
\par\medskip
{\bf 2.2 Baer subplanes.} {\it Each $8$-dimensional $closed$ subplane $\cB$ of   $\cP$ is a Baer subplane\/}, i.e., each point of $\cP$ is incident with a line of 
$\cB$ (and dually, each line of $\cP$ contains a point of $\cB$), see \cite{sz3} \S~\hskip-3pt3 or 
\cite{cp} 55.5 for details.  By a result of L\"owen \cite{lw2}, any two closed Baer subplanes of $\cP$ have a point and a line in common. This is remarkable because a finite Pappian plane of order $q^2$ is a  union of  $q^2{-}q{+}1$ disjoint Baer subplanes, see \cite{sz3} 2.5. Generally, 
$\langle \cM\rangle$ will denote the smallest {\it closed\/} subplane of $\cP$ containing the  set $\cM$ of points and lines. We write $\cB{\,\,\ldot}\cP$ if $\cB$ is a Baer subplane and 
$\cB{\,\ledot\!}\cP$ if $\,8\vert\dim\cB$. 
\par\medskip
{\bf 2.3 Groups.} {\it Any {\rm connected} subgroup $\Delta$ of $\Sigma{\,=\,}\Aut\cP$ with 
$\dim\Delta{\,\ge\,}27$ is a Lie  group\/}, see \cite{psz}. For $\dim\Delta{\,\ge\,}35$, the result has been proved in \cite{cp} 87.1. If $\dim\Sigma{\,\ge\,}29$, then $\Sigma$ itself is a Lie group, cf. 
\cite{sz4}. In particular, $\Delta$ is then either semi-simple, or $\Delta$ has a central torus subgroup or a minimal normal subgroup $\Theta{\,\cong\,}\RR^k$, see \cite{cp} 94.26.
\par\medskip
{\bf 2.4 Coordinates.} Let ${\fre}{\,=\,}(o,e,u,v)$ be a (non-degen\-erate) quadrangle in $\cP$. Then the    affine subplane  $\cP\sm uv$ can be coordinatized with respect to  the 
{\it frame\/}  $\fre$  by a so-called {\it ternary field\/} $K_\tau$, where $K$ is homeomorphic to an affine line, the affine point set is written as $K{\times}K$, lines  are given by an equation 
$y{\,=\,}\tau(s,x,t)$ or $x{\,=\,}c$ (verticals), they are denoted by $[\hskip1pt s,t\hskip1pt]$ or  
$[\hskip1pt c\hskip1pt]$, respectively. Parallels are  either vertical or both have the same 
{\it``slope''\/} $\hskip-2pt s$. The axioms 
of an affine plane can easily~be translated into properties of $\tau$,  see \cite{cp} \S\hskip2pt22. 
Write $e{\,=\,}(1,1)$ and put $x{\,+\,}t{\,=\,}\tau(1,x,t)$. In the classical case  
$\tau(s,x,t){\,=\,}s\hskip1pt x+t$. A~ternary field $K_\tau$ is called a {\it Cartesian field\/} if  $(K,+)$ 
is a group and if identically $\tau(s,x,t){\,=\,}s{\bullett}x+t$ with  multiplication 
$s{\bullett}x{\,=\,}\tau(s,x,0)$; equivalently, the maps $(x,y)\mapsto(x,y{+}t)$ form a transitive group of translations.
Any ternary field $K_\tau$ such that $K$ is homeomorphic to $\RR^8$ and $\tau{\,:\,} K^3{\,\to\,} K$ is continuous yields a compact projective plane $\cP$, 
 cf. \cite{cp} 43.6. 
\par\medskip
{\bf 2.5 Octonions.} The division algebra $\OO$ of the  octonions can be defined as the real 
vector space $\HH^2$ with the multiplication 
$(a,b)(x,y){\,=\,}(ax{-}y\overline b, xb{+}\overline ay)$. Put $\overline{(a,b)}{\,=\,}(\overline a,-b)$.
{\it Conjugation\/}  $c{\,\mapsto\,}\overline c$ is an anti-automorphism of $\OO$; 
as usual, {\it norm\/} and {\it orthogonality\/} are given by $\|c\|^2{\,=\,}c\overline c$ and 
$c{\,\perp\,}d{\,\Leftrightarrow\,}\|c{+}d\|{\,=\,}\|c{-}d\|$. Sometimes it is convenient to write 
$\ell{\,=\,}(0,1)$ and $\OO{\,=\,}\HH{+}\ell\HH$. 
Any two elements of $\OO$ are contained in an associative subfield: $\OO$ is {\it bi-associative\/}. 
The center of $\OO$ consists of the real field $\RR{\,\subset\,}\CC{\times\{0\}}$; 
if $c$ belongs to the orthogonal complement 
$\RR^\perp$, then $c$ is called a {\it pure\/} element. An arbitrary octonion $c$ will also be written 
in the form $c{\,=\,}c_0{\,+\,}\frz$ with $c_0{\,=\,}{\rm re}\,c{\,=\,}\frac{1}{2}(c{\,+\,}\overline c)$.
The group $\Gamma{\,=\,}\Aut\OO$ is the 
$14$-dimensional compact simple Lie group $\Gtwo$, it is transitive on the unit sphere $\sS_6$ in
$\RR^\perp$, and $\Gamma_{\hskip-1.5pt i}$ is transitive on $\sS_5{\,\subseteq\,}\CC^\perp$\ 
(\cite{cp} 11.31--35). 
Obviously, $\lambda{\,:\,}(x,y){\,\mapsto\,}(x,-y){\,\in\,}\Gamma$, and easy verification shows~that 
$$\Cs\Gamma\lambda{\,=\,}\{(x,y){\,\mapsto\,}(x^a,y^ab)\mid 
a,b{\,\in\,} \HH{\,\and\,} a\overline a{\,=\,}b\overline b{\,=\,}1\}{\,\cong\,}\SO4\RR\,.$$
\break
Each involution in $\Gamma$ is conjugate to $\lambda$ \ (see \cite{cp} 11.31d). Note that the
commutator subgroup of $\HH^{\times}$ is  $\HH'{\,=\,}
\{c{\,\in\,}\HH\mid c\overline c{\,=\,}1\}{\,\cong\,}\Spin3\RR{\,\cong\,}\SU2\CC{\,\cong\,}\U1\HH$.
 \par\medskip
{\bf 2.6 Stiffness} refers to the fact that the dimension of the stabilizer $\Lambda$ of a (non-degenerate) quadrangle $\fre$ cannot be very large. The group $\Lambda$ can be interpreted 
as the automorphism group of the ternary field $K_\tau$ defined with respect to $\fre$. The classical plane  is coordinatized by the octonion algebra $\OO$; in this case $\Lambda$ is  
the $14$-dimensional exceptional compact simple Lie group $\Gtwo{\,\cong\,}\Aut\OO$ for 
{\it each\/} choice of $\fre$, cf. \cite{cp} 11.\,30--35. In the general case, let $\Lambda$ denote the connected component of  $\Delta_\fre$. Then  the following holds: \par
\quad (a) {\it $\Lambda{\,\cong\,}\Gtwo$ and $\cF_\Lambda$ is flat, or $\dim\Lambda{\,<\,}14$ \/}\ (\cite{sz6},\,\cite{cp} 83.\,23,\,24), \par
\quad (\^a) {\it $\Lambda{\,\cong\,}\Gtwo$ or $\dim\Lambda{\,\le\,}11$\/}\ (\cite{B1} 4.1), \par
\quad (b) {\it if $\cF_\Lambda$ is a Baer subplane, then $\Delta_\fre$ is compact and 
$\dim\Lambda{\,\le\,}7$\/}\  (\cite{cp} 83.6), \par
\quad (\^b) {\it if $\cF_\Lambda{\;\ldot}\cP$ and $\Lambda$ is a Lie group, then 
$\Lambda{\,\cong\,}\Spin3\RR$ or $\dim\Lambda{\,\le\,}1$\/}\ (\cite{cp} 83.22), \par
\quad (c) {\it if there exists a subplane $\cB$ with  $\cF_\Lambda{\;\ldot}\cB{\;\ldot}\cP$, then 
$\Lambda$ is  compact\/} \ (\cite{sz6} 2.2, \cite{cp} 83.9), \par
\quad (\^c) {\it if $\Lambda$ contains a pair of commuting involutions, then $\Delta_\fre$ is 
compact\/}\ (\cite{cp} 83.10), \par
\quad (d) {\it if $\Lambda$ is compact or semi-simple or if  $\cF_\Lambda$ is connected, then 
$\Lambda{\,\cong\,}\Gtwo$ or $\dim\Lambda{\,\le\,}10$\/}\par
\hskip30pt (\cite{gs} XI.9.8, \cite{sz6} 4.1, \cite{B1} 3.5), \par
\quad (e) {\it if $\Lambda$ is compact, or if $\dim\cF_\Lambda{\,=\,}4$, or if  
$\Lambda$ is a Lie group and $\cF_\Lambda$ is connected,  then \par 
\hskip30pt $\Lambda{\,\cong\,}\Gtwo$ or $\Lambda{\,\cong\,}\SU3\CC$ or $\dim\Lambda{\,<\,}8$\/}\ 
 (\cite{sz6} 2.1; \cite{sz7}, \cite{B1} 3.5; \cite{B2}), \par
\quad (\^e) {\it if $\Lambda$ is a compact Lie group, then $\Lambda{\,\cong\,}\Gtwo,\; \SU3\CC,\; \SO4\RR$, or $\dim\Lambda{\,\le\,}4$\/}\,\ (\cite{sz6} 2.1), \par
\quad (f) {\it  if $\Alpha$ is the automorphism group of a locally compact $8$-dimensional 
{\emph double loop}, then \par\hskip30pt  $\dim\Alpha{\,\le\,}16$\/}\ (\cite{bd}; for the notion of a double loop cf. \cite{gs} XI. \,\S\S\,1,\,8, and 9).
\par\smallskip
{\tt Remark.} The proof of \cite{cp} 83.10 refers to 83.9, which has been stated only for 
{\it connected\/} groups, but 83.8 can be used instead.
\par\medskip

{\bf 2.7  Dimension formula.} By \cite{hal} or  \cite{cp} 96.10, the following holds for the action of 
$\Delta$ on  $P$ or on any closed $\Delta$-invariant subset $M$ of $P$, and for any point 
$a{\,\in\,}M$:
$$\dim\Delta{\,=\,}\dim\Delta_a{\,+\,}\dim a^\Delta \quad {\rm or} Ê\quad 
\dim a^\Delta{\,=\,}\Delta{\,:\,}\Delta_a\,.$$
\par\medskip

{\bf 2.8 Corollary: Fixed configurations.} {\it If $\dim\Delta{\,\ge\,}27$, then $\Delta$ fixes at most a triangle or, up to duality, some collinear points and one further line\/}. \\
Thus, there are only $9$ possibilities for $\cF_\Delta$, again up to duality: \\
(a) $\emptyset$, \quad  (b) $\{W\}$, $W{\in\,}\frL$, \quad (c) {\it flag\/} $\{v,W\}$, $v{\hskip1pt\in\hskip1pt}W$, \quad
(d) $\{o,W\}$, $o{\hskip1pt\notin\hskip1pt}W$, \quad (e) $\{u,v,uv\}$,\\ (f) $\{u,v,w,...,uv\}$,\quad
(g) {\it double flag\/} $\{u,v,ov,uv\}$,\quad (h) $\{u,v,w,...,uv,ov\}$,\quad (i) {\it triangle\/}. \\
These cases will be treated separately.
\par\medskip
{\bf 2.9  Reflections and translations.} 
{\it Let $\sigma$ be a reflection with axis $W$ and center~$c$ in the connected group $\Delta$, and let $\Tau$ denote the  group of translations in $\Delta$ with axis~$W$. If  $W^\Delta{\,=\,}W$  and   $\dim c^\Delta{\,=\,}k{\,>\,}0$, then  
$\,\sigma^\Delta\sigma{\,=\,}\Tau{\,\cong\,}\RR^k$, 
$\,\tau^\sigma{\,=\,}\tau^{-1}$ for each $\tau{\,\in\,}\Tau$, and $k$ is even\/}.
This improves \cite{cp} 61.19b; a detailed {\tt proof} can be found in \cite{sz2} Lemma 2.
\par\medskip
{\bf 2.10 Involutions.} Each involution $\iota$ is either a reflection or it is {\it planar\/} 
($\cF_{\hskip-1.5pt\iota}\ldot\cP$), see \cite{cp} 55.29. Commuting involutions with the same 
fixed point set are identical (\cite{cp} 55.32).
Let $\ZZ_2^{\hskip3pt r}{\,\cong\,}\Phi{\,\le\,}\Aut\cP$ and 
$\dim P{\,=\,}2^m$.  Then $r{\,\le\,}m{+}1$; if $\Phi$ is generated by reflections, then $r{\,\le\,}2$, 
see \cite{cp} 55.34(c,b). If $r{\,\ge\,}m \ ({\,>\,}2)$, then  $\Phi$ contains a  reflection (a planar 
involution), cf. \cite{cp} 55.34(d,b). 
Any torus group in $\Aut\cP$ has dimension at most $m$,\, see \cite{cp} 55.37, for $m{\,=\,}4$ in particular, $\rk\Delta{\,\le\,}4$,  see also~\cite{cs}.
The orthogonal group $\SO5\RR$ cannot act non-trivially on $\cP\,$ (\cite{cp} 55.40). 
\par\medskip
{\bf 2.11 Lemma.} {\it If $\cF_\Delta{\,\ne\,}\emptyset$, then $\Delta$ has no subgroup 
$\Phi{\,\cong\,}(\SO3\RR)^2$\/}.
\par\smallskip
{\tt Proof.} Each of the two simple factors of $\Phi$ contains $3$ pairwise commuting conjugate
involutions. Let them be denoted by $\alpha_\mu$ and $\beta_\nu$, respectively. If one of these 
triples consists of reflections, their centers $c_\mu$ form a triangle. Each fixed line of $\Delta$ is 
incident with one of the centers $c_\mu$;  on the other hand the $c_\mu$ are contained in a $\Phi$-orbit. This contradiction shows that $\Delta$ fixes no line, dually, $\Delta$ has no fixed point. Hence
the $\alpha_\mu$ and the $\beta_\nu$ are planar. Similarly,  all products $\alpha_\mu\beta_\nu$ are planar,  but the group $\langle \alpha_\mu,\beta_\nu\mid \mu,\nu{\,=\,}1,2,3\rangle{\,\cong\,}(\ZZ_2)^4$ must contain  a reflection by \cite{cp} 55.34(d) or 2.10. \qed
\par\smallskip
{\tt Remark.} The assumption $\cF_\Delta{\,\ne\,}\emptyset$ in Lemma 2.11 is indispensable, even in the case of $8$-dimensional planes. In fact, the stabilizer of the real subplane of the classical  quaternion plane $\cP_\HH$  contains 
$\Aut\cP_\RR{\times}\Aut\HH{\,\cong\,}\SL3\RR{\times}\SO3\RR$\: (cf. \cite{cp} 13.6). 
By \cite{cp}~ 86.31, the stabilizer of a quaternion subplane of $\cO$ has a subgroup $\SL3\HH$. 
Hence the stabilizer of a real subplane contains $\SL3\RR$ and $\Aut\OO{\,=\,}\Gtwo$ and their 
direct product (use 2.5). 
\par\medskip
{\bf 2.12 Homologies.} {\it If $W$ is the translation axis of a $16$-dimensional {\rm translation} plane, 
then any group $\Delta_{[z,A]}$ of homologies with center $z{\,\in\,}W$ has dimension at most $1\,;$ more precisely, $\Delta_{[z,A]}$ is discrete or an extension of $e^{\RR}$ by a finite group\/} \ 
(\cite{bh} Th.\,2 and \cite{cp} 25.4). 
{\it In general a connected homology group of a compact connected projective plane is two-ended\/}: 
{\it it is either compact or a direct product of  $e^{\RR}$ with a compact group\/}, cf. \cite{cp}~61.2.

\par\medskip
{\bf 2.13 Lemma.} {\it Each commutative connected subgroup $\Eta{\,\le\,}\Aut\cP$ fixes a point or 
a line$;$ if $\cP$ is flat, then $\Eta$ fixes a point {\rm and} a  line\/}.
\par\smallskip
{\tt Proof by induction on $\dim\cP$.} Let $\1{\,\ne\,}\eta{\,\in\,}\Eta$ and $a^\eta{\,=\,}a$. Then 
$a^\Eta{\,=\,}a$ or $a^\Eta$ is contained in a fixed line of $\Eta$ or 
$\langle a^\Eta\rangle{\,=\,}\cE$ is a subplane of $\cP$. In the last case, $\eta|_\cE{\,=\,}\1$ and
$\cE{\,<\,}\cP$. By induction, $\Eta$ fixes some element of $\cE$, cf. also \cite{gs} XI.10.19. 
If $\dim P{\,=\,}2$, then the lines are homeomorphic to a circle (\cite{cp} 42.10 or 51.29). 
It suffices to show that $\Eta$ has a fixed point, by duality there is also a fixed line. Suppose that 
$\forall_{\hskip-1pt x{\in}P}\;x^\Eta{\,\ne\,}x$. Then $\eta|_{a^\Eta}{\,=\,}\1$,\; $a^\Eta$ is contained 
in some line $L$,  and $\Eta$ is transitive on $L$\; (or the endpoint of an orbit in $L$ would be 
fixed). Thus $\eta$ is an axial collineation with axis $L$ and the unique center of $\eta$ is 
$\Eta$-invariant. \qed
\par\medskip
From \cite{cp} 71.8,10 and 72.1--4 we need the following results:
\par\smallskip
{\bf 2.14 Semi-simple groups of a $4$-dimensional subplane.} 
{\it A semi-simple group $\Upsilon$ of automorphisms of a $4$-dimensional plane $\cD$ is almost 
simple. If $\Upsilon$ is compact and $\dim\Upsilon{\,=\,}3$, the action of $\Upsilon$ on the 
point space is equivalent to the standard action on the classical complex plane, and 
$\cF_\Upsilon{\,=\,}\emptyset$ or $\cF_\Upsilon{\,=\,}\{a,W\}$ with $a{\,\notin\,}W$. 
If $\dim\Upsilon{\,>\,}3$, then $\cD$ is classical and the action is the standard one$;$ 
in particular,  $\Upsilon$ does not fix a flag or exactly one line\/}.
\par\medskip
{\bf Richardson's classification} \cite{cp} 96.34 of compact groups acting on the sphere $\sS_4$ plays a fundamental  r\^ole in the study of $8$-dimensional planes:
\par\smallskip
($\dagger$) {\it  If a compact connected group $\Phi$ acts effectively on the $4$-sphere $S$,
and if $\Phi$ is a Lie group or if $\Phi$ has an orbit of dimernsion ${>\,}1$, then $(\Phi,S)$ is 
equivalent to the obvious standard action of a subgroup of $\SO5\RR$ on $\sS_4$ or 
$\Phi{\,\cong\,}\SO3\RR$  has no fixed point on $S$\/}. 
\par\medskip
Only fragmentary results are known for
compact transformation groups on $\sS_8$. The following will be needed, 
see \cite{cp} 96.\,35 and 36.
\par\smallskip
{\bf 2.15 Compact groups on $\sS_8$.} (a)  {\it Any non-trivial action of a group 
$\Gamma{\,\cong\,}\Gtwo$  on the sphere $\sS_8$ is equivalent to the action of $\Gtwo$ on 
$\OO{\smcup}\{\infty\}$ as  group of  automorphisms of $\OO\,;$ in particular, $\Gamma$ fixes a circle, all other orbits are $6$-spheres\/}. \\
(b) {\it Any non-trivial action of a group $\Omega{\,\cong\,}\Spin7\RR$ on $\sS_8$ has a fixed point 
$\infty$ and is linear on $\RR^8{\,=\,}\sS_8\sm\{\infty\}\,;$ either the action is faithful, $\Omega$ has 
exactly $2$ fixed points,  all other orbits are $7$-spheres, and stabilizers $\Omega_x$ are 
isomorphic to $\Gtwo$,  or  $\Omega$ induces the group $\SO7\RR$, there is a circle of fixed 
points, and all other orbits are $6$-spheres\/}.
\par\medskip
{\bf 2.16 The complex group of type G$\hskip-1pt_2$.} {\it If $\cF_\Delta{\,\ne\,}\emptyset$, then 
$\Delta{\,\not\hskip1pt\cong\,}\Gtwo^{\CC}$\/}.
\par\smallskip
{\tt Proof.}  Suppose that 
$\Delta{\,\cong\,}\Gtwo^{\CC}{\,=\,}\Aut(\OO{\hskip1pt\otimes\hskip1pt}\CC)$ 
is the complexification  of its maximal compact subgroup $\Gamma{=\,}\Gtwo$. 
Then all involutions of $\Delta$ are conjugate by  \cite{cp} 11.31(d) and 93.10(a), 
and no involution can act trivially on the fixed line $W$ since $\Delta$ is strictly simple. 
The centralizer of a given involution 
$\lambda$ contains a subgroup $\Phi{\,\cong\,}\SO3\RR$,  see 2.5. If $\Phi$ would contain a 
reflection, then even three pairwise commuting conjugate ones whose centers form a triangle in some 
$\Phi$-orbit, and $\cF_\Delta{\,\subseteq\,}\cF_\Phi{\,=\,}\emptyset$.
Consequently, each involution in $\Delta$ is planar and the involutions in $\Phi$ induce 
planar involutions on  $\cF_\lambda$. We conclude that $W{\,\approx\,}\sS_8$ is a manifold, 
cf. \cite{cp} 55.6. 
According to 2.15, the action of $\Gamma$ on $W$ is equivalent to~the action of 
$\Aut\OO$ on $\OO{\smcup}\{\infty\}$, and the fixed points of $\Gamma$ on $W$ form 
a circle $C$.  Note that $\Gamma$ is even a maximal subgroup of $\Delta$\: (\cite{cp} 94.34).
If $c{\,\in\,}C$, then $\dim\Delta_c{\,\ge\,}20$ and $\Gamma{\,<\,}\Delta_c{\,=\,}\Delta$. 
The same argument, applied to the pencils $\frL_c$ with $c{\,\in\,}C$ instead of $W$, 
shows that $\cF_\Delta$ is a flat subplane, but this contradicts Stiffness.  \qed
\par\medskip
The next lemma will be needed repeatedly:
\par\smallskip
{\bf 2.17 Stabilizer of a triangle.} {\it Suppose that the connected subgroup $\nabla$ of $\Delta$ 
fixes the triangle $a,u,v$ and that $\nabla$ is transitive on $S{\,:=\,}uv\sm\{u,v\}$. For $w{\,\in\,}S$, 
let $\Psi$ and $\Phi$ be maximal compact subgroups of $\nabla$ and of 
$\Omega{\,=\,}\nabla_{\hskip-2pt w}$, respectively, such that $\Phi{\,\le\,}\Psi$. Then\/{\rm:}}  \par
\quad(a) {\it $\Phi$ and $\Psi$ have equally many almost simple factors\/}, \par
\quad(b) {\it $0{\,<\,}\Psi'{:}\Phi'{\,<\,}8$ and  $\Phi'{\,\ne\,}\1$\/},\: $\dim\Phi'{\,\ne\,}10$,   \par
\quad(c) {\it If $\dim\Phi'{\,=\,}3$, then $\Psi'{\,\cong\,}\U2\HH{\,\cong\,}\Spin5\RR$\/}, \par
\quad(\^c) {\it If $\dim\Phi'{\,=\,}6$, then $\Psi'{\,\cong\,}\U2\HH{\,\cdot}\U1\HH$\/}, \par
\quad(d) {\it If $\dim\Phi'{\,=\,}8$, then $\Psi'{\,\cong\,}\SU4\CC{\,\cong\,}\Spin6\RR$\/}, 
see also \cite{sz5} Lemma (5), \par
\quad(e) {\it If $\Phi'{\,\cong\,}\Gtwo$, then $\Psi'{\,\cong\,}\Spin7\RR$\/}, \par
\quad(f) {\it If $\Phi'$ is almost simple and $\dim\Phi'{\,>\,}14$, then $\Psi'{\,\cong\,}\Spin8\RR$\/}.
\par\smallskip
{\tt Proof.} Note that $S$ is homotopy equivalent ($\simeq$) to $\sS_7$.   By  the Mal'cev-Iwasawa theorem \cite{cp} 93.10,  $\nabla{\,\simeq\,}\Psi$ and
$\Omega{\,\simeq\,}\Phi$; therefore the exact homotopy
sequence (\cite{cp} 96.12) for the action of $\nabla$ on $S$ can be written in the form
$$\dots\to\pi_{q{+}1}S\to\pi_q\Phi\to\pi_q\Psi\to\pi_qS\to\pi_{q{-}1}\Phi\to\dots,\quad 
q{\,\ge\,}1\,.\leqno{(*)}$$ 
From \cite{cp} 94.31(c) it follows that $\pi_q\Psi{\,\cong\,}\pi_q\Psi'$ and $\pi_q\Phi{\,\cong\,}\pi_q\Phi'$ for $q{\,>\,}1$.
We have \\ $\pi_qS{\,=\,}0\; (q{\,<\,}7),\: \pi_7S{\,\cong\,}\ZZ$ (\cite{sp} 7.5.6 or \cite{br} II.16.4),\: 
$\pi_8S$ is finite (${\cong\,}\ZZ_2$),\: 
see \cite{sp} 9.7.7. For the same range of the dimension $q$ (and beyond), the homotopy groups 
$\pi_q$ of the al\-most  simple compact Lie groups $\Kappa$ are known, cf. the references leading 
up to Theorem  94.36 in \cite{cp}. This theorem states that $\pi_2\Kappa{\,=\,}0$ and 
$\pi_3\Kappa{\,\cong\,}\ZZ$  for all compact simple Lie groups. For $\dim\Kappa{\,\le\,}52$ in particular,   $\pi_4\U n\HH{\,\cong\,}\ZZ_2$,\: $\pi_4\Kappa{\,=\,}0$ in all other cases,  and 
$\pi_6\Kappa{\,\le\,}\ZZ_{12}$. Furthermore $\pi_5\SU n\CC{\,\cong\,}\ZZ$ for $n{\,>\,}2$,\: 
$\pi_5\Kappa{\,\cong\,}\ZZ_2$ for $\Kappa{\,\cong\,}\Spin n\RR\hskip6pt (n{=}3,5)$ or 
$\U n\HH\hskip6pt (n{\ge}1)$,\: $\pi_5\Kappa{\,=\,}0$   for  all other groups $\Kappa$.   Note also   
$\pi_7\SO3\RR{\,\cong\,}\ZZ_2$,\:  
$\pi_7\Kappa{\,\cong\,}\ZZ$ for  $\Kappa{\,\cong\,}\SU n\CC\hskip6pt  (n{>}3)$,\: 
$\Spin n\RR\hskip6pt (n{>\hskip1pt}4,\, n{\ne}8)$, or $\U n\HH\hskip6pt (n{>}1)$, 
$\pi_7\Spin8\RR{\,\cong\,}\ZZ^2$, and $\pi_7\Kappa{\,=\,}0$ for all other groups $\Kappa$. 
Recall from the last part of 2.10 that each non-trivial action of $\U2\HH{\,\cong\,}\Spin5\RR$ is 
faithful. \\
(a) is an immediate consequence of $\pi_3\Phi{\,\cong\,}\pi_3\Psi$. \\
(b) If $\Phi'{\,=\,}\Psi'$, then $(*)$ yields $\pi_7\Phi'{\,\cong\,}\pi_7\Psi'\to\ZZ\to\pi_6\Phi'$, but this contradicts finiteness of  $\pi_6\Phi'$. In the case  $\Phi'{\,=\,}\1$, step (a) implies  $\Psi'{\,=\,}\1$, which is 
impossible by the first part of~(b). Because   $\pi_1\Phi{\,\cong\,}\pi_1\Psi$, again  by  $(*)$, the torus factors of $\Phi$ and $\Psi$ have the same dimension (cf. the structure theorem \cite{cp} 31(c) for 
compact Lie groups). From $w^\nabla{\,\cong\,}S$ follows 
$\Psi'{:}\Phi'{\,=\,}\dim\Psi{-}\dim\Phi{\,<\,}\nabla{:\hskip1pt}\Omega{\,=\,}\dim w^\nabla{\,=\,}8$. 
If $\dim\Phi'{\,=\,}10$, then $\Phi'{\,\cong\,}\U2\HH$, and $(*)$ implies $\pi_4\Psi'{\,\cong\,}\ZZ_2$, but the groups  $\U n\HH$ are the only groups $\Kappa$ with $\pi_4\Kappa{\,\cong\,}\ZZ_2$ and 
$\dim\U3\HH{\,=\,}21$ is too large. \\
(c) Part of $(*)$ reads \ $0\to\pi_4\Phi'{\,=\,}\pi_4\U1\HH{\,=\,}\ZZ_2\to\pi_4\Psi'\to0$. 
The groups $\U n\HH$ are the only groups $\Kappa$ with $\pi_4\Kappa{\,\ne\,}0$.
The dimension bounds in (b) and 2.10 imply $\Psi'{\,\cong\,}\U2\HH$. \\
(\^c) By assumption $\Phi'{\,\circeq\,}(\SU2\CC)^2$. Hence $\Psi'$ is a product of two almost simple factors. From $(*)$ we obtain  $\pi_5\Psi'{\,\cong\,}\pi_5\Phi'{\,\cong\,}\ZZ_2^{\hskip3pt2}$.
It follows that each  factor of $\Psi'$ is a group $\U n\HH$, $n{\,\ge\,}1$. We have 
$6{\,<\,}\dim\Psi'{\,\le\,}13$, thus  $\Psi'{\,\cong\,}\U2\HH{\hskip2pt\cdot\hskip1pt}\U1\HH$, because a factor $\SO3\RR$ would contain planar involutions and $\U2\HH$ cannot fix a triangle in a Baer subplane.  \\
(d) The homotopy sequence $(*)$ yields $\pi_5\Psi'{\cong\,}\pi_5\SU3\CC{\,\cong\,}\ZZ$. 
The groups $\SU n\CC$ are the only groups $\Kappa$ with $\pi_5\Kappa{\,\cong\,}\ZZ$, and 
$\SU5\CC$ is too large. \\
(e) We have $\pi_5\Psi'\cong\pi_5\Gtwo{\,=\,}0$ by $(*)$ and hence $\Psi'{\,\cong\,}\Spin7\RR$ 
by (b), the fact that $\Gtwo{\,<\,}\Psi$,  and the last part of  2.10. \\
(f) follows from $\pi_5\Phi'{\,\cong\,}\pi_5\Psi'$,\: $\pi_5\SU4\CC{\,\cong\,}\ZZ$,\: 
$\pi_5\Spin7\RR{\,\cong\,}0$,\: $\pi_5\U3\HH{\,\cong\,}\ZZ_2$,\:
$\pi_5\SU5\CC{\,\cong\,}\ZZ$,\: $\pi_5\Spin8\RR{\,=\,}0$,\:  $\SU5\CC{\,:\,}\SU4\CC{\,>\,}7$, and 
$\dim\nabla{\,\le\,}30$. 
\par\medskip
{\bf 2.18 Corollary.} {\it If $\nabla$ is transitive on the complement 
$D{\,=\,}P\sm(au\smcup av\smcup uv)$ of a triangle, then $\nabla'{\,\cong\,}\Spin8\RR$\/}.
\par\smallskip 
{\tt Proof.} We use the same notation as in 2.17. If $\nabla$ is transitive on $D$, then 
$\nabla$ is a Lie group by \cite{cp} 53.2, and $\Omega{\,=\,}\nabla_{\hskip-1pt w}$ is transitive 
on $V{\,=\,}av\sm\{a,v\}$. Choose $c{\,\in\,}V$, and let $\Lambda$ be a maximal compact 
subgroup of $\Omega_c$. Stiffness implies that $\Lambda'$ is one of the groups 
$\Spin3\RR$, $\SO4\RR$, $\SU3\CC$, or $\Gtwo$. Applying 2.17 to the action of $\Omega$ on $V$, 
the possible groups $\Phi'$ can be determined. If $\dim\Lambda'{\,\ge\,}8$, then 
$\Phi'{\,\cong\,}\Spin n\RR$ with $n{\,\in\,}\{6,7\}$, and then 2.17(f) shows $\Psi'{\,\cong\,}\Spin8\RR$.
In the case $\dim\Lambda'{\,=\,}3$ we obtain $\dim\Phi'{\,=\,}10$, which is excluded by 2.17(b). 
If $\Lambda'{\,\cong\,}\SO4\RR$, then 2.17(\^c) shows
 $\Phi'{\,\cong\,}\U2\HH{\hskip2pt\cdot\hskip1pt}\U1\HH$. 
Now $13{\,<\,}\dim\Psi'{\,\le\,}20$, and $\Psi'$ is an almost direct product $\Xi{\hskip1pt\cdot}\Ypsilon$, where $\Xi{\,\cong\,}\Ypsilon{\,\cong\,}\U2\HH$ and $\Xi{\,\triangleleft\,}\Phi'$. According to 2.10, the 
torus rank of $\nabla$ is $4$; hence $\Psi'{\,=\,}\Psi$,\, $\Phi{\,=\,}\Phi'$, 
and $\dim\Phi{\,=\,}13$. By construction $\Phi{\,=\,}\Psi_w$, and 
$\dim w^\Ypsilon{\,=\,}\dim w^\Psi{\,=\,}\Psi{\hskip1pt:\hskip1pt}\Phi{\,=\,7}$. 
Therefore $\Xi|_{w^\Ypsilon}{\,=\,}\1$ and $\cF_{\Xi_c}{\,=\,}\langle a,c,w^\Ypsilon\rangle{\,=\,}\cP$,
but $\dim\Xi_c{\,>\,}2$, a contradiction. \qed 
\par\bigskip

{\Bf 3. No fixed elements}
\par\medskip
The Hughes planes form a one-parameter family of planes admitting a fixed element free group action. 
They  can be characterized as  follows: $\cP$ has a $\Sigma$-invariant Baer subplane  $\cH$ such that 
$\Sigma$ induces on $\cH$~the full automorphism  group $\PSL3\HH$, see \cite{cp} \S\hskip2pt86 or 
\cite{sz3} 3.21. If $\cP$ is a proper Hughes plane, then $\dim\Sigma{\,=\,}36$  and  $\Sigma$ is 
transitive on the set of flags of the {\it outer\/} subgeometry consisting of the points and lines not belonging to $\cH$, see \cite{cp} 86.5.  Sometimes it is convenient to consider~the~classical plane together with the  stabilizer of a Baer subplane also as a Hughes  plane. Other characterizations of the Hughes planes  are given in \cite{pw1} and \cite{sz5}.
\par\medskip
{\bf 3.0.} {\it If $\cF_\Delta{\,=\,}\emptyset$ and $\dim\Delta{\,\ge\,}23$, then $\Delta$ is a Lie 
group\/}. 
\par\smallskip
{\tt Proof.} (a) Suppose that $\Delta$ is not a Lie group. Then $\dim\Delta{\,<\,}27$ and there 
are arbitrarily small compact, $0$-dimensional  central subgroups $\Nu$ such that $\Delta/\Nu$ 
is  a Lie group, cf. \cite{cp} 93.8. Choose a fixed point~$x$  of some element 
$\zeta{\,\in\,}\Nu\sm\{\1\}$. 
From $\zeta|_{x^\Delta}{\,=\,}\1$ and  $\cF_\Delta{\,=\,}\emptyset$ it follows that 
$\cD{\,=\,}\langle p^\Delta\rangle$ is a proper connected subplane. 
Put $\Delta|_\cD{\,=\,}\Delta/\Kappa$. By Stiffness, $\dim\Kappa{\,\le\,}14$,\; 
$\cD$ is not flat, $\dim\cD{\,\ge\,}4$ 
and $\dim\Kappa{\,\le\,}8$. Hence $\Delta{\hskip1pt:\hskip1pt}\Kappa{\,\ge\,}15$, 
and $\Delta/\Kappa$ is a Lie group by  \cite{cp}  71.2 and  \cite{pw3}. 
We may assume, therefore, that $\Nu{\,\le\,}\Kappa$. If $\cD{\;\ldot}\cP$, then $\dim\Kappa{\,<\,}8$ 
by Stiffness,   and $\dim\Delta|_\cD{\,\ge\,}16$. The same is true, if $\dim\cD{\,=\,}4$: in this case 
$\cD$ is classical (\cite{cp} 72.8), $\Aut\cD$ has no subgroup of dimension $15$, and 
$\Delta/\Kappa{\,\cong\,}\PSL3\CC$. If $\dim\Kappa{\,=\,}8$, then 
$\Kappa^1{\,\cong\,}\SU3\CC$,\,  $\Delta/\Kappa^1$ is a finite covering 
of $\Delta/\Kappa$, and $\Delta$ would be a Lie group.  \\
(b) Suppose first that $\cD{\,=\,}\cF_\Nu{\;\ldot}\cP$. Then  $\Kappa$ is compact and
\cite{sz1} 2.1 applies to $\Delta|_\cD$: 
either  $\cD$ is the classical quaternion plane, $\Delta/\Kappa$ is the elliptic or hyperbolic motion 
group of $\cD$, and $2{\,\le\,}\dim\Kappa{\,\le\,}5$, or $\cD$ is a Hughes plane, 
$\Delta'/\Kappa{\,\cong\,}\SL3\CC$, and $\dim\Kappa{\,\ge \,}6$.  \\
(c)  In the first case, the motion group $\Delta/\Kappa$ is covered by a subgroup 
$\Upsilon{\hskip1pt\circeq\,}\U3(\HH,r)$ of $\Delta$, see \cite{cp}~94.27. The representation 
of $\Ypsilon$ on $\Kappa/\Nu$ shows that $\Kappa{\,\le\,}\Cs{}\Ypsilon$. We may assume that 
$\Kappa$ is connected, since $\dim\Kappa^1\Ypsilon{\,=\,}\dim\Delta$ and $\Delta$ is connected. Moreover, the center $\Zeta$ of $\Kappa$ has positive dimension: because  $\Kappa$ is not a Lie group, the claim follows from the structure of compact groups as described in \cite{cp} 93.11. The stabilizer $\Pi{\,=\,}\Delta_p$ of a non-absolute point $p{\,\in\,}\cD$ fixes also the polar $L$ of $p$, and $\Pi/\Kappa{\,\circeq\,}\U2\HH{\times}\U1\HH$. In particular, $\Pi$ is compact and 
$\dim\Pi{\,\ge\,}15$. Choose $s{\,\in\,}L\sm\cD$ and $z{\,\in\,}ps\sm\{p,s\}$, and note that 
$\Pi_s\smcap\Kappa{\,=\,}\1$ because $\langle \cD,s\rangle{\,=\,}\cP$. Hence $\Pi_s$ and 
$\Lambda{\,=\,}(\Pi_{s,z})^1$  are Lie groups.  From \cite{cp} 96.11 it follows that 
$\Pi{\hskip1pt:\hskip1pt}\Pi_s,\, \Pi_s{\hskip1pt:\hskip1pt}\Lambda{\,\le\,}7$, and 
$\dim\Lambda{\,>\,}0$. As $\Nu{\,\le\,}\Zeta{\,\le\,}\Cs{}\Delta{\,\le\,}\Cs{}\Lambda$, 
we conclude that 
$\langle s^\Zeta,p,z\rangle{\,\le\,}\cF_\Lambda$, and $\cF_\Lambda$ is   a connected subplane; 
in fact, $\cF_\Lambda{\;\ldot}\cP$\; (otherwise $\Zeta$ would be a Lie group by \cite{cp} 32.21 and 71.2). Stiffness implies $\Lambda{\,\le\,}\Spin3\RR$. The involution $\iota{\,\in\,}\Lambda$ induces 
an involution $\overline\iota{\,=\iota|_\cD\,}$ on $\cD$, and $\cF_\Lambda{\,=\,}\cF_\iota$. 
If $\overline\iota$ is planar, then $\cD\smcap\cF_\iota{\,\ldot}\cF_\iota$, and the lines of 
$\cF_\iota$ are $4$-spheres (see \cite{cp} 53.15 and \cite{sz3} 3.7). If $\overline\iota$ is a reflection, then the axis $A$ of   $\overline\iota$ is a $4$-sphere contained in a line of $\cF_\iota$. 
 From \cite{cp} 51.20 it follows that the lines of $\cF_\iota$ are manifolds, and then 
 $H{\,=\,}L\smcap\cF_\iota{\,\approx\,}\sS_4$.  \\
(d) The structure theorem \cite{cp} 93.11 for compact Lie groups together with \cite{cp} 93.19 shows  that $\Kappa$ is commutative or a product of a commutative central group $\Alpha$ with a factor 
$\Omega{\,\cong\,}\Spin3\RR$. In the latter case $\dim\Pi{\,\ge\,}17$ and
$\Lambda{\,\cong\,}\Spin3\RR$. As $\Aut\Omega{\,\cong\,}\SO3\RR$, the involution 
$\iota{\,\in\Lambda\,}$ centralizes $\Kappa$, and $H^\Kappa{\,=\,}H$. Therefore the compact connected group $\Kappa$ acts freely  on $s^\Kappa{\,\subset\,}H{\,\approx\,}\sS_4$, and Richardson's theorem (\cite{cp} 96.34 or $(\dagger)$) implies that $\Kappa$ is a Lie group contrary to the assumption. 
\\
(e) In  the second case mentioned in step (b) or if $\cD{\,<\,}\cF_\Nu$, the plane $\cF_\Nu$ is a 
Hughes plane, $\Delta$ has a subgroup $\Ypsilon{\cong\,}\SL3\CC$, and $\dim\Kappa{\,\in\,}\{6,7\}$. By \cite{cp} 86.33, there exists a $\Delta$-invariant $4$-dimensional classical complex  subplane 
$\cC$. If $\dim\cD{\,=\,}4$, we will change the notation and put $\cD{\,=\,}\cF_\Nu$. By Stiffness, 
$\Kappa$ is compact, and $\Kappa$ acts freely on the set of points outside~$\cD$. The structure of 
the compact Lie group $\Kappa/\Nu$ implies that $\Ypsilon{\,\le\,}\Cs{}\Kappa$. As proved in 
\cite{cp}~86.34, the center $\langle \zeta\rangle{\,\cong\,}\ZZ_3$ of $\Ypsilon$ acts effectively on 
$\cD$ and fixes each element of $\cC$. Consider a point $z$ which is not incident with a line of $\cC$. Because $\cD$ is a Baer subplane, there is a unique line $L$ of $\cD$ with $z{\,\in\,}L$, and 
$L^\zeta{\,\ne\,}L$. The stabilizer $\Ypsilon_{\!z}$ fixes each point of $z^\Kappa{\,\subset\,}L$ 
and of $z^{\zeta\Kappa}{\,\subset\,}L^\zeta$. Both orbits of $z$ are homeomorphic to $\Kappa$; 
their dimension is at least $6$. Hence $\langle z^\Kappa,z^{\zeta\Kappa}\rangle{\,=\,}\cP$ and 
$\Ypsilon_{\!z}{\,=\,}\1$. It follows that $z^{\!\Ypsilon}$ is open in $P$, and $\Delta$ would be a Lie group by \cite{cp} 53.2.
\par
(f) Only the following possibility remains:\,  $\cD{\,=\,}\cF_\Kappa{\,=\,}\cF_\Nu$ is the classical 
complex plane, $\Delta{\hskip1pt:\hskip1pt}\Kappa{\,=\,}16$,  and $\Delta$ has a subgroup 
$\Ypsilon{\,\circeq\,}\SL3\CC$ inducing the full automorphism group on $\cD$\: (see \cite{cp} 94.27).  
All involutions in $\Upsilon$ are conjugate, they are either planar or reflections of~$\cP$. There are 
lines which do not intersect the point set of $\cD$; dually some points are not incident with a line 
of $\cD$.\\
(g) If $\iota{\,\in\,}\Upsilon$ is a reflection of $\cP$, then for $L{\,\in\,}\cD$ the translation group with 
axis $L$ is isomorphic to $\RR^4$. Each translation $\tau$ is uniquely determined by its axis and a 
pair $p,p^\tau$ with $p^\tau{\,\ne\,}p{\,\in\,}\cD$. As $\Kappa|_\cD{\,=\,}\1$, it follows that 
$\Kappa{\,\le\,}\Cs{}\tau$. The almost simple group $\Ypsilon$ is generated by all these translations. 
Hence $\Ypsilon{\,\le\,}\Cs{}\Kappa$. Choose a point $z$ such that 
$\forall_{\hskip-1pt L{\hskip1pt\in\hskip1pt}\cD}\;z{\,\notin\,}L$. If $\dim\Ypsilon_{\!z}{\,=\,}0$, then 
$z^\Delta$ is open in $P$, and $\Kappa$ would be  Lie group by \cite{cp} 53.2. Therefore 
$\dim\Ypsilon_{\!z}{\,>\,}0$. Moreover, $\langle z^\Kappa\rangle$ is a subplane, or else 
$z^\Kappa$ is contained in some line $L{\,=\,}L^\Kappa{\,\in\,}\cD$. 
From $\Kappa{\,\le\,}\Cs{}\Ypsilon$ 
it follows that $\Ypsilon_{\hskip-2pt z}|_{\langle z^\Kappa\rangle}{\,=\,}\1$ and 
$\langle z^\Kappa\rangle{\,<\,}\cP$, in fact, $\langle z^\Kappa\rangle\ldot\cP$\: (if not, then
$\Kappa|_{\langle z^\Kappa\rangle}$ is a Lie group, and we may assume that 
$\Nu|_{\langle z^\Kappa\rangle}{\,=\,}\1$, but then $z^\Kappa{\,\subseteq\,}\cF_\Nu$ contrary the 
choice of $z$). By Stiffness, $\Ypsilon_{\hskip-2pt z}$ is a compact Lie group, and an involution 
$\iota{\,\in\,}\Ypsilon_{\!z}$ would be planar. \\
(h) Suppose next that $\iota$ is planar and that $\Kappa{\,\le\,}\Cs{}\Ypsilon$. Then 
$\langle \cD,\cF_\iota\rangle{\,=\,}\cP$ because $\iota|_\cD{\,\ne\,}\1$ and $\cD{\,\not\le\,}\cF_\iota$. 
Hence $\Kappa$ acts effectively on $\cF_\iota$. We have 
$\Gamma{\,=\,}\Cs\Ypsilon\iota{\,\cong\,}\GL2\CC$ and $\dim\Gamma\Kappa{\,=\,}15$. Stiffness 
implies $\dim\Gamma\Kappa|_{\cF_\iota}{\,=\,}h{\,\ge\,}12$. According to a theorem of Priwitzer 
\cite{pw3},\, $\Gamma\Kappa$ induces a Lie group on $\cF_\iota$, but 
$\Kappa{\,\cong\,}\Kappa|_{\cF_\iota}$ and $\Kappa$ is not a Lie group. Therefore $\Ypsilon$ has a 
non-trivial representation on the Lie algebra of $\Kappa/\Nu$, and \cite{cp} 95.10 shows that 
$\Ypsilon{\,\cong\,}\SL3\CC$. \\
(i) Again by Stiffness, the almost simple commutator group $\Gamma'{\,\cong\,}\SL2\CC$ induces 
on $\cF_\iota$ a group $\PSL2\CC{\,\cong\,}\SO3\CC$. Note that $\iota|_\cD$ is a reflectiom of 
$\cD$. The line $L$ of $\cP$ such that $L\smcap\cD$ is the axis of $\iota|_\cD$ intersects 
$\cF_\iota$ in a line $S{\,=\,}L\smcap\cF_\iota$, and $S^{\Gamma'}{\,=\,}S$. An involution of a compact subgroup $\Phi{\,\cong\,}\SO3\RR$ of $\Gamma'|_{\cF_\iota}$ cannot have the axis 
$S$\: (use \cite{cp} 55.35  and simplicity of $\Phi$). Hence such an involution fixes a Baer 
subplane of $\cF_\iota$. Consequently, $S{\,\approx\,}\sS_4$,  see \cite{cp} 53.15 and 
\cite{sz3} 3.7. The set $L\smcap\cD{\,\subseteq\,}S$  is $\Gamma'$-invariant. From Richardson's theorem \cite{cp}~96.34 it follows that the fixed points of $\Phi$ in $S\sm\cD$ form a circle $C$. 
The group $\Nu$ is in the center of $\Delta$, we have $\Nu{\,\le\,}\Cs{}\Phi$ and $C^\Nu{\,=\,}C$. 
By \cite{MZ} 6.1 Th.\,3, the induced group $\Nu|_C$ is a Lie group, and we may assume that 
$\Nu|_C{\,=\,}\1$. Now $\cD{\;\ldot}\cF_\Nu{\;\ldot}\cP$, and Stiffness implies that $\Kappa$ is compact, but then $\Ypsilon$ acts trivially on $\Kappa/\Nu$\: (use \cite{cp} 93.19 and 94.31(c)\,).  This contradicts step (h).   \qed
\par\medskip
{\bf 3.1 Semi-simple groups.} {\it If  $\cF_\Delta{\,=\,}\emptyset$ and $\Delta$ is semi-simple, 
and if  $\dim\Delta{\,>\,}24$, then $\cP$  is a Hughes plane\/}.
\par\smallskip
{\tt Proof.} (a) {\it The assertion is true, if the center $\Zeta{\,=\,}\Cs{}\Delta$ contains an element 
$\zeta{\,\ne\,}\1$\/}. \\  There is some point $x$ such that $x^\zeta{\,=\,}x\,$  (see 2.1) and 
$\zeta|_{x^\Delta}{\,=\,}\1$. As $\,\cF_\Delta{\,=\,}\emptyset$, the orbit $x^\Delta$ is not contained 
in a line, and $\cD{\,=\,}\langle x^\Delta\rangle$ is a connected proper subplane. Put 
$\Delta^{\hskip-2pt\ast}{\,=\,}\Delta|_\cD{\,=\,}\Delta/\Kappa$ and  $\Lambda{\,=\,}\Kappa^1$. If $\cD$ 
is flat,  then  $\Delta{\hskip1pt:\hskip1pt}\Lambda{\,\le\,}8$,  $\dim\Lambda{\,\le\,}14$, and 
$\dim\Delta{\,\le\,}22$. Suppose that  $\dim\cD{\,=\,}4$. Then \cite{cp} 71.8 shows that 
$\dim\Delta^{\hskip-2pt\ast}{\,\le\,}8$ or $\Delta^{\hskip-2pt\ast}{\,\cong\,}\PSL3\CC$. In the first case, Stiffness 
would imply $\dim\Delta{\,<\,}22$. Therefore  $\Delta{\hskip1pt:\hskip1pt}\Lambda{\,=\,}16$, and 
$\dim\Lambda{\,>\,}8$. For any point $z$ outside $\cD$ on a line of $\cD$ it follows that 
$\dim z^\Lambda{\,\le\,}8$ and $\Lambda_z{\,\ne\,}\1$. Consequently, $\langle \cD,z\rangle$ 
is a Baer  subplane, and 2.6(c,e) imply that $\Lambda$ is compact and that $\dim\Lambda{\,\le\,}8$, a contradiction.  
Hence $\cD{\;\ldot}\cP$, $\,\Lambda$ is a compact connected normal, hence semi-simple group, 
the structure theorem \cite{cp} 93.11 shows that $\Lambda$ is a Lie group,  Stiffness implies
$\Lambda{\,\le\,}\Spin3\RR$, and  $\dim\Delta^{\hskip-2pt\ast}{\,>\,}21$.   Now 
\cite{sz1} 2.1 shows that $\dim\Delta^{\hskip-2pt\ast}{\,=\,}35$ and $\Delta^{\hskip-2pt\ast}{\,\cong\,}\PSL3\HH$ as claimed. \\
(b) Thus we may assume that $\Delta$ has  trivial center. If a proper simple factor $\Gamma$ of 
$\Delta$ contains a reflection with a center $c$, then the complement $\Lambda$ of $\Gamma$ 
in $\Delta$ acts trivially on the subplane $\cD{\,=\,}\langle c^\Gamma\rangle{\,=\,}\cD^\Delta$. 
As before, $\cD$  is a Baer subplane, and then $\Lambda{\,\cong\,}\Spin3\RR$, and $\Delta$ would have a center  $\Zeta{\,\ne\,}\1$. \\
(c) Assume that $\Zeta{\,=\,}\1$ and that $\beta$ is a planar involution in a proper simple factor 
$\Gamma$ of $\Delta$ of minimal dimension. Let $\Delta{\,=\,}\Gamma{\times}\Psi$ and consider the action of $\Psi$ on the Baer plane $\cB{\,=\,}\cF_\beta$. If $\Psi$ fixes some point $x$ in $\cB$, then
$x^\Gamma{\,\ne\,}x$, $\,\Psi|_{x^\Gamma}{\,=\,}\1$, $\,x^\Gamma$ is not contained in a line, and
$\cE{\,=\,}\langle x^\Gamma\rangle{\,=\,}\cE^\Delta{\,<\,}\cP$. Now $\dim\Psi{\,\le\,}14$, 
$\,\dim\Gamma{\,>\,}8$, $\,\dim\cE{\,>\,}2$, and $\dim\Gamma{\,\le\,}\dim\Psi{\,\le\,}8$ by Stiffness,
a contradiction. Hence $\Psi$ has no fixed elements in $\cB$ and  $\Psi{\,\cong\,}\Psi|_\cB$, since 
$\Zeta{\,=\,}\1$ and a non-trivial kernel would be isomorphic to $\Spin3\RR$; moreover, 
$\Psi{\,\not\cong\,}\SL3\CC$, and  
\cite{sz1} 2.1 implies that $\cB{\,=\,}(B,\frB)$ is the classical quaternion plane and that $\Psi$ contains one of the 
$3$  motion groups. Each orbit of $\Psi$ on $B$ is either open or a sphere of dimension $7$ or $5$, see \cite{cp} 13.17 and 18.32 or  
\cite{sz1} 1.15. As $\beta$ is not central, there is some conjugate 
$\beta^\gamma{\,\ne\,}\beta$ in $\Gamma$. By \cite{cp} 55.38 or \cite{lw2}, the Baer subplanes 
$\cB$ and $\cB^\gamma$ have a common point $c$, and 
$c^\Psi{\,\subseteq\,}B{\smcap}B^\gamma$. Consequently, $B{\,=\,}B^\gamma$ for each 
$\gamma{\,\in\,}\Gamma$. The simple group $\Gamma$ 
is generated by its conjugacy class $\beta^{\hskip1pt\Gamma}$, and $\Gamma|_\cB{\,=\,}\1$. 
From the stiffness result 2.6(\^b) it follows that $\Gamma{\,\cong\,}\Spin3\RR$ has a non-trivial center contrary to the assumption. Hence $\Delta$ is strictly simple. \\
(d) Even without  assuming  $\cF_\Delta{\,=\,}\emptyset$, the theorem is true for semi-simple groups $\Delta$ of dimension ${>}28\,$ (see \cite{pw1}). Therefore it suffices to exclude the 
groups $\POpr8(\RR,r)$, $0{\,\le\,}r{\,\le\,}4$, and the complex group  
$\Gtwo^{\CC}{\,=\,}\Aut(\OO{\hskip1pt\otimes\hskip1pt}\CC)$. For $r{\,<\,}4$, the orthogonal groups contain a  subgroup $\SO5\RR$ and cannot act on $\cP$ by 2.10. \\
(e) A maximal compact subgroup $\Omega$ of $\POpr8(\RR,4)$ is isomorphic to 
$(\SO4\RR)^2/\langle -\1\rangle$. Each factor $\SO4\RR$ has   a subgroup 
$\SO3\RR{\,>\,}\ZZ_2^{\hskip2pt2}$, such that the elements of the latter group are represented by 
diagonal matrices with entries $\pm1$. By 2.10, the product 
$\Phi{\,=\,}\ZZ_2^{\hskip2pt2}{\times\ZZ_2^{\hskip2pt2}}$ contains a reflection $\alpha$, and each conjugate of $\alpha$ in $\Phi$ is also a reflection. Consequently, there are exactly $3$ reflections in $\Phi$, and 2.10 implies that the central involution $\omega{\,\in\,}\Omega$ is planar. 
Stiffness 2.6(\^b) shows that 
$\Omega/\omega{\,\cong\,}(\SO3\RR)^4{\,>\,}\ZZ_2^{\hskip2pt8}$ acts effectively on $\cF_\omega$, 
but this contradicts 2.10. \\
(f) The real form $\Gamma{=\,}\Gtwo$ is a maximal compact subgroup of 
$\GGamma{\,=\,}\Gtwo^{\CC}$.    By 2.5, 
the centralizer of any involution contains a subgroup $\SO4\RR$. Hence $\Gamma$ has a subgroup 
$\ZZ_2^{\hskip2pt3}$, and then 2.\,10 and 5 imply that each involution $\beta{\,\in\,}\GGamma$ is planar, and that $\Cs\Gamma\beta{\,\cong\,}\SO4\RR$. Consequently, the centralizer of $\beta$ in $\GGamma$ induces on $\cF_\beta$ a group $\Chi{\,\cong\,}(\SO3\CC)^2$, and this action is effective by Stiffness. 
Inspection of \cite{sz1} 7.3 shows that $\Chi$ fixes a non-incident pair $\{o,W\}$ in $\cF_\beta$. 
It follows that each involution in $\Chi$ is planar. On the other hand, $\Chi$ has a subgroup 
$\ZZ_2^{\hskip2pt4}$ and must contain a reflection by 2.10. \qed
\par\smallskip
{\tt Remark.} {\it If $\Delta$ is semi-simple, $\cF_\Delta{\,=\,}\emptyset$ and $\dim\Delta{\,\ge\,}21$, then $\Delta$ is a Lie group\/}.
\par\smallskip
{\tt Proof.} It suffices to show that a maximal compact subgroup $\Phi$ of $\Delta$ is semi-simple,
cf. \cite{cp} 93.\hskip1pt10,11, or that the center $\Zeta$ of $\Delta$ is a Lie group. Suppose that 
$\1{\,\ne\,}\zeta{\,\in\,}\Zeta$, and let 
$x$ be a fixed point of $\zeta$. Then $x^\Delta$ is not contained in a line, and
$\cD{\,=\,}\langle x^\Delta\rangle{\,\le\,}\cF_\zeta{\,<\,}\cP$. Put 
$\Delta^{\hskip-1.5pt\ast}{\,=\,}\Delta|_\cD{\,=\,}\Delta/\Kappa$ and $\Kappa^1{\,=\,}\Lambda$;
both $\Delta^{\hskip-1.5pt\ast}$ and $\Lambda$ are semi-simple. 
Note that $\cD$ is connected, and write $\dim\cD{\,=\,}d$. If $d{\,=\,}2$, then 
$\Delta^{\hskip-1.5pt\ast}{\le\,}\SL3\RR$ and $\Lambda{\,\cong\,}\Gtwo$ or 
$\dim\Kappa{\,\le\,}10$\: (see 2.6(d)\,). Hence $\dim\Delta^{\hskip-1.5pt\ast}{\,=\,}8$ and  
$\Phi$ is semi-simple. If $d{\,=\,}4$, then $5{\,\le\,}\dim\Lambda{\,\le\,}8$,\: 
$\dim\Delta^{\hskip-1.5pt\ast}{\,\ge\,}13$, and then  $\Delta^{\hskip-1.5pt\ast}{\,\cong\,}\PSL3\CC$,\: (see  \cite{sz15} 5.6, \cite{cp} 71.8, or 2.14 above). By Stiffness, $\Lambda{\,\cong\,}\SU3\CC$ or 
$\dim\Lambda{\,=\,}6$. If $\cD{\,\ldot}\cB{\,<\,}\cP$ for some (closed) subplane $\cB$, then $\Lambda$ is compact by Stiffness, and $\Phi$ is semi-simple. We may assume, therefore, that 
$\cD$ is a maximal  subplane of $\cP$, in particular, $\cD{\,=\,}\cF_\zeta$.   Choose a line 
$L{\,\in\,}\cD$ and a point $p{\,\in\,}L\sm\cD$. Then $\Kappa_p{\,=\,}\1$, so that   $\Delta_p$ 
acts faithfully on $\cD$. Moreover, $\dim\Delta{\,\ge\,}16{+}6$ and $\dim\Delta_p{\,\ge\,}10$. 
Hence $\cD$ is the cšassical complex plane, see \cite{cp} 72.8. For two distinct points 
$a,b{\,\in\,}\cD\sm L$, the stabilizer $\Gamma{=\,}\Delta_{p,a,b}$   fixes each point  of the orbit $p^{\langle \zeta\rangle}$ and of the line $ab\smcap\cD$, and $\Gamma$ is  contained in  a 
group $\CC{\rtimes}\CC^{\times}$. As $\dim\Gamma{\,\ge\,}2$, the group $\Gamma$ is 
not compact. It follows that $\cF{\,=\,}\cF_\Gamma$ is a $4$-dimensional subplane 
${\ne\,}\cD$. Note that both $\cD$ and $\cF$ are $\Zeta$-invariant. The induced groups
$\Zeta|_\cD{\,=\,}\Zeta/\Nu$ and $\Zeta|_\cF{\,=\,}\Zeta/\Xi$ are Lie groups by \cite{cp} 71.2, and
$\Nu\smcap\Xi$ acts trivially on $\langle \cD,\cF\rangle$, which coincides with $\cP$ by 
maximality of $\cD$. Consequently, 
$\Nu\smcap\Xi{\,=\,}\1$ and $\Zeta$ itself is a Lie group, cf. \cite{cp} 94.3.
Finally let $\cD{\;\ldot}\cP$. By Stiffness, $\Lambda$ is 
compact and semi-simple, hence $\Lambda{\,\cong\,}\Spin3\RR$ or $\Lambda{\,=\,}\1$. 
Thus $\dim\Delta^{\hskip-1.5pt\ast}{\,\ge\,}18$, and we  may assume that $\dim\Delta{\,<\,}27$.  
Then \cite{sz1} 2.1 implies $\Delta^{\hskip-1.5pt\ast}{\,\cong\,}\PU3(\HH,r)$, and $\Phi$ is  semi-simple. \qed 

\par\medskip
{\bf 3.2 Normal torus.} {\it Assume that  $\cF_\Delta{\,=\,}\emptyset$ and that 
$\dim\Delta{\,\ge\,}23$, so that $\Delta$ is a Lie group. If $\Delta$ has a normal torus 
subgroup $\Theta$, then  $\cP$  is a Hughes plane\/}.
\par\smallskip
{\tt Proof.} The torus $\Theta$ is even central (\cite{cp} 93.19), hence it cannot contain a reflection, and there is a planar involution $\iota$ in the center of $\Delta$. Let 
$\Delta^{\hskip-2pt\ast}{\,=\,}\Delta|_{\cF_\iota}{\,=\,}\Delta/\Kappa$. 
By Stiffness, $\dim\Kappa{\,\le\,}3$. If $\Theta|_{\cF_\iota}{\,\ne\,}\1$, then $\cF_\iota$ has a $\Delta$-invariant Baer subplane and $\dim\Delta^{\hskip-2pt\ast}{\,\le\,}17$\: (see \cite{cp} 86.35 or use \cite{cp} 83.11). Hence 
$\Theta{\,\trianglelefteq\,}\Kappa^1$ and $\Kappa^1{\,\not\cong\,}\Spin3\RR$.  Therefore 
$\Kappa^1{\,=\,}\Theta$ and $\dim\Delta^{\hskip-2pt\ast}{\,>\,}21$; moreover, $\Delta$ fixes no 
element in  $\cF_\iota$. Now  
\cite{sz1} 2.1 implies that $\dim\Delta^{\hskip-2pt\ast}{\,=\,}35$ and $\Delta^{\hskip-2pt\ast}{\,\cong\,}\PSL3\HH$  as claimed. \qed
\par\medskip
{\bf 3.3 Normal vector group.} {\it If $\cF_\Delta{\,=\,}\emptyset$ and if $\Delta$ has a 
{\rm(}minimal\,{\rm)} normal vector subgroup $\Theta$, then $\dim\Delta{\,\le\,}23$\/}.
\par\smallskip
{\tt Proof.} By Lemma 2.13, the group $\Theta$ fixes some element, say a point $a$, and the 
assumption $\cF_\Delta{\,=\,}\emptyset$ implies that $\langle a^\Delta\rangle{\,=\,}\cD$ is a 
connected subplane of $\cP$; moreover, $\Theta|_\cD{\,=\,}\1$ and $\cD{\,<\,}\cP$. Put 
$\Delta^{\hskip-2pt\ast}{\,=\,}\Delta|_\cD{\,=\,}\Delta/\Kappa$ and note that 
$\Theta{\,\le\,}\Lambda{\hskip1pt:=\,}\Kappa^1$; in particular, $\Lambda$ is not compact, 
$\dim\cD\hskip1pt|\hskip1pt4$, and $\dim\Lambda{\,\le\,}10$ by 2.6(d). 
If $\cD$ is flat, then $\dim\Delta^{\hskip-2pt\ast}{\,\le\,}8$ and $\dim\Delta{\,\le\,}18$.
Hence we may assume that $\dim\cF_\Lambda{\,=\,}4$. Stiffness 2.6(e) yields 
$\dim\Lambda{\,\le\,}7$. As $\Delta^{\hskip-2pt\ast}$ fixes no element in $\cD$, it follows from
\cite{cp} 71.\,4 and 8 or 2.14  that $\Delta^{\hskip-2pt\ast}$ is semi-simple, even almost simple, and that $\dim\Delta^{\hskip-2pt\ast}{\,\le\,}8$ or $\Delta^{\hskip-2pt\ast}{\,\cong\,}\PSL3\CC$. \qed
\par\smallskip
{\tt Remarks.} Cf. also \cite{gs} XI.10.19.  Under the assumptions of 3.3, the invariant subplane 
$\cD$  is uniquely determined.  $\langle \cD,z\rangle{\,=\,}\cP$ for each point $z{\,\notin\,}\cD$ \ 
(or $\Lambda$ would be compact  by 2.6(c)\,). Hence $\Lambda$ acts freely on the set of points outside $\cD$. There is a covering group  $\Psi$ of $\Delta^{\hskip-2pt\ast}$ in $\Delta$, see \cite{cp} 94.27. As in \cite{sz2} Th.\,1(e) it follows that the involutions in $\Psi$ are planar or 
$\Lambda{\,\le\,}\Cs{}\Psi$. Problem: is there a plane with a group 
$\Delta{\,\cong\,}\SL3\CC{\ltimes}\CC^3$\, and $\cF_\Delta{\,=\,}\emptyset$\,?
\par\bigskip
{\Bf 4. Exactly one fixed element}
\par\medskip
Throughout this section, assume that $\Delta$ fixes a unique line $W$ and no point. In the classical octonion plane, $\Sigma_W$ has a subgroup $\Delta{\,\cong\,}\Spin{10}(\RR,1){\ltimes}\RR^{16}$ of 
dimension $61$,\, cf. \cite{cp} 15.6; in all other cases $\dim\Delta{\,<\,}40$ by \cite{cp} 87.7. 
In \cite{sz8} it has been shown that $\cP$ is a translation plane whenever $\dim\Delta{\,\ge\,}35$.
\par\medskip
{\bf 4.0.} {\it If $\cF_\Delta{\,=\,}\{W\}$ and if $\dim\Delta{\,\ge\,}23$, then $\Delta$ is a Lie group\/}.
\par\smallskip
{\tt Proof.} (a) There exist arbitrarily small compact central  subgroups $\Nu{\,\le\,}\Delta$ of dimension $0$ such that $\Delta/\Nu$ is a Lie group, cf. \cite{cp} 93.8. 
If $\Nu$ acts freely on  $P{\ssm}W$, then each stabilizer $\Delta_x$ with 
$x{\,\notin\,}W$  is a Lie group because $\Delta_x{\smcap}\Nu{\,=\,}\1$. By \cite{cp} 
51.\hskip1pt6 and 8 and 52.12, the one-point compactification $X$ of $P{\ssm}W$ is 
homeomorphic to the quotient space $P/W$, and $X$ is a Peano continuum 
(i.e., a continuous image of the unit interval); moreover, $X$ is homotopy equivalent to 
$\sS_{16}$, and  $X$ has Euler characteristic $\chi(X){\,=\,}2$.  According to a theorem of 
L\"owen \cite{lw3}, these properties suffice for $\Delta$ to be a Lie group. \\
(b) Suppose now that  $x^\zeta{\,=\,}x$ for some $x{\,\notin\,}W$ and some 
$\zeta{\,\in\,}\Nu{\ssm}\1$.  
By assumption, $x^\Delta$ is not contained in a line and hence generates a $\Delta$-invariant 
subplane $\cD{\,=\,}\langle x^\Delta\rangle{\,\le\,}\cF_\zeta{\,<\,}\cP$. Put 
$\Delta^{\hskip-1pt*}{\,=\,}\Delta|_\cD{\,=\,}\Delta/\Kappa$. If $\cD$ is flat, then 
$\dim\Delta{\,\le\,}6{+}14$ by Stiffness; similarly,  $\dim\cD{\,=\,}4$ implies 
$\dim\Delta{\,\le\,}12{+}8$. Hence $\cD{\;\ldot}\cP$,\: $\Kappa$ is compact, 
$\dim\Kappa{\,<\,}8$,\: $\dim\Delta^{\hskip-1pt*}{\,\ge\,}16$,\: $\Delta^{\hskip-1pt*}$ is a 
Lie group\, (\cite{pw3}),
and we may assume that $\1{\,\ne\,}\Nu{\,\le\,}\Kappa$. If even $\dim\Delta^{\hskip-1pt*}{\,>\,}16$, 
then \cite{sz1} 1.10 shows that $\cD$ is the classical quaternion plane, because $\Delta^{\hskip-1pt*}$ does not fix a flag. \\
(c) If $\Delta$ is not transitive on $S{\,=\,}W\smcap\cD$, then there are
points $u,v,w{\,\in\,}S$ such that $\Delta{\hskip1pt:\hskip1pt}\Delta_{u,v,w}{\,\le\,}9$. Choose a line
$L$ of $\cD$ in the pencil $\frL_v\sm\{W\}$ and a point $z{\,\in\,}L\sm\cD$, so that 
$\Delta_z\smcap\Kappa{\,=\,}\1$ and $\Delta_z$ acts faithfully on $\cD$. In particular, 
$\Delta_z$ fixes $L$ and $\Delta_z$ is a Lie group. Moreover, $\dim z^{\Delta_L}{\,<\,}8$\:
(or else $\Kappa$ would be a Lie group by \cite{cp} 53.2). Therefore
$\Delta{\hskip1pt:\hskip1pt}\Delta_z{\,\le\,}3{+}4{+}7{\,=\,}14$ and 
$\Lambda{\,=\,}(\Delta_{z,u,w})^1$ satisfies $\dim\Lambda{\,\ge\,}3$. 
Equality is possible only if $\Delta$ is triply transitive on $V{\,=\,}v^\Delta$.
Next, it will be shown that $\cF_\Lambda$ is connected. If $\dim\Kappa{\,<\,}7$, then $\cD$ is 
classical  and $\Lambda$ fixes all points of a circle in $S$ containing $u,v,w$.
By \cite{cp} 55.32, the compact group $\Kappa$ has no subgroup $\TT^2$. Hence $\Kappa^1$ is commutative or an almost direct product of $\Omega{\,\circeq\,}\SO3\RR$ with a commutative group.
From \cite{cp} 93.19 it follows that $\Lambda$ acts trivially on the commutative factor $\Alpha$ 
of $\Kappa^1$. If $\dim\Kappa{\,>\,}3$, then $z^\Alpha{\,\subseteq\,}\cF_\Lambda$ has positive
dimension. In any case, $\dim\cF_\Lambda{\,\ge\,}2$. As $\Nu$ acts freely on $L\sm\cD$ and 
 groups of planes of dimension ${\le\,}4$ are Lie groups (cf.  \cite{cp} 32.21  and 71.2), it
 follows that $\cF_\Lambda{\,\ldot}\cP$. Stiffness implies $\Lambda{\,\cong\,}\Spin3\RR$. 
Consequently, $\Delta_z$ is doubly transitive on $V\sm\{v\}{\,\approx\,}\RR^3$, the group 
$\nabla{\,=\,}(\Delta_{z,u})^1$ is $6$-dimensional, and
$\nabla|_V{\,\cong\,}e^\RR{\hskip1pt\cdot\hskip1pt}\SO3\RR$. 
Now $\dim\Lambda|_V{\,=\,}1$, but the almost simple group $\Lambda$ does not have a $1$-dimensional factor groi§p. \\
(d) The previous step shows that $\Delta$ is transitive on $S{\,=\,}W\smcap\cD$. Hence 
$S{\,\approx\,}\sS_4$, and $\Delta$ has a subgroup {\cyss Yu}${\,\cong\,}\U2\HH{\,\cong\,}\Spin5\RR$, see  \cite{cp} 53.2,\;96.19--22,\;55.40, and 94.27. The central involution $\sigma$ of {\cyss Yu} induces on $\cD$ a reflection with axis $S$. As its center is not fixed by $\Delta$, it follows from \cite{cp} 61.13 that $\Delta^{\hskip-1.5pt*}$ contains a transitive translation group with axis~$S$. Therefore $\Delta{\hskip1pt:\hskip1pt}\Kappa{\,\ge\,}18$ and $\cD$ is classical. By the last part of 2.10, 
the involution $\sigma$ is even a reflection of $\cP$, and 
$\sigma^\Delta\sigma{\,=\,}\Tau{\,\cong\,}\RR^8$ is normal in $\Delta$. 
The group  $\Delta|_\cD$ is contained in 
$\SL2\HH{\hskip1pt\cdot\hskip1pt}\HH^{\times}{\ltimes}\RR^8$. 
Either  $\Delta|_S$ is compact, or $\Delta$ has a subgroup $\Ypsilon{\,\cong\,}\SL2\HH$.
In the first case, $\Delta|_S{\,=\,}${\cyss Yu}$|_S{\,\cong\,}\SO5\RR$. For $v,w{\,\in\,}S$, 
the stabilizer $\Delta_v$ fixes a second point $u{\,\in\,}S$, and 
$\Delta{\hskip1pt:\hskip1pt}\Delta_{u,v,w}{\,=\,}4{+}3$. Choose  a line $L$ of $\cD$ in the pencil 
$\frL_v\sm\{W\}$, a point $z{\,\in\,}L\sm\cD$, and put $\Gamma{\,=\,}\Delta_{L,w}$. 
Then $\Delta{\hskip1pt:\hskip1pt}\Gamma{\,\le\,}11$,
and $\Lambda{\,=\,}(\Gamma_{\hskip-2pt z})^1$ satisfies $\dim\Lambda{\,\ge\,}4$.   
Again  $\Lambda$ is a Lie group since $\Nu$ acts freely on $L\sm\cD$. 
As $\cD$ is the classical quaternion plane, any collineation which fixes $3$ collinear points 
of $\cD$ even fixes all points of a circle. Hence $\cF_\Lambda$ is connected. 
By \cite{cp} 32.21  and 71.2, we  have $z^\Nu{\,\subseteq\,}\cF_\Lambda{\;\ldot}\cP$. 
Now Stiffness would imply $\dim\Lambda{\,\le\,}3$ contrary to what has been stated before. \\
(e) In the case  $\SL2\HH{\,\cong\,}\Ypsilon{\,\le\,}\Delta$, 
choose $z{\hskip1pt\in\hskip1pt}W\sm\cD$ and $p{\,\notin\,}\cD\smcup W$\hskip-2pt. As $\Kappa$ acts freely outside~$\cD$,\; \cite{cp} 53.2 implies $\dim z^\Delta{\,<\,}8$ and 
$\dim p^\Delta{\,<\,}16$. Consequently, $\Lambda{\,=\,}(\Delta_{p,z})^1$ has positive dimension, 
and $\Lambda$ is a Lie group because $\Lambda\smcap\Kappa{\,=\,}\1$. 
Moreover, $\langle p^\Nu,z^\Nu\rangle{\,\le\,}\cF_\Lambda{\,<\,}\cP$. The orbit $p^\Nu$ is contained in a unique  ``inner'' line $L$ of $\cD$ intersecting $W$ in a point $v$. Hence the line 
$pz$ does not belong to $\cD$; it intersects $\cD$ in a unique inner point $a$, and $av{\,\ne\,}L,W$. By definition, $\Lambda$ fixes $L,v,$ and $a$ and hence $3$ distinct  lines of the 
quaternion plane $\cD$ in the  pencil $\frL_v$. Therefore $\cF_\Lambda$ is connected. 
Again $\dim\cF_\Lambda{\,>\,}4$,\;  $\cF_\Lambda{\,\ldot}\cP$, and $\Lambda$ is isomorphic to a subgroup of $\Spin3\RR$. In particular, $\Lambda$ is compact and there exists an  involution 
$\lambda{\,\in\,}\Lambda$ with $\cF_\lambda{\,=\,}\cF_\Lambda$. If $\dim\Delta{\,=\,}23$, then a maximal compact subgroup is isomorphic to {\cyss Yu} and $\Delta$ is a Lie group. Therefore 
$\dim\Delta{\,>\,}23$,\; $\dim\Lambda{\,>\,}1$, and $\Lambda{\,\cong\,}\Spin3\RR$. Obviously
$\overline\lambda{\,=\,}\lambda|_\cD{\,\ne\,}\1$,\; $\overline\lambda$ is not conjugate to the 
central involution $\sigma$ of {\cyss Yu},  and Richardson's theorem $(\dagger)$ 
shows that $\Lambda$ does not fix a $2$-sphere in $S{\,=\,}W\smcap\cD$. Therefore 
$\overline\lambda$ is a reflection with center $v$ and some axis $au\smcap\cD$. 
Up to conjugation, $a^{\hskip-1.5pt\Upsilon}{=\,}a$. As $\Upsilon_{\hskip-2pt v}$ is (doubly) transitive on $S\sm\{v\}$, the dual of \cite{cp} 61.19 shows that the elation group 
$\Delta^{\hskip-1.5pt*}_{[v,av]}$ is transitive. Consequently 
$\Delta{\hskip1pt:\hskip1pt}\Kappa{\,=\,}27$, and then $\dim\Lambda{\,\ge\,}5$, 
a contradiction. \qed
\par\medskip

{\bf 4.1 Theorem.} {\it  If $\Delta$ is transitive on $W$, then $\cP$ is classical\/}.
\par\smallskip
{\tt Proof.}  By \cite{cp} 42.8 and 44.3, the space $W$ is locally contractible and $\Delta|_W$ is a locally
compact transformation group with a countable basis. Hence  
\cite{HK} Cor.\,5.5 applies and $W$ is a topological manifold, in fact, $W{\,\approx\,}\sS_8$, see \cite{cp} 52.3. 
From \cite{cp}  96.\,19,\,21, and~22 it follows that $\Delta|_W$ has a subgroup $\SO9\RR$;  
 it is covered by an almost simple subgroup $\Phi{\,\le\,}\Delta$\, (cf. \cite{cp} 94.27). The last part 
 of 2.10 implies that $\Phi{\,\cong\,}\Spin9\RR$. The central involution $\sigma{\,\in\,}\Phi$ is a 
 reflection with axis $W$ and some center $c$. By assumption $c^\Delta{\,\ne\,}c$, and 2.9 shows 
that the translation group $\sigma^\Delta\sigma{\,=\,}\Tau$ has dimension $\dim c^\Delta{\,>\,}0$. 
As $\Delta$ is transitive on $W$, 
each group $\Delta_{[z,W]}$ with $z{\,\in\,}W$ has the same dimension, and then
$\Tau{\,\cong\,}\RR^{16}$ by \cite{sz11} or \cite{cp} 61.13, and $\dim\Delta{\,\ge\,}52$. \qed \\
For a related result see \cite{cp} 81.19.
\par\medskip
Without the restriction $|\cF_\Delta|{\,=\,}1$, Priwitzer \cite{pw2}, \cite{pw1} has treated
semi-simple groups of dimension $\ge29$. She proved in particular:
\par\smallskip
{\bf 4.2 Theorem.} {\it If $\cP$ is not classical,  if $\Delta$ is semi-simple, $\cF_\Delta{\,\ne\,}\emptyset$, and  $\dim\Delta{\,>\,}28$, then $\Delta{\,\cong\,}\Spin9(\RR,r)$ with $r{\,\le\,}1$, 
the central involution of $\Delta$ is a reflection, and $\cF_\Delta$ is a non-incident point-line pair\/}.
\par\medskip
{\bf 4.3 Semi-simple groups.} {\it If  $\cF_\Delta{\,=\,}W$ and if $\Delta$ is semi-simple, then  
$\dim\Delta{\,\le\,}16$ or $\dim\Delta{\,=\,}21$ and $\Delta$ is almost simple 
or $\Delta$ may be a  product  of at most $7$ factors isomorphic to the simply connected covering 
group of $\SL2\RR$. If $\dim\Delta{\,>\,}13$, then $\Delta$ is a Lie group\/}, 
see \cite{szp} for a {\tt proof}.
\par\medskip
{\bf 4.4 Normal torus.} {\it Suppose that $\Delta$ is a Lie group, that $\dim\Delta{\,\ge\,}18$, and 
that $\Delta$ has a normal torus subgroup $\Theta$. Then the fixed elements of the involution 
$\iota{\,\in\,}\Theta$ form a $\Delta$-invariant classical quaternion plane $\cB{\,=\,}\cF_\iota$\/}.
\par\smallskip
{\tt Proof.} $\cF_\Delta{\,=\,}\{W\}$ implies that $\iota$ is not a reflection.
Put $\Delta^{\hskip-2pt\ast}{\,=\,}\Delta|_\cB{\,=\,}\Delta/\Kappa$, and note that 
$\cF_{\Delta^{\hskip-1pt\ast}}{\,=\,}\{W\}$. Stiffness shows $\Kappa^1{\,\cong\,}\Spin3\RR$ or 
$\dim\Kappa{\,\le\,}1$. If $\Delta^{\hskip-2pt\ast}$ is semi-simple or has a normal torus, in particular, 
if $\Theta^*{\,\ne\,}\1$, then $\Delta{\hskip1pt:\hskip1pt}\Kappa{\,\le\,}13$ by 
\cite{sz1} 3.\hskip1pt1,2, 
and $\dim\Delta{\,\le\,}16$. Hence $\Delta^{\hskip-2pt\ast}$ has a normal vector subgroup, 
$\Theta^*{\,=\,}\1$, $\,\Theta{\,\trianglelefteq\,}\Kappa$, and $\Theta{\,=\,}\Kappa^1$. 
From \cite{sz1} 3.3 it follows that $\cB$ is classical or $\dim\Delta^{\hskip-2pt\ast}{\,\le\,}16$. 
In the second case, $\dim\Delta{\,<\,}18$. \qed
\par\smallskip
{\tt Remark.} Even if $\Delta$ induces on $\cB$ the full affine group or if $\Delta$ is an
{\it affine Hughes group\/} (i.e., the $28$-dimensional stabilizer of an interior line in a proper 
Hughes plane), it seems to be difficult to describe all the corresponding planes, 
cf. \cite{sz1} 3.2 Remark. 
\par\medskip
{\bf 4.5 Normal vector group.} {\it If $\Theta{\,\cong\,}\RR^t$ is a minimal normal subgroup of 
$\Delta$ and if $\dim\Delta{\,\ge\,}33$, then the connected component $\Tau$ of the
translation group $\Delta_{[W,W]}$ is a vector group of dimension at least $2$. Moreover, 
$\Theta{\,\le\,}\Tau$, or $t{\,=\,}8$ and $\dim\Tau{\,\ge\,}7$, or $t{\,=\,}12$\/}\, (\cite{sz10}).
\par\smallskip
{\tt Remarks on the proof.} (a) {\it If a one-parameter subgroup $\Pi{\,\le\,}\Theta$ is straight, 
then $\Pi{\,\le\,}\Tau$\/}. As $\Pi$ is not compact,  $\cF_\Pi$ is not a Baer subplane, and
$\Pi{\,\le\,}\Delta_{[z]}$ for some center $z$. In fact, $z{\,\in\,}W$ and $\Pi{\,\le\,}\Delta_{[z,W]}$, 
because $z^\Delta{\,\ne\,}z$ and $\Pi^\Delta$ is  contained in the commutative group $\Theta$. 
From $z^\Delta{\,\ne\,}z$ it follows also that $\Tau$ is commutative, $\Pi{\,<\,}\Tau$ and
$\dim\Tau{\,>\,}1$. \\
(b) {\it If some one-parameter subgroup $\Pi{\,\le\,}\Theta$ is not straight and if  $t{\,<\,}8$, then
$\dim\Delta{\,\le\,}32$\/}. Some orbit $b^\Pi$ is not contained in a line, and $\langle b^\Pi\rangle$
is a subplane of $\cP$. Choose $\rho{\,\in\,}\Pi\sm\{\1\}$, put $\Lambda{\,=\,}(\Delta_{b,\rho})^1$, 
and  note that $\Lambda{\,\le\,}\Cs{}\Pi$. Hence $\Lambda|_{\langle b^\Pi\rangle}{\,=\,}\1$, and the 
stiffness result 2.6(e) shows that $\Lambda{\,\cong\,}\Gtwo$ or $\dim\Lambda{\,\le\,}8$. 
The dimension formula yields 
$$\dim\Delta{\,=\,}\dim b^\Delta{\,+\,}\dim\Delta_b{\,=\,}\dim b^\Delta{\,+\,}\dim\rho^{\Delta_b}{\,+\,}
\dim\Lambda{\,\le\,}16{\,+\,}t{\,+\,}\dim\Lambda\,.$$
In the second case $\dim\Delta{\,<\,}32$. If $\Lambda{\,\cong\,}\Gtwo$, then 
$\Lambda{\,\le\,}\Cs{}\Theta$, since $\Pi^\Lambda{\,=\,}\Pi$ and $\Lambda$ has irreducible representations only in dimensions $7,14$ or larger. It follows that $\langle b^\Theta\rangle$
is contained in the flat plane $\cF_\Lambda$. Note that $\Theta_b{\,\le\,}\Delta_{b,\rho}$,  that
$\Theta_{\hskip-.5pt b}^{\hskip1pt1}{\,\le\,}\Lambda$, and that $\Lambda$ is compact. Hence 
$\Theta_{\hskip-.5pt b}^{\hskip1pt1}{\,=\,}\1$, 
$\,\dim\Theta_b{\,=\,}0$,  $\,\dim b^\Theta{\,=\,}t{\,\le\,}2$, and $\dim\Delta{\,\le\,}32$ as claimed.  \\
(c) For each $t{\,\ge\,}8$ the proof in \cite{sz10} uses the irreducible representation of $\Delta$ 
on $\Theta$. The arguments are different for distinct values of $t$. They are rather involved and 
shall not be reproduced here.
\par\medskip

{\bf 4.6 Translation group.} {\it If $\cF_\Delta{\,=\,}\{W\}$, if $\Delta$ has a normal vector subgroup, 
and if $\dim\Delta{\,\ge\,}35$, then the translation group $\Tau$ is transitive\/}\, (\cite{sz8}).
\par\smallskip
{\tt Remarks on the proof.} (a)  $\Tau$ contains a minimal normal subgroup $\Theta{\,\cong\,}\RR^t$ 
of $\Delta$. In a first step, we show that $t{\,\ge\,}6$. Choose some point $a{\,\notin\,}W$ and put
$\Gamma{\,=\,}(\Delta_a)^1$. For $\rho,\rho'$ on one-parameter groups in $\Theta$ with different centers,  the group $\Lambda{\,=\,}(\Gamma_{\hskip-2pt\rho,\rho'})^1$ fixes a connected subplane and Stiffness applies. We have $19{\,\le\,}\dim\Gamma{\,\le\,}2t{\,+\,}\dim\Lambda{\,\le\,}2t{\,+\,}14$ and $t{\,>\,}2$.   If $\Lambda{\,\cong\,}\Gtwo$, then $\Theta{\smcap}\Cs{}\Lambda{\,\cong\,}\RR^2$. Hence  $t{\,\ge\,}9$ or $\dim\Lambda{\,\le\,}8$ and $t{\,\ge\,}6$. \\
(b) Suppose that $t{\,=\,}6$. Then $\Lambda{\,\not\hskip1pt\cong\,}\SU3\CC$ and 
$\dim\Lambda{\,<\,}8$.  Consequently $\dim\rho^\Gamma{\,=\,}6$, $\,\rho^\Gamma$ is open in 
$\Theta$ for each $\rho$,  and $\Gamma$ is transitive on $\Theta\sm\{\1\}$. Similarly, $\Gamma$ is doubly transitive on the  $5$-dimensional projective space ${\rm P}\Theta$ consisting of the one-dimensional subspaces  of $\Theta$. Put $\tilde\Gamma{\,=\,}\Gamma|_\Theta{\,=\,}\Gamma/\Kappa$, and note that   $\langle a^\Theta\rangle{\ledot\hskip-3pt}\cP$. From \cite{cp} 96.17 it follows that 
$\tilde\Gamma$ is isomorphic to $\PSL3\CC$. Therefore $\dim\Kappa{\,=\,}3$ and  
$\langle a^\Theta\rangle{\,=\,}\cB$ is a Baer subplane. \cite{sz9} Th.~3.3 or 
\cite{sz1} 7.3 implies 
that $\cB$ is the classical quaternion plane.  Hence $\tilde\Gamma$ would be contained in 
$\HH^{\times}{\cdot\,}\SL2\HH$, which is impossible. \\
(c) Similarly, for each $t{\,<\,}16$ a contradiction can be derived from $\dim\Tau{\,<\,}16$ and  
a detailed analysis of the   representaion of $\Delta$ on $\Theta$ in combination with the consequences of the assumptions  on $\Delta$ as a group of automorphisms of $\cP$. The 
complicated proofs in \cite{sz8} shall not be repeated here.
\par\medskip
{\bf 4.7 Lemma.}  {\it Suppose that $\cF_\Delta{\,=\,}\{W\}$ and that $\dim\Delta{\,\ge\,}35$. 
If $\Theta{\,\cong\,}\RR^t$ is a minimal normal subgroup of $\Delta$ and if $\Theta{\,\le\,}\Tau$,
then   $8{\,\le\,}t\equiv0\bmod2$\/}.
\par\smallskip
{\tt Remark.} The claim is contained in steps (1)--(10) of the proof in \cite{sz8};  note that these 
steps do not use the assumption $\dim\Tau{\,<\,}16$.
\par\medskip
{\bf 4.8 Translation complement.} {\it Under the assumptions of {\rm 4.6}, the plane $\cP$ is the classical octonion plane or a complement $\Gamma{\,=\,}\Delta_a$ of the translation group 
$\Tau$ contains one of the groups $\SL2\HH{\,\cdot\,}\SU2\CC$ or $\SU4\CC{\,\cdot\,}{\rm H}$ with 
${\rm H}{\,=\,}\SU2\CC$ or $\SL2\RR$, and $\dim\Delta{\,=\,}35$\/}.
\par\smallskip
{\tt Remarks.} A proof of Theorem 4.8 is due to H\"ahl \cite{Ha1}, who also determined all planes 
such that $\Gamma$ has a factor $\SU4\CC$, see \cite{Ha2} and 4.12 below.  By different arguments we will show the following  weaker result:
\par\medskip
{\bf 4.9 Proposition.} {\it Suppose that $\cP$ is not classical. Under the assumptions of 
{\rm 4.7} and {\rm 4.8}, either $\Theta{\,\cong\,}\RR^8$ and $\Gamma'{\,\cong\,}\HH'{\cdot\hskip1pt}\SL2\HH$, or $\Theta{\,=\,}\Tau$ and $\Gamma'$ has a proper semi-simple factor $\Beta$ acting faithfully and irreducibly on some $8$-dimensional subspace $\Xi{\,<\,}\Tau$\/}.
\par\smallskip
{\tt Proof.} (a) We may assume that $\dim\Delta{\,<\,}38$ and $19{\,\le\,}\dim\Gamma{\,\le\,}21$.
In fact, coordinatizing quasi-fields of all {\it translation planes\/} with a group of dimension at least
$38$ have been determined by H\"ahl, see \cite{cp} 82.28 and \cite{Ha3}. There are two large 
families of such quasi-fields, the so-called generalized mutations and the perturbations of
the octonion algebra~$\OO$. In the first case, it follows from \cite{cp} 82.29 that the full 
automorphism group $\Sigma$ of the corresponding plane fixes a flag. Now let $\rho$ be a homeomorphism of the interval $[0,\infty)$ with $1^\rho{\,=\,}1$. The multiplication 
of the {\it perturbation\/} $\OO^{(\rho)}{\,=\,}(\OO,+,\circ)$ is defined by 
$$c{\,\circ\,}z{\,=\,}c\hskip1pt z_0{\,+\,}\|c\|^{\rho-1}\hskip1pt c\hskip1.5pt\frz \quad {\rm and}  \quad  
0{\,\circ\,}z{\,=\,}0\,.$$
Obviously, each automorphism of $\OO$ is also an automorphism of $\OO^{(\rho)}$. The plane 
over $\OO^{(\rho)}$ admits even a subgroup $\Phi{\,\cong\,}\Spin7\RR$ fixing the points
$(0)$ and $(\infty)$ on $W$, see \cite{cp} 82.5(a). It follows that $(0)$ and $(\infty)$ are fixed 
points of the full automorphism group if the plane is not classical. Thus in both cases 
$\cF_\Delta{\,\ne\,}\{W\}$. Because of 4.6, we may also assume that $\Gamma$ contains a 
one-parameter group $\Rho$ of homologies. \par
(b)  Of the  cases $t{\,\in\,}\{10,12,14\}$ of Lemma 4.7 the last one is quite easy: the commutator  
group $\Gamma'$ is semi-simple by \cite{cp} 95.6,  moreover,  $17{\,\le\,}\dim\Gamma'{\,\le\,}20$. According to \cite{cp} 95.10, the group  $\Gamma'$ is not almost simple, and   
$\Gamma'{\,=\,}\Alpha\Beta$ is a product of an almost simple factor 
$\Alpha$ of minimal dimension and the  semi-simple centralizer 
$\Beta{\,=\,}\Gamma'{\hskip-2pt\smcap}\Cs{}\Alpha$, and $\dim\Beta{\,\ge\,}10$. 
Clifford's Lemma \cite{cp} 95.5 implies that  $\Beta$ acts faithfully and 
irreducibly on $\RR^7$. Now $\Beta$ is a group of type $\Gtwo$ and dimension $14$,\: 
$\Gamma{:\,}\Gamma'{\,<\,}2$,\:   $\dim\Gamma'{\,>\,}17$, and then $\dim\Alpha{\,=\,}6$ and 
$\Alpha{\,\circeq\,}\SL2\CC$. Again by  Clifford's Lemma, $\Alpha$ acts faithfully and irreducibly 
on $\RR^s$ with $s|14$, but this contradicts \cite{cp} 95.10 and  shows that $t{\,\ne\,}14$. \\
(c) Similarly, $t{\,=\,}10$ is impossible: note that $\Gamma$ acts faithfully on $\Theta$. The only almost simple group in the dimension range for 
$\Gamma'$  admitting a $10$-dimensional irreducible representation is $\SO5\CC$, but this group
contains $\SO5\RR$ which cannot act on $\cP$, see 2.10. Again $\Gamma'{\,=\,}\Alpha\Beta$  is a product of proper factors  with $\dim\Beta{\,\ge\,}10$. Clifford's Lemma implies that  $\Beta$ acts faithfully and  irreducibly on $\RR^5$, and then $\Beta$ is a group $\Opr5(\RR,r)$ with $r{\,>\,}0$.
The fixed elements in $\Theta$ of a suitable subgroup  of $\Beta$ form an $\Alpha$-invariant 
$2$-dimensional vector space on which $\Alpha$ acts  faithfully, $\dim\Alpha{\,=\,}3$ and 
$\dim\Gamma'{\,=\,}13$ contrary to the hypothesis.  \\
(d) The case $t{\,=\,}12$ is more complicated. An almost simple group $\Gamma'$ would be 
locally isomorphic to $\Sp4\CC$ and has no irreducible representation on $\RR^{12}$. 
Therefore  $\Gamma'{\,=\,}\Alpha\Beta$, where $\Alpha$ is a factor of minimal dimension and 
$10{\,\le\,}\dim\Beta{\,\le\,}17$. By complete reducibility, $\Theta$ has a $\Gamma'$-invariant 
complement $\Chi$ in $\Tau$. Suppose that $a^\Chi$ is contained in a line $av, \ v{\,\in\,}W$. 
Then $v^{\Gamma'}{\hskip-2pt=\,}v{\,\ne\,}v^\zeta$ for some $\zeta$ in the center of $\Gamma$, 
and $\Gamma'$ fixes the triangle $a,v,v^\zeta$.  For $\tau{\,\in\,}\Chi$,   the centralizer 
$\Lambda{\,=\,}\Gamma'{\hskip-2pt\smcap}\Cs{}\tau$ fixes also the  point $a^{\tau^\zeta}$, 
hence a quadrangle,  $\dim\Lambda{\,\ge\,}13$, and Stiffness yields  
$\Lambda{\,\cong\,}\Gtwo$ and  $\Lambda|_\Chi{\,=\,}\1$, but $\dim\cF_\Lambda{\,=\,}2$. \\
(e) Consequently, $\langle a^\Chi\rangle{\,=\,}\cF$ is a subplane, and $\cF{\,<\,}\cP$\, (or 
$\Gamma'$  would act faithfully on $\Chi{\,\cong\,}\RR^4$). 
Put $\Gamma'|_\cF{\,=\,}\Gamma'/\Kappa$, and note that  $\Gamma'|_\cF{\,\le\,}\SL4\RR$. 
If $\dim\cF{\,=\,}4$, then Stiffness implies 
$\dim\Kappa{\,\le\,}8$ and $\Gamma'{\hskip1pt:\hskip1pt}\Kappa{\,>\,}8$, which is impossible. 
Hence $\cF{\,\ldot}\cP$,\, $\Alpha{\,=\,}\Kappa{\,\cong\,}\SU2\CC$,\, $\dim\Beta{\,\ge\,}14$, and  
$\Beta{\,\cong\,}\SL4\RR$. If $\tau{\,\in\,}\Tau$  is a translation of $\cP$ such that 
$a^\tau{\,\in\,}\cF$, then $a^{\tau\Chi}{\,\subseteq\,}\cF^\Chi{\,=\,}\cF$ and 
$\cF^\tau{\,=\,}\langle a^\Chi\rangle^\tau{\,=\,}\langle a^{\tau\Chi}\rangle{\,=\,}\cF$. 
Therefore $\cF$ is a translation plane,  $\dim\Aut\cF{\,\ge\,}23$, and  $\cF$ is the   classical 
quaternion plane by  \cite{cp} 81.9 or 84.27. (In fact, an $8$-dimensional plane is classical, 
if it has an automorphism group of dimension ${\ge\,}19$, see 
\cite{sz1} 1.10.) It follows that 
$\Beta{\,\cong\,}\SL2\HH$, a contradiction.
This proves that $\Theta{\,\cong\,}\RR^8$ or  $\Theta{\,=\,}\Tau$. \\
(f) {\it  If $\Gamma$\,acts {\rm faithfully} and irreducibly  on $\Theta{\,\cong\,}\RR^8$, then 
$\Gamma'{\,=\,}\Alpha\Beta$ with  $\Alpha{\,\cong\,}\HH'$ and $\Beta{\,\cong\,}\SL2\HH$\/}. 
{Proof:} The   group $\Gamma'$ is semi-simple by \cite{cp} 95.6,  and again $17{\,\le\,}\dim\Gamma'{\,\le\,}20$. If  $\Gamma'$ is almost simple,  
then $\dim\Gamma'{\,=\,}20$.  The list \cite{cp} 95.10 of irreducible representations shows 
that $\Gamma'{\,\cong\,}\Sp4\CC{\,\cong\,}\Spin5\CC$. The central involution in $\Gamma'$ 
is a reflection with center $a$ because $\Gamma$ has no fixed point on $W$. The involution 
$\omega{\,=\,}{\rm diag}\,(1,1,-1,-1)$ is either planar or a reflection  with some center $u{\,\in\,}W$, 
and $\Omega{\,=\,}\Gamma'\smcap\Cs{}\omega{\,\cong\,}(\Sp2\CC)^2{\,\cong\,}(\SL2\CC)^2$. 
If $\cF_\omega{\,\ldot}\cP$, then $\Omega$ induces on $W\smcap\cF_\omega{\,\approx\,}\sS_4$ 
a group $({\rm P})\SL2\CC{\times}\SO3\CC$, and a maximal compact subgroup 
$\Phi$ of $\Omega$  acts  as $\SU2\CC{\times}\SO3\RR$ or as $(\SO3\RR)^2$ in contradiction 
to Richardson's Theorem  \cite{cp} 96.34. If $\omega{\,\in\,}\Delta_{[u,av]}$, then $\Omega$ fixes the triangle $a,u,v$, and  \cite{cp} 81.8 implies that $\Omega$ has a compact subgroup of codimension at most $2$, but $\Omega{\,:\,}\Phi{\,=\,}6$. Hence  $\Gamma'{\,=\,}\Alpha\Beta$ is a product of an almost simple factor  $\Alpha$ of minimal dimension and the  semi-simple centralizer 
$\Beta{\,=\,}\Gamma'{\hskip-2pt\smcap}\Cs{}\Alpha$, and $\dim\Beta{\,\ge\,}10$. By Clifford's 
Lemma \cite{cp} 95.5,  $\Beta$ acts faithfully and irreducibly on a subspace $\Xi$ of $\Theta$ of 
dimension $4$ or~$8$. In the second case, $\Alpha{\,\cong\,}\HH'$ and $\Beta{\,\cong\,}\SL2\HH$. 
If $\Xi{\,\cong\,}\RR^4$, then $\Beta$ is isomorphic to $\Sp4\RR$ or $\SL4\RR$. A suitable subgroup  of  $\Beta$  fixes a one-parameter subgroup of $\Xi$ and of each $\Beta$-invariant  complement of  $\Xi$ in $\Theta$, and $\Alpha$ acts on a $2$-dimensional subspace of $\Theta$. 
Consequently, $\Alpha{\,\cong\,}\SL2\RR$ and $\Beta{\,\cong\,}\SL4\RR$. 
If $a^\Xi{\,\subseteq\,}av$ with $v{\,\in\,}W$, then $v^\Beta{\,=\,}v{\,\ne\,}v^\Alpha$, and $\Beta$ 
fixes a triangle $a,u,v$. By \cite{cp} 81.8, a maximal compact subgroup $\Phi$ of $\Beta$ has
codimension $\Beta{\,:\,}\Phi{\,\le\,}2$, but $\Phi{\,\cong\,}\SO4\RR$ is much too small. 
Thus $\Xi$ contains translations $\xi,\eta$ with different centers,  the stabilizer 
$\Lambda{\,=\,}(\Beta\Rho)_{\hskip-1pt a^\xi,\hskip1pt a^\eta}$  fixes a quadrangle, and Stiffness 
implies $\Lambda{\,\cong\,}\SU3\CC$, but then $\Lambda$ is not contained in $\SL4\RR$. \par
(g) Suppose now that $\tilde\Gamma{\,=\,}\Gamma|_\Theta{\,=\,}\Gamma/\Kappa$ with 
$\Kappa{\,\ne\,}\1$.  Then $\langle a^\Theta\rangle{\,\le\,}\cF_\Kappa$ is a Baer subplane since 
$z^\Gamma{\ne\,}z$  for each $z{\,\in\,}W$. Either $\Kappa{\,\cong\,}\Spin 3\RR{\,\cong\,}\HH'$ 
or $\dim\Kappa{\,\le\,}1$. In the second case, $\dim\tilde\Gamma{\,\ge\,}18$ and 
$\dim\Aut\cF_\Kappa{\,\ge\,}26$. By \cite{cp} 84.27 the plane  $\cF_\Kappa$ 
is the classical quaternion plane.  Hence $\tilde\Gamma$ is contained in 
$\HH^{\times}{\cdot\,}\SL2\HH$ and $\Gamma$ has a subgroup $\SL2\HH$.
Moreover, $\dim(\tilde\Gamma\smcap\HH'){\,\ge\,}2$ and $\HH'{\,\le\,}\Gamma$. Consequently 
$\Gamma'{\,\cong\,}\HH'{\cdot\,}\SL2\HH$ \par
(h) Only the case that $\Gamma$ acts irreducibly on $\Theta{\,=\,}\Tau{\,\cong\,}\RR^{16}$ remains.
If $\Gamma'$ is almost simple, then $\dim\Gamma'{\,=\,}20$, and $\Gamma'$ is isomorphic to 
$\SO5\CC$ or to its covering group $\Sp4\CC$. In the first case, $\Gamma$ would contain 
$\SO5\RR$ contrary to 2.10. In the second case, $\Tau$ is a direct sum of two subspaces on
which $\Gamma'$ acts in the natural way (see \cite{cp} 95.\,6 and 10), and
$\Gamma{\,=\,}\Gamma'\Rho$ would not be irreducible on $\Tau$. Therefore 
$\Gamma'{\,=\,}\Alpha\Beta$  is again a product of proper factors such that  $\Alpha$ has 
minimal dimension, and  $10{\,\le\,}\dim\Beta{\,\le\,}17$. By Clifford's  Lemma, $\Beta$ acts 
faithfully and irreducibly on a subspace $\Xi$ of $\Tau$ of dimension $k|16$. \\
(i)  In the case $k{\,=\,}4$, we have  either $\Beta{\,\cong\,}\Sp4\RR$ and  $\dim\Alpha{\,\ge\,}8$ 
or  $\Beta{\,\cong\,}\SL4\RR$ and  $\dim\Alpha{\,=\,}3$. First,  let $\dim\Beta{\,=\,}10$. 
The centralizer $\Omega$ of the involution  $\omega{\,=\,}{\rm diag}(1,1,-1,-1){\,\in\,}\Beta$ 
contains  $\Alpha{\hskip.5pt\cdot\hskip.5pt}(\SL2\RR)^2$. 
If $\omega$ is a reflection, then $\Omega$ fixes a triangle $a,u,v$. By \cite{cp} 81.8, a maximal 
compact subgroup $\Phi$ of $\Omega$   satisfies $\Omega{\,:\,}\Phi{\,\le\,}2$, but 
$\Phi{\,:\,}\Alpha{\,\le\,}2$ and $\Omega{\,:\,}\Alpha{\,=\,}6$. 
Hence $\cF_\omega{\,\ldot}\cP$ and $\Omega$ induces on $\cF_\omega$  a group of dimension 
at least $14$; moreover, $\Tau\smcap\Cs{}\omega{\,\cong\,}\RR^8$ and $\cF_\omega$ is a translation plane. Therefore $\cF_\omega$ is the classical quaternion plane, cf. \cite{cp} 81.9.  
It follows that $\Omega|_{\cF_\omega}{\,\cong\,}\SL2\HH$, but $\Alpha$ is a proper 
factor  of~$\Omega$.  Consequently, $\dim\Beta{\,=\,}15$ and $\Gamma$ has a subgroup 
$\Ypsilon{=\,}\Beta\Rho{\,\cong\,}\GL4\RR$   acting on $\Xi$. 
If $a^\Xi{\,\subseteq\,}av$ with $v{\,\in\,}W$, then $v^\Ypsilon{=\,}v{\,\ne\,}v^\Alpha$, and 
$\Ypsilon$ fixes a triangle $a,u,v$, but the maximal compact subgroup $\SO4\RR$ of 
$\Ypsilon$ is too small (cf. \cite{cp} 81.8). Thus $\Xi$ contains tranlations $\xi,\eta$ with 
different centers,  the stabilizer 
$\Lambda{\,=\,}\Ypsilon_{\hskip-1pt a^\xi,\hskip1pt a^\eta}$  fixes a quadrangle, and Stiffness 
implies $\Lambda{\,\cong\,}\SU3\CC$, an obvious contradiction. Hence $k{\,\ge\,}8$. \\ 
(j) Suppose next that $k{\,=\,}16$. As $\Alpha{\,\le\,}\Cs{}\Beta$, it follows from \cite{cp} 95.4 
that $\Alpha{\,\cong\,}\HH'$. Therefore $14{\,\le\,}\dim\Beta{\,\le\,}17$, and   $\Beta$ acts 
irreducibly on $\HH^4$. 
The list \cite{cp} 95.10 shows that $\Beta$ is not almost simple. Clifford's Lemma implies that 
each proper factor $\Nu$ of $\Beta$ is contained in $\SL2\HH$, in particular, $\Nu$ has torus 
rank $\rk\Nu{\,\le\,}2$, and $\Beta$ is composed of factors $\HH'$,\, \ $\U2(\HH,r)$, and 
$\SL2\HH$, but then the conditions on rank and dimension cannot be satisfied simultaneously. \qed
\par\medskip
Based on his previous papers \cite{Lh1} and \cite{Lh2}, L\"owe \cite{Lh} explicitly determined 
not only the planes of the first case in 4.9, but all translation planes admitting a group  locally isomorphic to $\SL2\HH$. 
\par\smallskip
{\bf 4.10 Construction.} Choose a {\it continuous\/} map $\mu{\,:\,}\HH'{\,\to\,}e^{\RR}$. 
Write $\eta^*{\,=\,}(\overline\eta^{\hskip1.5pt t})^{-1}$ and put
$\Psi{\,=\,} \bigl\{\bigl({\eta^* \atop }{ \atop \eta}\bigr){\,\big|\,}\eta{\,\in\,}\SL2\HH \bigr\}$. 
For $h{\,\in\,}\HH$, let $L_h{\,=\,}\{(x,hx)\mid x{\,\in\,}\OO\}$. 
Then  $E{\,=\,}\{(\HH,\ell\HH){\,\subset\,}\OO^2\}$ and 
 the {\it ``lines''\/} $L_{h^\mu h}$ with $h{\,\in\,}\HH'$ together with all their $\Psi$-images  
 form a spread $\cS_\mu$. Obviously $\Psi$ acts on the translation plane $\cP_\mu$
 defined by this spread, see \cite{Lh} for proofs.
\par\smallskip
Let again $\Gamma{\,=\,}\Delta_a$ denote a complement of the translation group, and consider the subgroup $\Omega{\,=\,}{\rm S}\Gamma{\,=\,}\{\gamma{\,\in\,}\Gamma{\,\le\,}\GL{16}\RR\mid
\det\gamma{\,=\,}1\}$.
In our terminology, L\"owe's main result \cite{Lh} 3.8 can be stated as follows:
\par\smallskip
{\bf 4.11 \,$\SL2\HH$-planes.} {\it Let $\cP$ be a non-classical translation plane. If $\Delta$ 
contains a subgroup  locally isomorphic to $\SL2\HH$, then $\cP$ is isomorphic to $\cP_\mu$, where 
$\mu{\,:\,}\HH'{\,\to\,}e^{\RR}$ is continuous and not constant. There are $3$ 
possibilities\/}: 
\par\smallskip
(1) {\it If $\mu$ depends only on the real part of its argument, then $\Omega$ is 
an almost direct product of  $\Psi{\,\cong\,}\SL2\HH$ and
$\Eta{\,=\,}\{(x,y){\,\to\,}(hx,hy)\mid h{\,\in\,}\HH'\}$. In particular,  $\dim\Sigma{\,=\,}35$\/}.
\par\smallskip
(2)  {\it  If $\mu$ depends only on the absolute value of the complex part of its argument, then
 $\Omega$ is an almost direct product of  $\Psi{\,\cong\,}\SL2\HH$ and  
 $\Alpha{\,=\,}\{(x,y){\,\to\,}(bx,cy)\mid b,c{\,\in\,}\HH'{\smcap\CC}\}{\,\cong\,}\TT^2$. 
In particular,  $\dim\Sigma{\,=\,}34$\/}.
\par\smallskip
(3) {\it In all other cases,  $\Omega$ is isomorphic to $\Psi{\,\cong\,}\SL2\HH$ or to  an almost direct product of  $\Psi$ and $\SO2\RR$. Hence  $32{\,\le\,}\dim\Sigma{\,\le\,}33$\/}.
\par\medskip
{\bf 4.12 \,$\SU4\CC$-planes.} As mentioned above, H\"ahl \cite{Ha2} has described explicitly 
all translation planes admitting a semi-simple group $\Ypsilon{\,=\,}\Phi\Psi$ with 
$\Phi{\,\cong\,}\SU4\CC$ and $\Psi{\,\in\,}\{\SU2\CC,\SL2\RR\}$. Only  some facts shall be 
stated; for the full theorem the reader is referred to \cite{Ha2}. \par\smallskip
(1) {\it Let the translation group $\Tau$ be written in the form $\CC^4{\,{\otimes}_{\CC}\,}\CC^2$. 
The group $\Phi$ acts in the standard way on the left factor and trivially on the right one; $\Psi$ acts  only on the right factor. The subspaces  
$\CC^4{\,{\otimes}_{\CC}\,}z$,\: $z{\,\in\,}\CC^2\sm\{0\}$,  are affine lines\/}. 
\par\smallskip
(2) {\it With $\CC^2{\hskip1pt=\,}\HH$, the map 
 $h{\,{\otimes}_{\CC}\,}(x,y){\,\mapsto\,}h{\,{\otimes}_{\HH}\,}(x,y)$ yields an isomorphism 
$\Tau{\,\cong\,}\HH^2{{\otimes}_{\HH}\,}\HH^2$. Let $\U2\HH{\,=\,}${\cyss Yu}${\,<\,}\Phi$. 
The {\cyss Yu}-invariant subspaces of $\Tau$ are exactly 
$U_q{\,=\,}\HH^2{\,{\otimes}_{\HH}\,}(1,q)$ and 
$U_\infty{\,=\,}\HH^2{\,{\otimes}_{\HH}\,}(0,1)$\/}. \par\smallskip
(3) {\it Let $\sigma{\,:\,}(-1,1){\,\to\,}\{c{\,\in\,}\CC\mid c\overline c{\,=\,}1\}$ be any continuous map 
{\rm(}which satisfies $\sigma(-s){\,=\,}\sigma(s)$ if $\Psi$ is compact{\rm)}. Then $\sigma$ 
selects a $\Ypsilon\!$-invariant spread among the subspaces described in\/} (1) {\it and\/} (2) 
{\it containing the affine lines $U_q$, where $q$ has the form $si{\,+\,}tj\sigma(s){\,\in\,}\HH$, 
$s^2{\,+\,}t^2{\,=\,}1$ and $s,t{\,\in\,}(-1,1)$,  their $\Ypsilon\!$-images, and the lines 
$U_c$, $c{\,\in\,}\CC\smcup\{\infty\}$\/}.
\par\smallskip
(4) {\it  Each translation plane admitting an action of the group $\Ypsilon$ can be obtained in this 
way. The plane is classical if, and only if, $\sigma$ is constant; in all other cases 
$\dim\Sigma{\,=\,}35$\/}.
\par\bigskip	
{\Bf 5. Exactly two fixed elements}
\par\medskip
For comparison with later results, \cite{cp} 87.4 shall be restated here:
\par\medskip
{\bf  Proposition.}  {\it Assume that $\cP$ is not classical.  If $\dim\Delta{\,=\,}h{\,\ge\,}36$, and if $\Delta$ has a minimal normal vector subgroup $\Theta{\,\cong\,}\RR^t$,  then there are only the following possibilities\/}:
\par\medskip
$$\begin{cases} \Theta{\,\le\,}\Delta_{[a,W]}, \ a{\,\notin\,}W 
\begin{cases} \Delta'{\,\cong\,}\Spin9(\RR ,r),\; r{\,\le\,}1 \hskip8pt and \hskip8pt  h{\,=\,}37 \cr
v^\Delta{\,=\,}v{\,\in\,}W \hskip8pt  and \hskip8pt  h{\,\le\,}38 \cr \end{cases} \cr
\Theta{\,\le\,}\Tau{\,=\,}\Delta_{[W,W]}
\begin{cases} t{\,>\,}8 \begin{cases} h{\,=\,}36 \cr
\Tau{\,\cong\,}\RR^{16} \hskip8pt  and \hskip8pt  h{\,=\,}37 \cr \end{cases}\cr
t{\,=\,}8 \Rightarrow \Theta{\,=\,}\Delta_{[v,W]} \hskip8pt  or \hskip8pt  \Tau{\,\cong\,}\RR^{16} \cr
t{\,<\,}8 \Rightarrow \Theta{\,\le\,}\Delta_{[v,W]},\; v{\,\in\,}W\;. \cr \end{cases}
\end{cases}$$
\par\medskip
The case $h{\,\ge\,}40$ has already been described in the introduction, for $h{\,=\,}39$ it follows from  the Proposition that $\cF_\Delta$ is a flag. If a semi-simple group $\Delta$ of the octonion plane $\cO$ fixes a flag, then $\Delta$ is  contained in $\Spin8\RR$ and hence fixes even a triangle. 
Priwitzer's Theorem 4.2 shows that  $\cF_\Delta{\,=\,}\{a,W\}$ with $a{\,\notin\,}W$, if $\Delta$ is a semi-simple group of a non-classical plane, $\cF_\Delta{\,\ne\,}\emptyset$, 
and  $\dim\Delta{\,>\,}28$.
\par\medskip
{\Bf A. Fixed flag}
\par\medskip
{\bf 5.0.} {\it If $\cF_\Delta{\,=\,}\{v,W\}$ is a flag, and if $\dim\Delta{\,\ge\,}22$, then $\Delta$ is a Lie group\/}. 
\par\smallskip
This can be {\tt proved} in a similar way as 4.0. Steps (a) and the first part of step (b) are the 
same as before. In the present case the weaker assumption on $\dim\Delta$ yields 
$\dim\Delta|_\cD{\,\ge\,}15$, and $\cD$ is no longer known explicitly. 
Choose a line $L$ of $\cD$ in $\frL_v\sm\{W\}$ and points $p{\,\in\,}L,\,z{\,\in\,}W$ outside $\cD$.
Put $\Gamma{\,=\,}\Delta_L$ and $\Lambda{\,=\,}(\Gamma_{\hskip-2pt p,z})^1$. 
Recall that $\Delta|_\cD{\,=\,}\Delta/\Kappa$. If $\Kappa$ is not a Lie group, then 
$\dim p^\Gamma,\, \dim z^\Gamma{<\,}8$ by \cite{cp} 53.2, and 
$\dim\Lambda{\,\ge\,}22{-}4{-}2{\cdot}7{\,>\,}3$. The line $pz$ intersects $\cD$ in a unique point $a$. If $\dim\Kappa{\,\le\,}3$, then $\dim\Delta|_\cD{\,\ge\,}19$ and $\cD$ is classical 
(cf. \cite{cp} 84.28).
The group $\Lambda$ fixes the distinct lines $L,W,av$ and hence a connected set of lines in the 
pencil $\frL_v$. In particular, $\cF_\Lambda$ is connected. If $\dim\Kappa{\,>\,}3$, then 
$\Kappa^1$ has a connected commutative factor $\Alpha$ and $\Lambda$ acts trivially on 
$\Alpha$\: (see the arguments in step (c) in the proof of 4.0). Hence 
$p^\Alpha{\,\subseteq\,}\cF_\Lambda$ and again $\cF_\Lambda$ is connected. 
As $\Nu$ acts freely on $L\sm\cD$, it follows from \cite{cp} 32.21 and 71.2 that 
$\dim\cF_\Lambda{\,>\,}4$ and $\cF_\Lambda{\;\ldot}\cP$. Now $\Lambda{\,\cong\,}\Spin3\RR$
contrary to the statement above. \qed
\par\medskip
{\bf 5.1 Lie.} {\it If $\Delta$ fixes a flag $\{v,W\}$  and no line other than $W\hskip-3pt$, if $\Delta$ is 
semi-simple,  and if  $\dim\Delta{\,>\,}14$, then $\Delta$ is a  Lie group\/}.
\par\smallskip
The first part of the {\tt Proof} is identical with step (a) and the beginning of step (b) in the proof of 4.0.  Let  $\Delta|_\cD{\,=\,}\Delta/\Kappa$ and $\Lambda{\,=\,}\Kappa^1$\hskip-2pt. Note that 
$\Delta/\Kappa$ and $\Lambda$ are semi-simple. If $\dim\cD{\,=\,}2$, then $\Delta{\,=\,}\Lambda$ 
by  \cite{cp}~33.8,  and  $\dim\Lambda{\,\le\,}14$ by Stiffness;  if $\dim\cD{\,=\,}4$, then 
$\Delta{\,:\,}\Kappa{\,\le\,}3$ and $\dim\Lambda{\,\le\,}8$ by 2.14 and 2.6(e). 
Finally, if $\cD{\,\ldot}\cP$, then   $\Delta{\,:\,}\Kappa{\,\le\,}11$ by 
\cite{sz1} 7.3,  
$\,\Lambda$ is compact and semi-simple,  hence a Lie group (cf. \cite{cp} 93.11 and 94.20), 
and $\dim\Lambda{\,\le\,}3$ by 2.6(\^b).  Thus $\dim\Delta{\,\le\,}14$ or $\Nu{\,=\,}\1$. \qed
\par\medskip
{\bf 5.2 Semi-simple groups.} {\it If $\cF_\Delta$ is a flag $\{v,W\}$ and if $\Delta$ is semi-simple,
then $\dim\Delta{\,\le\,}21${\rm;} equality is possible only if each compact subgroup of $\Delta$ is 
trivial.   If $\Delta$ is almost simple, then $\dim\Delta{\,\le\,}16$\/}. 
\par\smallskip 
The following  {\tt proof} does not make full use of the assumption $\cF_\Delta{\,=\,}\{v,W\}$ but works also in the case  $\cF_\Delta{\,=\,}\{u,v,W{=\hskip2pt}uv\}$, see 6.2. \\
(a) {\it Each involution in $\Delta$ is planar\/}: if $\sigma$ is a reflection, then 
$\sigma^\Delta\sigma$ would be a normal vector subgroup by 2.9 or its dual. \\
(b) {\it If there exists a proper $\Delta$-invariant subplane $\cE$, then $\dim\Delta{\,\le\,}14$\/}. 
This follows in the same way as in the proof of 5.1; for discrete $\cE$, it is an immediate consequence of Stiffness. \\
(c)  Suppose now and in the next steps that $\dim\Delta{\,>\,}14$. Then  $\Delta$ does not contain a central involution, and $\Delta$ is a Lie group by 5.1. \\
(d) Let $\Gamma$ be an almost simple factor of $\Delta$ and denote the product of the other 
factors by~$\Psi$. In steps (e--h),\, $\Gamma$ is assumed to be a proper factor of minimal dimension. \\
(e) If $\Gamma$ contains an involution $\iota$, then $\Psi|_{\cF_\iota}{\,=\,}\Psi/\Kappa$ has dimension  $\Psi{\hskip1pt:\hskip1pt}\Kappa{\,\le\,}11$ by 
\cite{sz1} 7.3, the semi-simple 
group $\Kappa$ is compact, and   $\dim\Kappa{\,\le\,}3$ by Stiffness. 
Either $\dim\Gamma{\,=\,}3$ and $\dim\Delta{\,\le\,}17$, or $\dim\Kappa{\,=\,}0$,\,  
$\dim\Psi{\,\ne\,}11$, and  $6{\,\le\,}\dim\Gamma{\,\le\,}\dim\Psi{\,\in\,}\{6,8,10\}$. 
In any case,\,   $\dim\Delta{\,\le\,}20$. \\
(f)  If $\Gamma$ does not contain an involution, in particular, if $\dim\Delta{\,>\,}20$, then 
$\Gamma$ is isomorphic to the simply connected covering group  $\Omega$ of $\SL2\RR$, since $\Omega$ is the only almost simple Lie group without involution. \\
(g) {\it If $\Gamma$\hskip2pt is straight and if $\dim\Delta{\,>\,}14$, then 
$\Gamma$\hskip2pt is a group of translations and the center is a fixed point of $\Delta$\/}:
By Baer's theorem, either $\Gamma$ is a group of axial collineations and center and axis of 
$\Gamma$ are fixed elements, or $\cD{\,:=\,}\cF_\Gamma{\,=\,}\cD^\Delta{\;\ldot}\cP$, and $\dim\Delta{\,\le\,}14$ by step (b). \\
(h) {\it If no factor of $\Delta$ has a compact subgroup other than $\{\1\}$, then $\dim\Delta{\,\le\,}21$\/}. \\
In fact, $\Delta$ is a product of factors $\Gamma_{\hskip-1.5pt\nu}{\,\cong\,}\Omega$. At most 
two factors $\Gamma_{\hskip-1.5pt\nu}$ consist of translations (or they would be commutative
by \cite{cp} 23.13). Either $\langle x^\Delta\rangle$ is a proper subplane for some point $x$ and 
$\dim\Delta{\,\le\,}14$ by (b), or there are factors $\Gamma_{\hskip-1.5pt\nu}$ such that
$\langle x^{\Gamma_1}\rangle$ is a non-flat subplane (use \cite{cp} 33.8) and 
$\cH{\,=\,}\langle x^{\Gamma_1\Gamma_2}\rangle{\;\ledot\!}\cP$. In the latter case, put 
$\Chi{\,=\,}\Psi_1\smcap\Psi_2$. Then $\Chi_x|_\cH{\,=\,}\1$ and $\Chi_x$ is compact 
by 2.6(b), hence $\Chi_x{\,=\,}\1$ and $\dim\Chi{\,\le\,}15$\, ($\Chi$ being a product of
$3$-dimensional factors).\\
(i) From now on assume that $\dim\Delta{\,>\,}20$ and that $\Gamma$ contains an involution
$\iota$, but is not necessarily of minimal dimension. Step (e) implies 
$\Psi{\hskip1pt:\hskip1pt}\Kappa{\,\le\,}11$ and $\dim\Gamma{\,\ge\,}8$, and this is true for
each possible choice of $\Gamma$. Consequently $\dim\Kappa{\,=\,}0$, and then 
$\dim\Gamma{\,\ge\,}10$. In the case $\dim\Psi{\,=\,}11$, no factor of $\Psi$ contains an
involution and each factor of $\Psi$ is isomorphic to $\Omega$, but then $\dim\Psi{\,|\,}9$.
Hence $\dim\Gamma{\,>\,}10$ and $\dim\Gamma{\,\ge\,}14$. Again $\Psi$ is a product of factors isomorphic to $\Omega$ and $\dim\Psi{\,\le\,}9$.    
The kernel of the action of $\Psi$ on $\cF_\iota$ is compact and hence trivial (cf. 2.6(b)\,). 
As  $\rk\Gamma{\hskip1pt>\hskip1pt}1$, some other involution $\iota'{\,\in\,}\Gamma$ commutes 
with $\iota$, and  $\Psi$ acts also effectively on $\cF_{\hskip-1pt\iota'}$. Suppose that 
$\cF_{\hskip-1pt\iota,\iota'}$  is a subplane.  Then $\Psi$ is even almost effective on  
$\cF_{\hskip-1pt\iota,\iota'}$ by the stiffness result 2.6(c), and  2.14 implies 
$\dim\Psi{\,\le\,}3$,\, $\dim\Gamma{\,\ge\,}18$, and then  $\dim\Gamma{\,\ge\,}20$.
If $\overline{\iota'}{\,=\,}\iota'|_{\cF_\iota}$ is a reflection, then $\Psi$ fixes the center $a{\,\notin\,}W$ 
of $\overline{\iota'}$\, (up duality) or  its axis $av$ and center $u{\,\in\,}W$. In the first case,  
$\Psi|_{a^\Gamma}{\,=\,}\1$,\, $\langle a ^\Gamma\rangle{\,=\,}\cP$ by step~(b), and $\Psi{\,=\,}\1$; 
in the second case, $\Psi{\hskip1pt:\hskip1pt}\Psi_a{\,\le\,}4$,\, $u^\Gamma{\,=\,}u^\Delta{\,\ne\,}u$,\,
$\cE{\,=\,}\langle a^\Gamma,u^\Gamma\rangle{\,=\,}\cE^\Gamma{\ledot\!}\cP$ by 2.14, 
$\Psi_{\hskip-1pt a}|_\cE{\,=\,}\1$, $\Psi_{\hskip-1pt a}$ is compact, 
$\dim\Psi_{\hskip-1pt a}{\,=\,}0$, and $\dim\Psi{\,\le\,}3$. Again, $\dim\Gamma{\,\ge\,}20$. \\
(j) If $\dim\Gamma{\,\ge\,}30$,  Priwitzer's theorem \cite{pw1} shows the following: 
$\dim\Gamma{\,=\,}36$ and $\Gamma$ fixes a non-incident point-line pair, 
or $\Gamma{\,\cong\,}\SL3\HH$ and $\cF_\Gamma{\,=\,}\emptyset$, or  the plane is classical; in the last case $\cF_\Gamma$ is not a flag, since a maximal semi-simple subgroup in the stabilizer of a 
flag is isomorphic to $\Spin8\RR$. The case $|\cF_\Delta|{\,=\,}3$ is also covered by Theorem 6.2. Therefore, only groups $\Gamma$ with $\dim\Gamma{\,<\,}30$ have to be discussed. 
 2.16 shows that $\Gamma{\,\not\hskip1pt\cong\,}\Gtwo^{\hskip1pt\CC}$.
According to  step (c),\, $\Gamma$  does not contain  a central involution. 
By 2.10 and 2.11, neither $\SO5\RR$ nor $(\SO3\RR)^2$ can be a subgroup of~$\Gamma$.
These facts  exclude the group $\SO5\CC$ of dimension $20$ and leave only the following 
$9$ possibilities:  $\PU3(\HH,r)$ and $\PSp6\RR$ of type ${\rm C}_3$,\,  
${\rm(\hskip-.5pt P\hskip-1pt)}\SU5(\CC,r)$  of type ${\rm A}_4$.   
For each of these groups a contradiction will be derived.  \\
(k)  The symplectic group contains $3$ pairwise commuting conjugate  involutions. Such an involution  $\iota$ must be planar, its centralizer $\Ypsilon$ is locally isomorphic to 
$\SL2\RR{\times}\Sp4\RR$.  Stiffness implies $\dim\Ypsilon|_{\cF_\iota}{\,=\,}13$, which 
contradicts \cite{sz1} 7.3, see also  \cite{cp} 84.18. Similarly, the compact group $\PU3\HH$ 
contains a planar involution $\iota$ such that  $\Ypsilon{\,=\,}\Cs{}\iota$ is locally isomorphic to 
$\U2\HH{\times}\HH'{\,=\,}\Spin5\RR{\times}\Spin3\RR$. 
By \cite{cp} 96.13 or by Richardson's classification of compact (Lie) groups on $\sS_4$\, (see 
\cite{cp} 96.34), the $10$-dimensional factor of $\Ypsilon$  is transitive on 
$W{\smcap}\cF_\iota{\,\approx\,}\sS_4$ and cannot fix the point~$v$. Finally, let 
$\Gamma{\,\cong\,}\PU3(\HH,1)$. Then 
$\Gamma$ has $3$ pairwise commuting involutions $\alpha,\beta$, and $\gamma{\,=\,}\alpha\beta$ such that  $\alpha$ and $\beta$ are conjugate and hence are planar.  The centralizer  
$\Ypsilon{\,=\,}\Cs{}\beta$ is locally isomorphic to  $\U2(\HH,1){\times}\HH'$, and 
$\dim\Ypsilon|_{\cF_\beta}{\,\le\,}11$ by 
\cite{sz1} 4.1 or 7.3. Therefore the compact factor of 
$\Ypsilon$ acts trivially on $\cF_\beta$, and $\gamma$ is also planar. The group 
$\U2\HH{\,<\,}\Cs{}\gamma$ acts faithfully on $\cF_\gamma$; for the same reason as before, 
it cannot fix a flag. \\
($\ell$) Lastly, let $\dim\Gamma{\,\ge\,}24$.
The compact group ${\rm(P\hskip-.8pt)}\SU5\CC$ has two conjugacy classes of pairwise 
commuting  (diagonal) involutions consisting of $5$ and of $10$ elements respectively, but this contradicts 2.10. If $r{\,>\,}0$, then some diagonal planar involution $\iota$ has a centralizer 
$\Upsilon{\,\cong\,}\SU4(\CC,r{-}1)$; on the other hand, $\dim\Upsilon|_{\cF_\iota}{\,\le\,}11$ by 
\cite{sz1} 7.3. This contradiction shows that $\Gamma$ is not of type   ${\rm A}_4$. \\
(m) {\tt Conclusion.} In steps (h--$\ell$), the following has been shown: If $\dim\Delta{\,>\,}20$ and if 
$\Gamma$ is an almost simple factor of $\Delta$, then $\Gamma$ does not contain an involution, 
and $\dim\Delta{\,=\,}21$; if $\dim\Delta{\,\le\,}20$ and $\Delta$ is almost simple, then 
$\dim\Delta{\,\ne\,}20$ and hence $\dim\Delta{\,\le\,}16$. \qed
\par\medskip
{\bf 5.3 Normal torus.} {\it If $\cF_\Delta{\,=\,}\{v,W\}$ is a flag, if $\Delta$ is a Lie group with a normal  subgroup $\Theta{\,\cong\,}\TT$, and if $\dim\Delta{\,\ge\,}18$, then there exists a 
$\Delta$-invariant Baer subplane which is a known translation plane\/}.
\par\smallskip
{\tt Proof.} Note that $\Theta{\,\le\,}\Cs{}\Delta$\, (see \cite{cp} 93.19). The involution 
$\iota{\,\in\,}\Theta$ cannot be a reflection, and $\cF_{\hskip-1.5pt\iota}{\,\ldot}\cP$. Write
$\Delta^{\hskip-2pt\ast}{\,=\,}\Delta|_{\cF_{\hskip-1pt\iota}}{\,=\,}\Delta/\Kappa$. 
Obviously, $\Delta^{\hskip-2pt\ast}$ fixes exactly $v$ and $W$. If $\Theta$ acts non-trivially on
$\cF_\iota$, then $\Delta{\,:\,}\Kappa{\,\le\,}11$ by  
\cite{sz1} 4.2,  Stiffness gives 
$\dim\Kappa{\,\le\,}3$, and $\dim\Delta$ would be too small. Hence $\Theta{\,\le\,}\Kappa$,\: 
$\Kappa^1{\,\not\hskip1pt\cong\,}\Spin3\RR$,\: $\dim\Kappa{\,=\,}1$, and 
$\dim\Delta^{\hskip-2pt\ast}{\,\ge\,}17$. 
As $v$ and $W$ are fixed, $\cF_\iota$ is not a Hughes pane, and the claim follows from 
\cite{sz1} 1.10.  \qed
\par\medskip
Recall from 2.3 that we put $\Sigma{\,=\,}\Aut\cP$. The following has been proved in L\"uneburg's dissertation \cite{Lb}:
\par\smallskip
{\bf 5.4  Large groups.} {\it Assume that $\dim\Delta{\,\ge\,}38$, that $\cF_\Delta{\,=\,}\{v,W\}$ is
 a flag, and that $\Ypsilon$ is a maximal semi-simple subgroup of $\Delta$. Then there are $3$ possibilities\/}:\par
\quad(a) {\it $\cP$ or its dual is a translation plane\/}, \par
\quad(b) {\it  $\Ypsilon{\,\cong\,}\Spin7\RR$,\: $\dim\Sigma{\,=\,}38$, and 
$\Delta_{[v,W]}{\,\cong\,}\RR^8$\/} \ (\cite{Lb} V Satz), \par
\quad(c) {\it $\Ypsilon{\,\cong\,}\Gtwo$ is a maximal compact subgroup of $\Delta$\/}; 
{\it if  $\dim\Delta{\,=\,}39$, then  $\Delta_{[v,W]}{\,\cong\,}\RR^8$ and,\par  \hskip30pt  up to duality, 
$\Delta$ is transitive on $\frL\sm\frL_v$ and hence on $W\sm\{v\}$\/} \par \hskip30pt  
(\cite{Lb} VI Kor.\:5, Satz\:5, and  Kor.\:3).
\par\smallskip
{\tt Remarks.} All {\it translation\/} planes with $\dim\Sigma{\,\ge\,}38$ have been determined 
by H\"ahl,  see \cite{cp} 82.\hskip1pt28,\hskip1pt29 and \cite{Ha6}, \cite{Ha3}.  
No substantial progress has been made in the meantime.
\par\smallskip
A few additional results can be obtained in case (c) of the preceding theorem:
\par\medskip
{\bf 5.5 Proposition.}  {\it  Suppose that  $\dim\Delta{\,=\,}39$,  that 
$\cF_\Delta{\,=\,}\{v,W\}$ is a flag, and that  neither $\cP$ nor its dual is a translation plane.  Let 
$\Ypsilon$ denote  a Levi complement  of $\sqrt\Delta$, choose  $a{\,\notin\,}W$, and
put $\Gamma{\hskip1pt=\,}(\Delta_a)^1$. Then\/} \par
\quad(a) {\it $\cF_\Ypsilon$ is isomorphic to the classical real plane $\cP_\RR$\/}, \par
\quad(b) {\it if $a^{\hskip-1pt\Ypsilon\hskip-2.5pt}{\,=\,}a$, then $\dim a^\Delta{\,\ge\,}15$\/}, \par
\quad(c) {\it $\Delta$ is transitive on $P\sm W$ and, dually, on $\frL\sm\frL_v$, or 
 $\Gamma_{\hskip-1pt[u]}{\,=\,}\1$  for each  $u{\,\in\,}W\sm\{v\}$\/}; \par
\hskip30pt {\it if $\Delta_{[a]}{\,\ne\,}\1$, then $\dim a^\Delta{\,=\,}15$ and $\Delta_{[a]}$ does not contain an involution\/},  \par
\quad(d) {\it if $a^{\hskip-1pt\Ypsilon\hskip-2.5pt}{\,=\,}a$, then $a^\Delta$ is open in $P$ or 
$\dim\Delta_{[v,av]}{\,=\,}1$\/}.
\par\smallskip
{\tt Proof.} (a) The representation of $\Ypsilon{\,\cong\,}\Gtwo$ on the Lie algebra 
$\frl\hskip1pt\sqrt\Delta$ together with \cite{cp} 95.10 shows that 
$\dim\Cs\Delta{\hskip-2pt\Ypsilon}{\,=\,}4$. By \cite{cp} 96.35 the fixed points of $\Ypsilon$ on $W$ 
form a circle $S$. Dually, $\Ypsilon$ fixes a one-dimensional set of lines in each pencil $\frL_s$ with 
$s{\,\in\,}S$. Therefore $\cF_\Ypsilon$ is a flat subplane. Stiffness implies that the connected component of $\Cs{}\Ypsilon$ acts effectively on $\cF_\Ypsilon$. According to \cite{cp} 33.9(a), 
either $\cF_\Ypsilon$ is classical, or $\Cs{}\Ypsilon$ fixes exactly a non-incident point-line pair. \\
(b)  Note that $\Ypsilon$ is a Levi complement of $\Rho{\,=\,}\sqrt\Gamma$. Choose 
$b{\,\notin\,}W\hskip-2pt\smcup av$, and consider a minimal $\Gamma_{\hskip-1pt b}$-invariant subgroup $\Xi{\,\cong\,}\RR^s$ of $\Delta_{[v,W]}{\,\cong\,}\RR^8$. As $\Gamma$ does not 
act irreducibly on $\Delta_{[v,W]}$, we have $s{\,<\,}8$. For 
$c{\,\in\,}a^\Xi\sm\{a\}$, the action of $\Lambda{\,=\,}(\Gamma_{\hskip-1pt b,c})^1$ on $\Xi$ shows 
$\Lambda{\,\not\cong\,}\Gtwo$. The stiffness result 2.6(e) implies $\dim\Lambda{\,\le\,}8$ and 
$\dim\Gamma_{\hskip-1pt b}{\,\le\,}15$. Now let $\Delta{\,:\,}\Gamma{\,=\,}\dim a^\Delta{\,=\,}k$, 
and choose $b{\,\in\,}a^\Delta$. Then 
$39{\,=\,}\dim\Delta{\,\le\,}2k{\,+\,}\dim\Gamma_{\hskip-1pt b}{\,\le\,}2k{\,+\,}15$. Consequently  
$k{\,\ge\,}12$ and $\dim\Gamma{\,\le\,}27$. We have $\Gamma{\,=\,}\Ypsilon\Rho$ with 
$\dim\Rho{\,<\,}14$. The representation of $\Ypsilon$ on $\frl\hskip1pt\Rho$ shows 
$\dim\Rho{\,=\,}7{\,+\,}\dim(\Rho\smcap\Cs{}\Ypsilon){\,\le\,}7{\,+\,}3$ and 
$\dim\Gamma{\,\le\,}24$. \\ 
(c) Suppose that $\dim a^\Delta{\,=\,}15$ and that $\Gamma_{\hskip-1pt[u]}{\,\ne\,}\1$.
Put $L{\,=\,}av$ and $\Epsilon{\,=\,}\Delta_{[v,L]}$. 
By \cite{cp}  61.20, we have  $e{\,=\,}\dim\Epsilon{\,=\,}\dim u^{\Delta_L}$ and 
$h{\,=\,}\dim\Tau_{\hskip-1pt[u]}{\,=\,}\dim L^{\Delta_u}$. 
We may assume that  $\Upsilon{\,\le\,}\nabla{\,:=\,}\Gamma_{\hskip-1.5pt u}$, since
$u^\Delta{\,=\,}W\sm\{v\}$ according to  5.4(c). It will turn out  that  $e{\,>\,}0$. 
If $\Delta_L{\,\le\,}\Delta_u$ and if $u^\delta{\,\ne\,}u$, then $L^\delta{\,\ne\,}L$, and 
$\Gamma_{\hskip-2pt L^\delta}$ fixes $a, u, L^\delta, u^\delta$ and, hence, a quadrangle, but 
$\dim\Gamma{\,=\,}24$ and $\dim\Gamma_{\hskip-2pt L^\delta}{\,>\,}14$, a contradiction to Stiffness.
Therefore $u^{\Delta_L}{\,\ne\,}u$. Analogously, $L^{\Delta_u}{\,\ne\,}L$ and $h{\,>\,}0$. By assumption $e,h{\,<\,}8$, and the action of $\Ypsilon$ implies $e,h{\,\in\,}\{1,7\}$.  Moreover, 
$\dim u^\Gamma{\,\le\,}e{\,=\,}\dim u^\Epsilon{\,\le\,}\dim u^{\Delta_{[L]}}{\,\le\,}\dim u^\Gamma$. 
Therefore $24{\,=\,}\dim\Gamma{\,=\,}\dim\nabla{\,+\,}e$ and   $\dim\nabla{\,\ge\,}17$.
On the other hand, $\dim\nabla{\,\le\,}e{\,+\,}h{\,+\,}14$ and $10{\,\le\,}2e{\,+\,}h$. In the same way 
it follows that $10{\,\le\,}e{\,+\,}2h$. Therefore $e{\,=\,}h{\,=\,}7$ and  $\dim\nabla{\,=\,}17$.  
Note that $\Ypsilon$ and $\nabla$ act irreducibly on $\Epsilon$ and on $\Tau_{\hskip-1pt[u]}$. 
The action of $\Ypsilon$ on $\lie\nabla$ shows that  $\Mu{\,=\,}\nabla\smcap\Cs{}\Upsilon$ 
satisfies $\dim\Mu{\,=\,}3$. On the other hand, $\Mu$  induces  on $\Epsilon$ and on 
$\Tau_{\hskip-1pt[u]}$ groups of dimension at most $1$,\: see \cite{cp} 95.10. Hence $\Mu$ 
has a non-trivial subgroup which acts trivially on $\langle a^{\Tau_{[u]}},L^\Epsilon\rangle{\,=\,}\cP$.
This contradiction proves  the first part of the claim. The second part follows from \cite{cp} 61.20 and 2.9 above. \\
(d) If $a^\Delta$ is not open, then $\dim\Gamma{\,=\,}24$. By step (c), $\Gamma$ does not contain 
a homology with axis $av$. Representation of $\Ypsilon$ on the Lie algebra $\lie\Gamma$ shows 
that $\dim\Cs\Gamma{\hskip-2pt\Ypsilon}{\,\ge\,}3$.  Put again $\Epsilon{\,=\,}\Delta_{[v,av]}$ and 
$e{\,=\,}\dim\Epsilon$.  As has been stated in 2.15(a), $\Ypsilon$ acts on 
$\Xi{\,=\,}\Delta_{[v,W]}{\,\cong\,}\RR^8$ in the same way  as $\Aut\OO$. In this action the centralizer 
of $\Ypsilon$ in $\Aut\Xi$ is $2$-dimensional. We have $\Gamma\smcap\Cs{}\Xi{\,\le\,}\Epsilon$. 
If $e{\,=\,}0$, then $\Epsilon{\,=\,}\1$, and $\Gamma$embeds into $\Aut\Xi$, so that  
$\dim\Cs\Gamma{\hskip-2pt\Ypsilon}{\,\le\,}2$. Hence  $e{\,>\,}0$.
As $\dim\{z{\,\in\,}W\mid z^{\hskip-1.5pt\Ypsilon}{\,=\,}z\}{\,=\,}1$, it follows that $e{\,\in\,}\{1,7\}$.
If $e{\,=\,}7$, then $\Ypsilon$ acts irreducibly on $\Epsilon$ and 
$\Epsilon\smcap\Cs\Gamma{\hskip-2pt\Ypsilon}{\,=\,}\1$. Again this would imply 
$\dim\Cs\Gamma{\hskip-2pt\Ypsilon}{\,\le\,}2$.    \qed
\par\medskip
{\bf Conjecture.} {\it If $\dim\Delta{\,=\,}39$ and if $\cF_\Delta$ is a flag, then $\cP$ or its dual is a 
translation plane\/}.
\par\bigskip
{\Bf B. Non-incident fixed elements} 
\par\medskip
The following has been proved in \cite{sz12}; for semi-simple groups see also  4.2:  
\par\smallskip
{\bf 5.6 Theorem.} {\it If  $\Delta$ fixes exactly one line $W$ and one point $a{\,\notin\,}W$, and if 
$\dim\Delta{\,\ge\,}35$, then $\Delta$ contains a group $\Spin9(\RR,r)$ with 
$r{\,\le\,}1$ and $\dim\Delta{\,\le\,}37$, or $\Delta$ is triply transitive on~$W$ and $\cP$ is the classical Moufang plane\/}.
\par\medskip
{\bf 5.7 Almost simple groups.} {\it If $28{\,\le\,}\dim\Delta{\,<\,}36$, 
$\,\cF_\Delta{\,=\,}\{a,W\},\: a{\,\notin\,}W$, and if $\Delta$ is almost simple, then $\Delta$ 
is isomorphic to one of the groups $\Spin8(\RR,r),\, r{\,=\,}1,3$ or to some covering of the anti-unitary group  {\cyss Ya}${=\,}\SaU4{\,\circeq\,}\Opr8(\RR,2)$ and $\dim\Delta{\,=\,}28$\/}.
\par\smallskip
{\tt Proof.} Because of Priwitzer's results 4.2, it suffices to consider groups of dimension~$28$. From 2.16 it follows that $\Delta{\,\not\hskip1pt\cong\,}\Gtwo^{\CC}$. Hence $\Delta$ is locally isomorphic to an orthogonal group $\Opr8(\RR,r)$. Recall from 2.10 and 2.11 that $\Delta$ has no subgroup $\SO5\RR$ or $(\SO3\RR)^2$. This excludes all groups ${\rm(P)}\Opr8(\RR,r)$. Each 
central involution $\zeta{\,\in\,}\Delta$ is a reflection (or $\dim\Delta|_{\cF_\zeta}{\,=\,}28$ and 
$a^\Delta{\,\ne\,}a$). If $r$ is even, then  $\Spin8(\RR,r)$ has a center $(\ZZ_2)^2$ and hence 
fixes a triangle; in particular, $r{\,\ne\,}0$. Note that the double covering  $\Spin8(\RR,r)$ of  
$\Opr8(\RR,r)$ is not simply 
connected for $r{\,\ge\,}2$.  In each  case, let $\tilde\Delta$ denote the simply connected covering group of $\Delta$. If $r{\,\ge\,}2$, the center of $\Delta$ consists of homologies with axis 
$W\,$\:  (or some stabilizer $\Delta_p$ would fix a quadrangle and $(\Delta_p)^1{\,\cong\,}\Gtwo$ 
by Stiffness,  but then $r{\,\le\,}1$). \\
(a)  First let $r{\,=\,}2$. The group  $\Opr8(\RR,2){\,\circeq}${\cyss Ya} has center 
{\rm Cs}\,{\cyss Ya}${\,\cong\,}  \ZZ_2$ by \cite{cp} 94.32(f),\hskip1pt33, and 95.10. A maximal 
compact subgroup $\Kappa$ of {\cyss Ya} is isomorphic to $\U4\CC$ and satisfies 
$\Kappa'{\,\cong\,}\Spin6\RR$. Thus {\cyss Ya} and its coverings are possible candidates for a group of automorphisms of a compact plane. \\
(b) If $r{\,=\,}3$, then a maximal compact subgroup of $\tilde\Delta$ is isomorphic to  
$\Spin5\RR{\times}\Spin3\RR$, and $\Cs{}{\tilde\Delta}{\,\cong\,}(\ZZ_2)^2$. Hence  $\Delta$ is not simply connected. Assume now that $\Delta$ has a compact subgroup 
{\cyss Yu}${\times}\Omega{\,=\,}\Spin5\RR{\times}\SO3\RR$, and choose commuting involutions 
$\alpha,\beta{\,\in\,}\Omega$, they are planar (as in the proof of 2.11).  By \cite{cp} 96.13, the factor {\cyss Yu} acts transitively on $S{\,=\,}W{\smcap}\cF_\beta$, and the central involution $\sigma$ of 
{\cyss Yu} induces a reflection on $\cF_\beta$, moreover, $\alpha|_S{\,=\,}\1$ because $\alpha$ fixes some  point $z{\,\in\,}S$ and {\cyss Yu} is transitive on $S$.
According to \cite{cp} 55.32, commuting  involutions with the same fixed point set are equal. 
Therefore  $\alpha|_{\cF_\beta}{\,=\,}\sigma|_{\cF_\beta}$, $\,\alpha\sigma|_{\cF_\beta}{\,=\,}\1$, 
$\,\alpha\sigma{\,=\,}\beta$, and $\sigma{\,=\,}\alpha\beta{\,\in\,}\Omega$, a contradiction. \\
(c) For $r{\,=\,}4$, a maximal compact subgroup of $\tilde\Delta$ is isomorphic to 
$(\Spin4\RR)^2{\,\cong\,}(\Spin3\RR)^4$, and $\Cs{}{\tilde\Delta}{\,\cong\,}(\ZZ_2)^4$.  
Consequently $\Delta$ is simple or a double covering of  ${\rm(P)}\Opr8(\RR,4)$, but then $\Delta$ has a subgroup $(\SO3\RR)^2$. \qed
\par\medskip
{\tt Remarks.} Suppose that $\Delta$ satisfies the conditions of 5.7 and that $\Delta$ acts on the classical octonion plane $\cO$. Then $\Delta$ is contained in the affine group 
$\Alpha{\,=\,}\Sigma_{a,W}$ of $\cO$. From \cite{cp} 15.6 it follows that 
$\Alpha'{\,\cong\,}\Spin{10}(\RR,1)$, and $\dim\Alpha'{\,=\,}45$.  Hence $\Delta$ intersects a 
maximal compact subgroup $\Spin9\RR$ of $\Alpha'$ in a group of dimension at least $19$. 
This excludes all orthogonal groups of index $r{\,>\,}1$\: (see\cite{cp} 94.33). 
Consequently $\Delta{\,\cong\,}\Spin8(\RR,1)$. For each $z{\,\in\,}W$ a Levi complement of 
$\sqrt{\Alpha_z}$ in $\Alpha$ is a compact group $\Spin8\RR$. Hence $z^\Delta{\,\ne\,}z$, and 
$\Delta$ has no fixed point other than $a$. No other example for the situation of 5.7 has been 
found as yet.
\par\medskip
{\bf 5.8 Semi-simple groups.} {\it If $\dim\Delta{\,=\,}28$,\: $\,\cF_\Delta{\,=\,}\{a,W\}$ with 
$a{\,\notin\,}W$, and if $\Delta$ is semi-simple,  then   $\Delta$ is almost simple\/}.
\par\smallskip
{\tt Proof.} Suppose first that $\Delta$ has exactly two almost simple factors. The list \cite{cp} 94.33 
shows that there are only the following two cases: \\
(a) $\:\Delta$ has a factor $\Psi{\,\cong\,}\Sp4\CC$ with a maximal compact subgroup 
$\Phi{\,\cong\,}\U2\HH$  and an $8$-dimensional factor $\Gamma$ of type ${\rm A}_2$. The central involution $\sigma{\,\in\,}\Psi$ is a reflection in $\Delta_{[a,W]}$ because 
$\Psi/\langle \sigma\rangle{\,\cong\,}\SO5\CC$ cannot  act on a Baer subplane. Let   
$\alpha{\,=\,}\bigl({1 \atop }{ \atop -1}\bigr){\,\in\,}\Phi$. 
The centralizer $\Cs{}\alpha$ contains a group $\Gamma{\cdot\hskip1pt}(\Sp2\CC)^2$ of dimension $20$. 
Stiffness and \cite{cp} 83.26 imply that $\alpha$ is not planar. Therefore $\alpha$ is a reflection, 
its center is some point $u{\,\in\,}W$, and $u^\Gamma{\,=\,}u{\,\ne\,}u^\Psi$ by the assumption on 
$\cF_\Delta$. In fact, $u^\Phi{\,\ne\,}u$ because $\alpha$ is conjugate to $\alpha\sigma$ in $\Phi$ 
and the center of $\alpha\sigma$ is incident with the axis of $\alpha$. Now $\dim u^\Phi{\,\ge\,}4$ 
by \cite{cp} 96.13. Note that $\Gamma|_{u^\Psi}{\,=\,}\1$. From 2.17(b) it follows that $\Gamma$ is 
not transitive on  $H{\,=\,}au\sm\{a,u\}$, and there is some $b{\,\in\,}H$ such that 
$\Gamma_{\hskip-1.5pt b}{\,\ne\,}\1$. Hence 
$\cF_{\Gamma_{\hskip-,5pt b}}{\,= \langle a,b,u^\Psi\rangle\,}{\,=\,}\cB$ is a Baer subplane of $\cP$, 
and $\cB^\Psi{\,=\,}\cB$ since $\Gamma{\,\le\,}\Cs{}\Psi$. Moreover, $\Psi$ acts effectively on $\cB$ 
and fixes $a$ and $W$, but $\dim(\Aut\cB)_{a,W}{\,\le\,}19$, see again \cite{cp} 83.26. Thus,
case (a) is impossible. \\
(b) $\:\Delta$ is a product of two $14$-dimensional factors of type ${\rm G}_2$. The compact form is excluded by Stiffness. A maximal compact subgroup of ${\rm G}{\,=\,}\Gtwo(2)$ is isomorphic to 
$\SO4\RR$. Lemma 2.11 implies that at least one of the two factors of $\Delta$ is the simply connected covering group  $\Gamma$ of ${\rm G}$, the other factor $\Psi$ is locally isomorphic to 
${\rm G}$.  Stroppel's theorem \cite{str} 4.5\, (cf. 
\cite{sz1} 2.2) shows that neither $\Gamma$ nor $\Psi$ contains a planar involution.  As both factors have torus rank $2$, the group $\Delta$ contains a group  $\ZZ_2^{\hskip2pt4}$ generated by reflections. This contradicts 2.10 and proves that $\Delta$ has more than two factors. \\ 
If $\Delta$ is properly semi-simple, only the following possibilities remain: \\
(c) $\:\Delta$ has a factor $\Psi{\,\circeq\,}\SL3\CC$ and some factors of dimension $3$ or $6$. 
There are $3$ pairwise commuting conjugate involutions in $\Psi$ corresponding to diagonal 
elements in $\SL3\CC$. If they were  reflections, their centers would be $\Psi$-conjugate, but one would have center~$a$, which is fixed. Hence there is a planar involution $\beta{\,\in\,}\Psi$
such that $\SL2\CC{\,\le\,}\Cs\Psi\beta$, and $\Cs\Delta\beta$ contains a semi-simple 
$18$-dimensional group  inducing on $\cF_\beta$ a semi-simple group $\Chi$ of 
dimension at least $15$\: (by Stiffness). Therefore, $\cF_\beta$ is isomorphic to the classical quaternion plane, 
see \cite{sz9} 3.3 or \cite{sz1} 4.4. The group $\Chi$ is contained in a maximal semi-simple 
subgroup $\Beta$ of $(\Aut\cF_\beta)_{a,W}$. We have $\dim\Beta{\,=\,}18$ and $\Beta$ has 
a $10$-dimensional compact subgroup \Yu${\,\cong\,}\U2\HH$. It follows that 
$\Phi{\,=\,}\Chi\smcap$\Yu\ has dimension $\ge\!7$. If $\Phi{\,<\,}$\Yu, then 
\Yu\  acts trans\-itively and almost effectively on the coset space \Yu$/\Phi$, and \cite{cp} 96.13(a) 
would imply \Yu${\hskip-1pt:\hskip0pt}\Phi{\,\ge\,}4$. Consequently, $\Chi$ has a subgroup 
$\U2\HH$, but this is impossible because $\Chi$ is a product of almost simple factors of dimension $3$ or $6$.   \\
(d) $\:\Delta$ has a factor $\Psi$ of type ${\rm A}_3$, a factor $\Ypsilon$ of type ${\rm C}_2$, 
and a factor $\Omega$ of dimension~$3$.  Several combinations are excluded by 2.10 and 11; 
the rank condition, in particular, implies that at least one of the factors has an infinite center, 
cf. \cite{cp} 94.32(e) and 33.  In any case,  $\rk\Psi{\,\ge\,}2$ and $\Psi$ contains a non-central involution $\alpha$ centralizing a $6$-dimensional semi-simple subgroup $\Gamma$ of $\Psi$. 
If $\alpha$ is a reflection, its center is some point $u{\,\in\,}W$\: (or else 
$\langle \alpha^\Psi\rangle{\,=\,}\Psi$ would act trivially on $W$), and 
$u^{\Ypsilon\Omega}{\,=\,}u{\,\ne\,}u^\Psi$. As $\Psi$ acts almost effectively on $u^\Psi$ and 
contains a compact subgroup of dimension ${\!\ge\,}6$, it follows from \cite{cp} 96.13 that 
$\dim u^\Psi{>\,}2$. Let $b{\,\in\,}au\sm\{a,u\},\ \Lambda{\,=\,}(\Ypsilon\Omega)_b$, and 
$\cF{\,=\,}\langle u^\Psi,a,b\rangle$. Then $\dim\Lambda{\,\ge\,}5$,\: $\Lambda|_\cF{\,=\,}\1$, 
and $\cF{\,\ldot}\cP$, which contradicts Stiffness. If $\alpha$ is planar, then 
$\Gamma\Ypsilon\Omega$ induces on $\cF_\alpha$ a semi-simple group $\Chi$ of dimension
 at least 16. Again $\cF_\alpha$ is the classical quaternion plane (cf. 
\cite{sz1} 4.4) and $\Chi$ is 
 contained in $\SL2\HH{\,\cdot}\Aut\HH$. As $\SL2\HH$ has no proper subgroup of codimension 
${\le\,}2$, it follows that $\SL2\HH{\,\le\,}\Chi$, but $\Chi$ is a product of factors of smaller dimension. Hence case (d) is impossible.   \\
(e) $\:\Delta$ has a factor $\Psi$ of type ${\rm G}_2$ and a factor $\Ypsilon$ of type ${\rm A}_2$. Stiffness implies that $\Psi$ is not compact. A maximal compact  subgroup $\Phi$ of $\Psi$ is isomorphic to $\SO4\RR$ or to $\Spin4\RR$, and there is an involution $\alpha{\,\in\,}\Phi$ such that 
$\Phi{\,\le\,}\Cs\Psi\alpha{\,<\,}\Psi$. As in the previous step, $\alpha$ is not planar, and   
$\alpha|_W{\,\ne\,}\1$ because $\langle \alpha^\Psi\rangle{\,=\,}\Psi{\,\not\le\,}\Delta_{[W]}$. 
Therefore $\alpha$ has some center $u{\,\in\,}W$. Choose $b{\,\in\,}au\sm\{a,u\}$ and let 
$\Lambda{\,=\,}(\Cs{}\Psi)_{\hskip-.5pt b}^{\,1}$. Then $\dim\Lambda{\,\ge\,}6$ and 
$\dim\cF_\Lambda{\,\le\,}4$, but then $\Psi$ cannot act on $\cF_\Lambda$. This excludes case (e). \\
(f) $\:\Delta$ has a factor $\Ypsilon$ of type ${\rm C}_2$,  all other factors have  dimension at most 
$10$. Put $\Psi{\,=\,}(\Cs\Delta{\hskip-2pt\Upsilon})^1$. Then $\dim\Psi{\,=\,}18$. If $\Psi$ fixes a point or a line other than $a$ or $W$, then there is some point $u{\,\in\,}W$ such that 
$u^\Psi{\,=\,}u{\,\ne\,}u^\Ypsilon$ and $\Psi|_{u^\Ypsilon}{\,=\,}\1$. If $b{\,\in\,}au\sm\{a,u\}$, 
we have $\dim\Psi_b{\,\ge\,}10$;   the fixed elements of  $\Psi_b$ form a  subplane $\cE{\,<\,}\cP$, 
and  $\Psi_b|_\cE{\,=\,}\1$. Hence $\cE$ is flat by Stiffness, $\cE^\Ypsilon{\,=\,}\cE$, 
and $\dim(\Ypsilon|_\cE){\,=\,}10$. This contradiction shows that that $\langle x^\Psi\rangle$ is a subplane  whenever $a{\,\ne\,}x{\,\notin\,}W$. The r\^oles of $\Ypsilon$ and $\Psi$ can be 
interchanged, and $u^{\hskip-1pt\Ypsilon}{\,=\,}u$ implies successively  $u^\Psi{\,\ne\,}u$,\: 
$\Ypsilon|_{u^\Psi}{\,=\,}\1$,\: $\dim\Ypsilon_{\hskip-1pt b}{\,\ge\,}2$,\: 
$\cF{\,=\,}\cF_{\Ypsilon_{\hskip-1pt b}}{\,<\,}\cP$,\: $\Ypsilon_{\hskip-1pt b}|_\cF{\,=\,}\1$, and 
$\cF^\Psi{\,=\,}\cF$. If $\cF$ is flat, then $\dim\Psi{\,\le\,}4{+}11$\:  (use Stiffness and 
$\Ypsilon_{\hskip-1pt b}|_\cF{\,=\,}\1$). Similarly, $\dim\cF{\,=\,}4$ implies $\dim\Psi{\,\le\,}8{+}8$.
Therefore $\cF{\,\ldot}\cP$, and $\dim\Psi|_\cF{\,=\,}18$\: (again because 
$\Ypsilon_{\hskip-1pt b}|_\cF{\,=\,}\1$). 
Now $\cF$ is classical by \cite{sz9} 3.3 or 
\cite{sz1} 4.4, and $\Psi|_\cF$ contains a group $\SL2\HH$ contrary to the assumption in case (f). It follows that $\langle x^\Ypsilon\rangle$ and 
$\langle x^\Psi\rangle$ are subplanes of $\cP$, and $\langle x^\Ypsilon\rangle{\,<\,}\cP$, since 
$\Psi_{\hskip-1pt x}|_{\langle x^\Ypsilon\rangle}{\,=\,}\1$ and $\dim\Psi_{\hskip-1pt x}{\,\ge\,}2$. Consequently, 
$\dim x^{\hskip-1pt\Ypsilon}{\,\le\,}8$ and $\dim\Ypsilon_{\hskip-1.5pt x}{\,\ge\,}2$. 
Hence $\langle x^\Psi\rangle$ is also a proper subplane, $\dim x^\Psi{\,\le\,}8$, and 
$\dim\Psi_{\hskip-1pt x}{\,\ge\,}10$. Stiffness implies that $\cE{\,=\,}\langle x^\Ypsilon\rangle$ is even flat, but
$\dim(\Aut\cE)_{a,W}{\,\le\,}4$.  \\
(g) Last case: all factors of $\Delta$ have dimension at most $8$. Then there are exactly two $8$-dimensional factors $\Ypsilon_{\hskip-1.5pt1},\!\Ypsilon_{\hskip-1.5pt2}$. An involution 
$\sigma{\,\in\,}\Cs{}\Delta$ is a reflection with axis $W$: if $\cF_\sigma{\,\ldot     }\cP$, then 
$\dim\Delta|_{\cF_\sigma}{\,<\,}20$ contrary to Stiffness. Moreover, each involution $\iota$ in one 
of the factors is a reflection: if $\iota$ is planar, denote the product of the other factors by $\Psi$ and 
put $\Psi|_{\cF_\iota}{\,=\,}\Psi/\Kappa$. Then $\Psi{\hskip0pt:\hskip1pt}\Kappa{\,<\,}20$,\: 
$\dim\Kappa{\,>\,}0$, and $\Kappa^1$ is semi-simple. Stiffness implies $\Kappa^1{\,\cong\,}\Spin3\RR$, and the involution  $\sigma{\,\in\,}\Kappa^1$ acts trivially on $\cF_\iota$; on the other hand, $\sigma$ is in the center of $\Delta$, but then $\sigma{\,\in\,}\Delta_{[a,W]}$. 
It follows that none of the factors has a subgroup isomorphic to 
$\SO3\RR$. Suppose that $\Delta$ has a factor $\Ypsilon{\circeq\,}\SU3(\CC,r)$. Then $\Ypsilon$ 
has a subgroup $\Phi{\,\cong\,}\SU2\CC$. The involution $\alpha{\,\in\,}\Phi$ is a reflection. As 
$\langle \alpha^\Ypsilon\rangle{\,=\,}\Ypsilon$, the axis of $\alpha$ is a line $av$, and 
$\dim v^\Ypsilon{\,>\,}1$. Choose $c{\,\in\,}av\sm\{a,v\}$ and put $\Lambda{\,=\,}(\Psi_c)^1$. Then 
$\cF_\Lambda{\,=\,}\langle a,c,v^\Ypsilon\rangle$ and $\dim\Lambda{\,\ge\,}12$. This contradicts 
Stiffness and proves that the $\Ypsilon_{\hskip-1.5pt\nu}$ are isomorphic to the simply connected covering group  of $\SL3\RR$. This possibility will be excluded in the next steps. \\
(h) Assume that $\Delta$ has a factor $\Phi{\,\cong\,}\SU2\CC$. Then $\Ypsilon_{\hskip-1.5pt1}$ 
and  $\Phi$ have the reflection $\sigma{\,\in\,}\Delta_{[a,W]}$ in common. Consequently, 
$\Ypsilon_{\hskip-1.5pt1}\Phi$ contains a subgroup $\Omega{\,\cong\,}\SO4\RR{\,>\,}\SO3\RR$.
Each involution $\beta{\,\in\,}\SO3\RR$ is planar. Put $\Psi{\,=\,}(\Cs{}\beta)^1$, note that 
$\Omega\smcap\Psi{\,\cong\,}\TT^2$, and consider 
$\overline\Psi{\,=\,}\Psi|_{\cF_\beta}{\,=\,}\Psi/\Kappa$. We have 
$\dim\Psi{\,\ge\,}2{+}17$,\: $\dim\overline\Psi{\,\ge\,}16$, and $\overline\Psi$ has a normal 
torus subgroup. According to 
\cite{sz1} 7.3,\: $\cF_\beta$ is the classical quaternion plane,
$\overline\Psi$ is contained in $\HH^{\times}{\cdot\,}\SL2\HH$, and 
$\dim(\overline\Psi\smcap\SL2\HH){\,\ge\,}12$, but such a group $\overline\Psi$ does not exist. 
Therefore $\Delta$ has no compact factor. If $\SL2\CC{\,\cong\,}\Gamma{\,\triangleleft\,}\Delta$, 
analogous arguments yield a group $\overline\Psi$ of dimension at least $13$. Again $\cF_\beta$ is 
classical, and $\Ypsilon_{\hskip-1pt2}$ maps onto a subgroup $\SL3\RR$ of $\SL2\HH$, since 
$\Ypsilon_{\hskip-1pt2}$ has no faithful linear representation, cf. \cite{cp} 95.9. In particular, 
$\sigma|_{\cF_\beta}{\,=\,}\1$, which is absurd. \\ 
(i) Steps (g,h) imply that there are $4$ factors $\Gamma_{\hskip-1.5pt\nu}{\,\circeq\,}\SL2\RR$. 
For similar reasons as before,   $\Ypsilon_{\hskip-1pt1}\Ypsilon_{\hskip-1pt2}$ has a subgroup 
$\Omega{\,\cong\,}\SO4\RR$ containing a planar involution $\beta$. Let again 
$\Psi{\,=\,}(\Cs{}\beta)^1$ ${\,\cong\,}\TT^2{\cdot}${\Large{$\!\times\!$}}$_\nu\Gamma_\nu$ and 
$\overline\Psi{\,=\,}\Psi|_{\cF_\beta}$. As before, $\cF_\beta$ is classical and $\overline\Psi$ is a subgroup of $\HH^{\times}{\cdot\,}\SL2\HH$. None of the groups $\overline\Gamma_{\hskip-1pt\nu}$ 
is simply connected, each $\Gamma_{\hskip-1pt\nu}$ maps onto a group of rank $1$. Hence the rank of $\overline\Psi$ is at least $5$, but this contradicts 2.10 or 
$\Psi{\,\le\,}\HH^{\times}{\cdot\,}\SL2\HH$. \qed
\par\medskip
{\bf 5.9 Normal torus.} {\it If $\dim\Delta{\,>\,}28$ and $\,\cF_\Delta{\,=\,}\{a,W\},\: a{\,\notin\,}W$, 
 if $\Delta$ has a normal subgroup $\Theta{\,\cong\,}\TT$, and if $\Delta$ is doubly transitive on 
 some proper subset $V{\,\subset\,}W$, then $V{\,\approx\,}\sS_6$,\, $\Delta|_V{\,\cong\,}\POpr8(\RR,1)$  is even triply transitive, and  $\dim\Delta{\,\le\,}30$\/}.
\par\smallskip
{\tt Proof.} If even $\dim\Delta{\,>\,}30$, it has been shown in \cite{sz5} that $\cP$ is a Hughes plane 
(including the classical plane). Groups of dimension ${\ge}35$ have been dealt with in 5.6. We may suppose, therefore, that $\dim\Delta{\,<\,}35$; in the cases $\dim\Delta{\,\ge\,}31$ a few arguments 
in \cite{sz5} will be simplified. Note that $v^\Delta{\,=\,}W$ implies $\dim\Delta{\,>\,}36$, see 
\cite{cp} 96.\,19\hskip1pt--22 and use $\Theta$. \\
(a) {\it $\Theta{\,\le\,}\Delta_{[a,W]}$ if, and only if, $\Delta$ has no subgroup 
$\Gamma{\,\cong\,}\Gtwo$\/}: in fact, $\cF_\Gamma$ would be a flat subplane, which does not admit 
a torus group of homologies. If, on the other hand, $\Theta|_W{\,\ne\,}\1$, then there is some point 
$x$ such that $\cF_{\Delta_x}$ is a connected subplane, $\dim\Delta_x{\,>\,}12$, and Stiffness yields 
$\Gamma{\,:=\,}(\Delta_x)^1{\,\cong\,}\Gtwo$. \\
(b) {\it If $v^\Delta{\,=\,}V{\,\subseteq\,}W$, then $\dim V{\,=:\,}k{\,>\,}4$\/}: let $\Lambda$ denote the connected component of the stabilizer of points $u,v,w{\,\in\,}V$ and $c{\,\in\,}S{\,:=\,}av\sm\{a,v\}$.  
The dimension formula yields 
$$28{\,<\,}\dim\Delta{\,\le\,}3k{\,+\,}8{\,+\,}\dim\Lambda\,,\hskip5pt 
21{\,\le\,}3k{\,+\,}\dim\Lambda{\,\le\,}3k{\,+\,}14\,,\hskip5pt{\rm and} \ k{\,\ge\,}3\,.$$
If $\Theta|_W{\,\ne\,}\1$, then $\Delta$ has a subgroup $\Gamma{\,\cong\,}\Gtwo$ by step (a), and 2.15 implies $z^\Gamma{\,\approx\,}\sS_6$ for some $z{\,\in\,}V$ and $k{\,\ge\,}6$. 
If $\Theta|_W{\,=\,}\1$, however, then $\Theta$ induces a group of homologies on $\cF_\Lambda$,\: 
$\dim\cF_\Lambda{\,\ge\,}4$,\: $\dim\Lambda{\,\le\,}8$ by Stiffness, and then $3k{\,>\,}12$.   \\
(c) Put $\tilde\Delta{\,=\,}\Delta|_V{\,=\,}\Delta/\Kappa$. {\it The kernel $\Kappa$ has dimension 
$\dim\Kappa{\,<\,}8$, and $22{\,\le\,}\dim\tilde\Delta{\,\le\,}34$\/}:  
note that $\langle a,c,V\rangle{\,=\,}\cP$, so that $\Kappa$ acts freely on $S$. If 
$\dim\Kappa{\,=\,}8$, then each orbit $c^\Kappa$ is open in $S$ by \cite{cp} 96.11(a), and 
$\Kappa$ is transitive on $S$. The space $S$ is homotopy equivalent to $\sS_7$.  A 
maximal compact subgroup $\Phi$ of $\Kappa$ is a product of a torus and some  almost simple 
Lie groups $\Psi_\nu$, and $\pi_3\Psi_\nu{\,\cong\,}\ZZ$, see \cite{cp} 94.\,31(c) and 36.  
The exact homotopy sequence \cite{cp} 96.12  shows that $\pi_q\Phi{\,\cong\,}\pi_q\sS_7$ for 
each  $q{\,\ge\,}1$. Consequently $\Phi$ is trivial, but $\pi_7\sS_7{\,\cong\,}\ZZ$. This contradiction proves the claim; see also 2.17(b). \\
(d)  {\it $V$ is compact\/}: if not, then $\tilde\Delta$ has a normal subgroup $\Nu{\,\cong\,}\RR^k$ acting sharply transitive on $V$, and $\tilde\Delta$ is an extension of $\Nu$ by a transitive linear group, cf. \cite{cp} 96.16. We use the notation  of step~(b), and we write $\nabla{\,=\,}\Delta_{v,u}$.  
Suppose first that $k{\,\le\,}6$. Then $\Gtwo$ is  not a subgroup of~$\Delta$,\, 
$\Delta{\hskip1pt:}\nabla{\,=\,}2k$,\, $\nabla$ fixes a line of the vector space $\RR^k$ or $\CC^3$,\, $\Delta{\hskip1pt:\hskip1pt}\Lambda{\,<\,}2k{+}2{+}8$,\, 
$\Lambda{\,\not\hskip1pt\cong\,}\SU3\CC$,\,  $\dim\Lambda{\,\le\,}6$ by 2.6(\^e), and 
$\dim\Delta{\,<\,}28$ contrary to the assumption. \\ Next, let $k{\,=\,}7$. Then $\Delta_v$ induces a transitive, hence irreducible, linear group  on $\Nu$, and $\Delta_v$ has a subgroup 
$\Chi{\,\cong\,}\Gtwo$ by \cite{cp}  96.~\hskip-3pt16\,--22 and 94.27.   By step (a) there is some point $x$ such that $\cE{\,=\,}\langle a,x^\Theta,W\rangle$  is a connected subplane, and 
$\Delta_x|_\cE{\,=\,}\1$. Stiffness implies $\dim\Delta_x{\,\le\,}14$ and  $\dim\Delta{\,\le\,}30$. 
There are two possibilities:  $\Delta_v|_V$ is isomorphic either to 
$e^{\RR}{\hskip1pt\cdot\hskip1pt}\SO7\RR$ or  to  $e^{\RR}{\hskip1pt\cdot\hskip1pt}\Gtwo$. 
As $k$ is odd, $\Theta{\,\le\,}\Kappa$. In the first case, 
$\Delta{\,:\,}\Kappa{\,=\,}29$,\, $\dim\Kappa{\,=\,}1$,\, $\Kappa^1{\,=\,}\Theta$, and 
$\nabla{\,\cong\,}\Spin6\RR{\hskip1pt\cdot\hskip1pt}\Theta$ is a compact subgroup of 
$\Delta$.  Consequently,  $\nabla$ fixes the points of a $1$-dimensional subspace of~$\RR^7$, 
the orbit $c^\nabla$ is not open in $av$, and $\Lambda{\,=\,}\nabla_{\hskip-2pt c}$ has dimension $>8$, but this contradicts 2.6(\^e). In the second case, $\Delta{\,:\,}\Kappa{\,=\,}22$,\, 
$\dim\Kappa{\,=\,}7$,\, $\Theta{\,\triangleleft\,}\Kappa$, the representation of $\Chi$ on $\lie\Kappa$ shows  $\Kappa{\,\le\,}\Cs{}\Chi$, and 
$\Kappa|_{\cF_\Chi}$ would be too large. \\
Finally, if $k{\,=\,}8$, then $\Delta_v|_V$ contains one of the groups   $\SO8\RR$, $\SU4\CC$, or
$\U2\HH$, again by \cite{cp}  96.~\hskip-3pt16\,--22. In the first case, $\rk\Delta{\,=\,}5$ contrary to
2.10. In the second case, there exists a subgroup  $\Phi{\,\cong\,}\SU4\CC{\,\cong\,}\Spin6\RR$ in 
$\Delta$, see  \cite{cp} 94.27. The central involution $\sigma{\,\in\,}\Phi$ is not planar\, (or 
$\Phi|_{\cF_\sigma}{\,\cong\,}\SO6\RR$), and $\sigma$ is a reflection with center $v$ or axis $av$. 
By double transitivity, the elation group $\Delta_{[v,av]}$ is transitive (cf. \cite{cp} 61.19(b)\hskip1pt), and so is each group $\Delta_{[z,az]}$ with $z{\,\in\,}W$. It follows that $\Delta$ is transitive on $W$ and has 
a subgroup $\Spin9\RR$, but then $\rk\Delta{\,\le\,}4$ implies that $\Delta$ does not have a normal torus. The third case can be excluded by the same arguments. \\
(e) We can now prove the assertions.  
All doubly transitive transformation groups have been determined by Tits 1955, cf. \cite{cp} 96.16\hskip1pt--18 for results and a sketch of the proof. Consider the effective transformation 
group $(\tilde\Delta,V)$ defined in step (c). The space $V$ is compact if, and only if, $\tilde\Delta$ is a strictly simple  Lie group; in this case $V$ is a projective space or a sphere, and $\tilde\Delta$ is the corresponding projective or hyperbolic group. In the given dimension range $5{\,\le\,}\dim V{\,<\,}8$ and  $22{\,\le\,}\dim\tilde\Delta{\,\le\,}34$ there are exactly the following possibilities: 
$(\PSL4\CC,{\rm P}_3\CC)$,\, $(\PSU5(\CC,1),\sS_7)$, and $(\POpr8(\RR,1),\sS_6)$. 
 In the first case, the assumption 
on~$\Delta$ implies $\dim\Delta{\,>\,}\dim\tilde\Delta{\,=\,}30$. Choose $c$ as before, 
$u{\,\in\,}V$, and $w$ on the line $uv$ in the projective space ${\rm P}_3\CC$. Let 
$\Lambda{\,=\,}(\Delta_{c,u,w})^1$. By step (a), the central torus $\Theta$ induces a group of homologies on $\cF_\Lambda$, and $\dim\cF_\Lambda{\,>2\,}$. Stiffness yields 
$\dim\Lambda{\,\le\,}8$, and the dimension formula gives 
$\dim\Delta{\,\le\,}2{\cdot}6{\,+\,}2{\,+\,}8{\,+\,}\dim\Lambda{\,\le \,}30$, a contradiction.  \\
If $V{\,\approx\,}\sS_7$, then $\Delta$ has a subgroup $\Ypsilon{\,\circeq\,}\SU5(\CC,1)$ and a compact subgroup $\Phi{\,\cong\,}\SU4\CC$. The kernel $\Kappa$ satisfies 
$5{\,\le\,}\dim\Kappa{\,\le\,}7$, and the representation of $\Ypsilon$ on $\frl\hskip1pt\Kappa$ 
shows that $\Kappa{\,\le\,}\Cs{}\Ypsilon$. There are $6$ conjugate planar (diagonal) involutions in 
$\Phi$. Let $\beta$ be one of these. At most $3$ of the others induce reflections on $\cF_\beta$. 
Hence $\cF_{\beta,\hskip.6pt\beta'}{\;\ldot}\cF_\beta$ for some $\beta'$. As $\Kappa$ acts freely on a non-trivial point orbit in $\cF_{\beta,\hskip.6pt\beta'}$, this would imply $5{\,\le\,}\dim\Kappa{\,\le\,}4$.
Therefore $V{\,\approx\,}\sS_6$, and $\Delta$ contains $\Spin8(\RR,1)$   and has a subgroup 
$\Gtwo$. As in step (d), it follows that  $\dim\Delta{\,\le\,}30$. According to \cite{cp} 96.18, the 
group $\Delta$ is triply transitive on $V$. \qed
\par\smallskip
{\tt Remark.} In the  classical plane, $\Sigma'_{a,W}{\,\cong\,}\Spin{10}(\RR,1)$ contains
$\Delta{\,\cong\,}\Spin8(\RR,1){\times}\SO2\RR$. For any $v{\,\in\,}W$, the compact group 
$\Spin8\RR$ is a maximal semi-simple subgroup of $\Sigma_{a,v,W}$. Therefore $\Delta'$ has 
no fixed point on $W$, and $\cF_\Delta{\,=\,}\{a,W\}$. 
No other plane with such an action of $\Delta$ is known.
\par\medskip
{\bf 5.10 Normal torus.} {\it If  $\,\cF_\Delta{\,=\,}\{a,W\}, a{\,\notin\,}W$, and
 if $\Delta$ has a normal subgroup $\Theta{\,\cong\,}\TT$,  then $\dim\Delta{\,\le\,}30$\/}. \\
For  $\dim\Delta{\,>\,}30$, it has been proved in \cite{sz5} that $\cP$ is a Hughes plane, and 
then $\cF_\Delta{\,=\,}\emptyset$. \qed 
\par\medskip
{\bf 5.11 Normal vector groups.} {\it Assume that $\Delta$ is not semi-simple and that $\cP$ is not a Hughes plane. If $\dim\Delta{\,\ge\,}33$, then $\Delta$ has a minimal normal subgroup 
$\Theta{\,\cong\,}\RR^t$. In the cases $1{\,<\,}t{\,\ne\,}8,\,12$, the group $\Theta$ consists of 
elations with common axis or common center. If $t{\,=\,}8$ or $12$, then either $\Theta$ or
some other minimal normal subgroup  $\tilde\Theta$ has this property\/}.
\par\smallskip
This has been proved in \cite{sz10} Theorem A and Propositions 8 and 9. The arguments are rather involved, they shall not be indicated here. 
\par\medskip
{\bf 5.12 Corollary.} {\it If $\,\cF_\Delta{\,=\,}\{a,W\}, a{\,\notin\,}W$, and  if $\Delta$ has a  minimal
normal vector subgroup $\Theta$ of dimension $t{\,>\,}1$,  then $\dim\Delta{\,\le\,}32$\/}.

\par\bigskip 
{\Bf 6. Two fixed points}
\par\medskip
In the classical octonion plane, $\Sigma_{u,v}$ is an extension of a transitive translation group by
$(e^{\RR})^2{\cdot\,}\Spin8\RR$, cf. \cite{cp} 17.13.
Throughout this section let $\cF_\Delta{\,=\,}\langle u,v\rangle$.
\par\medskip
{\bf 6.0.} {\it If $\cF_\Delta{\,=\,}\langle u,v\rangle$  and if $\dim\Delta{\,\ge\,}18$, then $\Delta$ is 
a Lie group\/}. 
\par\smallskip
{\tt Proof.} (a) Suppose  that $\Delta$ is not a Lie group.
As in the proof of 4.0, it follows from \cite{lw3} that $x^\zeta{\,=\,}x$ for some $x{\,\notin\,}W$ and some $\zeta{\,\ne\,}\1$ in a compact central subgroup $\Nu$ of $\Delta$ such that $\Delta/\Nu$ is a Lie group. By hypothesis,  the orbit $x^\Delta$ is not contained in a line. Therefore  
$\cD{\,=\,}\langle x^\Delta\rangle{\,\le\,}\cF_\zeta$ is a proper connected subplane.  
Put $\Delta^{\hskip-1.5pt*}{=\,}\Delta|_\cD{\,=\,}\Delta/\Kappa$. If $\cD$ is flat, then 
$\Delta{\hskip1pt:\hskip1pt}\Kappa{\,\le\,}4$,\: $\dim\Kappa{\,=\,}14$,\: 
$\Delta^{\hskip-1.5pt*}{\,\approx\,}\RR^4$, and $\Kappa$ is a maximal compact subgroup of 
$\Delta$. Hence $\Kappa$ is connected by the Mal'cev-Iwasawa theorem \cite{cp} 93.10, and  
$\Kappa{\,\cong\,}\Gtwo$ would be a Lie group. If $\dim\cD{\,=\,}4$, then 
$\Delta{\hskip1pt:\hskip1pt}\Kappa{\,\le\,}8$ and $\dim\Kappa{\,\le\,}8$ by Stiffness. 
Consequently   $\cD{\;\ldot}\cP$,\: $\Kappa$ is compact and $\dim\Kappa{\,\le\,}7$ by Stiffness. 
According to Lemma 6.0$'$ below, $\Delta^{\hskip-1.5pt*}$ is a Lie group. We may assume, therefore, 
that $\Nu{\,\le\,}\Kappa$. Choose elements $w,L$ in $\cD$ with $w{\,\in\,}uv\sm\{u,v\}$ and 
$L{\,\in\,}\frL_v\sm\{uv\}$, let $z{\,\in\,}L\sm\cD$,  put $\Gamma{\,=\,}\Delta_{L,w}$ and 
$\Lambda{\,=\,}(\Gamma_{\hskip-1.5pt z})^1$. Then $\Delta{\hskip1pt:\hskip1pt}\Gamma{\,\le\,}8$,\: $\dim\Gamma{\,\ge\,}10$ and $\Kappa{\,\le\,}\Gamma$. From \cite{cp} 53.2 it follows that 
$\Gamma{\hskip1pt:\hskip1pt}\Lambda{\,=\,}\dim z^\Gamma{\,<\,}8$ and $\dim\Lambda{\,\ge\,}3$.   Moreover, $\Lambda$ is a Lie group, since $\Lambda\smcap\Kappa{\,\le\,}\Kappa_z{\,=\,}\1$.
If $\dim\Kappa{\,=\,}0$, then  
\cite{sz1} 1.10 implies that $\dim\Delta^{\hskip-2pt*}{\,>\,}18$, and $\cD$ is the classical qaternion plane. It follows that a maximal compact subgroup of $\Delta$ is locally isomorphic to $(\Spin3\RR)^3$, and $\Delta$ would be a Lie group. \\
(b) Thus $\dim\Kappa{\,>\,}0$. 
We need to show that $\cF_\Lambda$ is connected. By \cite{cp} 55.32, the group $\Kappa$
does not contain a pair of commuting involutions. Hence $\Kappa^1$ is a product of a commutative connected group $\Alpha$ with at most one almost simple factor $\Omega{\,\cong\,}\Spin3\RR$, 
see the structure of compact groups as stated in \cite{cp} 93.11. Obviously, $\Alpha$ is normal in 
$\Delta$, and \cite{cp} 93.19 implies that $\Alpha$ is in the center of 
$\Lambda{\hskip.5pt\cdot\hskip.5pt}\Alpha$. 
Now $\Alpha{\,\approx\,}z^\Alpha{\,\subseteq\,}\cF_\Lambda{\,<\,}\cP$, and $\cF_\Lambda$ is 
indeed connected except possibly if $\Kappa^1{\,\cong\,}\Omega$. \\
(c) Suppose first that  $\Kappa^1{\,\not\cong\,}\Omega$. Then 
$z^\Nu{\,\subset\,}\cF_\Lambda{\;\ldot}\cP$ by \cite{cp} 32.21  and 71.2, and 
$\Lambda{\,\cong\,}\Spin3\RR$. It follows that $\dim\Gamma{\,=\,}10$,\; 
$\Delta{\hskip1pt:\hskip1pt}\Gamma{\,=\,}8$, and $\dim\Delta{\,=\,}18$. 
The involution $\iota{\,\in\,}\Lambda$ induces a reflection on $\cD$; in fact,   $\iota|_\cD{\,\ne\,}\1$ because $z^\iota{\,=\,}z{\,\notin\,}\cD$. Hence $\Lambda$ acts faithfully on  $\cD$, and 
$\iota|_\cD$ cannot be planar by the  stiffness result 
\cite{sz1} 1.5(4), see also \cite{cp} 83.11. The 
axis of $\iota|_\cD$ is $S{\,=\,}uv\smcap\cD$\; (since $w^\Lambda{\,=\,}w$), the center is some point $a$.
As $\Delta{\hskip1pt:\hskip1pt}\Gamma{\,=\,}8$, we have $\dim w^\Delta{\,=\,}4$ for each choice 
of $w$, and $w^\Delta{\,=\,}S\sm\{u,v\}$ is a manifold by 
\cite{HK} Cor.\,5.5. Therefore
$S{\,\approx\,}\sS_4$, and $uv{\,\approx\,}\sS_8$. \\
(d) Let $\Pi{\,=\,}(\Delta_z)^1{\,=\,}(\Delta_{L,z})^1$. Then 
$\Pi_w{\,=\,}\Lambda$,\: $\dim\Pi{\,=\,}7$, and $\Pi$ is transitive on $S\sm\{u,v\}$.
Consequently  $\Pi/\Lambda{\,\approx\,}\HH^{\times}$ and $\Pi'\cong\Spin4\RR$,\, (use the exact
homotopy sequence \cite{cp} 96.12 together with 94.36). 
One of the central involutions of $\Pi'$ induces on  $\cD$ a reflection with axis 
$L\smcap\cD$. Denote this involution by 
$\rho$.  We will show that $\rho$  is a reflection of $\cP$ with axis $L$ and center~$u$. If not, 
then $L\smcap\cD{\,\subset\,}\cF_{\hskip-2pt\rho}{\;\ldot}\cP$, and $L\smcap\cF_{\hskip-2pt\rho}$ 
is a manifold (${\approx\;}\sS_4$) by \cite{cp} 92.16. In fact, 
$L\smcap\cD{\,=\,}L\smcap\cF_{\hskip-2pt\rho}$, see \cite{cp} 51.20, \cite{du} XVI,\,6.2(2), but 
$\rho$ fixes $z{\,\notin\,}\cD$ and even each point of $z^\Nu$. As $\dim L^\Delta{\,=\,}4$, 
we conclude from 2.9 that $\rho^\Delta\rho{\,\cong\,}\RR^4$ consists of translations with center 
$u$. Interchanging the r\^oles of $u$ and $v$, we find that also the translation group
$\Tau_{\hskip-1pt[v]}$ is isomorphic to $\RR^4$. Therefore the translation group $\Tau$ induces a transitive group of translations on $\cD$. Its complement contains $\Pi$. Hence 
$\dim\Delta^{\hskip-2pt*}{\,\ge\,}8{+}7$ and $\dim\Kappa{\,\le\,}3$. \\
(e)  If $\dim\Kappa{\,=\,}1$, then $\cD$ is a near-field plane (cf. 
\cite{sz1} 1.10), its automorphism 
group is described in \cite{cp} 82.24. We conclude that 
$\Delta_a$ has a unique normal subgroup $\Phi{\,\cong\,}\Spin4\RR$. It follows that 
$\Pi'{\,=\,}\Phi$, and $z^\Phi{\,=\,}z$. Analogously, there is a point $y{\,\in\,}au\sm\cD$ such 
that $y^\Phi{\,=\,}y$, but the reflection in $\Pi'_{[L]}$ does not fix $y$. This contradiction shows that 
$\dim\Kappa{\,\ge\,}2$. \\
(f) Note that each involution of $\Pi'$ is central. Therefore $\Pi'{\,\le\,}\Delta_a$ and the radical 
$\Gamma{\,=\,}\sqrt{\Delta_a}$ is $4$-dimensional. If $\dim\Kappa{\,=\,}2$, then 
complete reducibility shows $\Pi'{\,\le\,}\Cs{}\Gamma/\Nu$ and hence $\Pi'{\,\le\,}\Cs{}\Gamma$.
The same is true in the case $\dim\Kappa{\,=\,}3$ because $\Kappa/\Nu$ is a torus.
Thus $\Pi'$ is the unique Levi complement of $\Gamma$, but this is impossible by the argument 
in step (e). \\
(g) Finally let $\Kappa^1{\,\cong\,}\Omega$. Again $\dim\Delta^{\hskip-2pt*}{\,\ge\,}15$, and 
\cite{sz1}  5.1-4 implies that $\Delta^{\hskip-2pt*}$ contains a transitive group of translations, 
their center will be chosen as the point $v$. In particular, $L\smcap\cD{\,\approx\,}\sS_4$.
If $\dim\Lambda{\,>\,}3$,  if $\Lambda$ is solvable, or if $\Lambda{\,\circeq\,}\SL2\RR$, then 
$\Rho{\,=\,}(\Cs\Lambda\Omega)^1{\,\trianglelefteq\,}\Lambda$ has positive dimension, 
since $\Lambda/\Cs\Lambda\Omega{\,\le\,}\Aut\Omega{\,\cong\,}\SO3\RR$. 
It follows that $z^\Rho{\subset\,}\cF_\Rho{\;\ldot}\cP$, and $\Rho{\,\cong\,}\TT$ or 
$\Rho{\,\cong\,}\Spin3\RR$. Hence there exists a planar involution $\rho$ in the center of
$\Lambda\Omega$. In the case $\cD\smcap\cF_\rho{\,=\,}\cC{\;\ldot}\cD,\cF_\rho$, the group 
$\Omega$ would act trivially on $\cC$, which contradicts Stiffness. Consequently $\rho$ induces 
on $\cD$ a reflection with axis $uv\smcap\cD{\,\subset\,}\cF_\rho$, and then   
$S{\,=\,}L\smcap\cF_\rho{\,\approx\,}\sS_4$, again by   \cite{cp} 92.16 and 51.20.
According to $(\dagger)$, the action of $\Omega$ on $S$ is equivalent to the  action of $\HH'$ 
on $\HH$ by multiplication, $S/\Omega{\,=\,}I$ is an interval, and $\Nu$ acts trivially on $I$.  
Therefore  $z^\Omega{\,=\,}z^{\Omega\Nu}$, but $(\Omega\Nu)_z{\,=\,}\1$, a contradiction. \\
(h) Only the cases $\Kappa'{\,\cong\,}\Omega{\,\circeq\,}\Lambda$ remain. If $\Lambda$ is simply 
connected, then $\Omega$ centralizes the involution $\rho{\,\in\,}\Lambda$. Again 
$\rho|_\cD$ is a reflection, and one may reason as in step (g). Similarly, 
$\SO3\RR{\,\cong\,}\Lambda{\,\le\,}\Cs{}\Omega$ is impossible. Therefore 
$\Lambda\Omega{\,\cong\,}\SO4\RR$ contains only one involution and a different argument is needed. We will show that $\Lambda\Kappa$ is a maximal compact subgroup of the connected 
group $\Xi{\,=\,}\Delta_{w,L}$. In fact, $\Delta{\hskip1pt:\hskip1pt}\Xi{\,\le\,}8$,\; 
$\dim\Xi{\,\ge\,}10$,\;$\Xi{\hskip1pt:\hskip1pt}\Lambda{\,=\,}\dim z^\Xi{\,<\,}8$,\; 
$\dim\Xi{\,=\,}10$,\; $\dim w^\Delta{\,=\,}\dim L^{\Delta_w}{\,=\,}4$ for each choice of $w$ and $L$, both orbits are simply connected, hence $\Xi$ is connected and so is $\Xi|_\cD$.
Recall from step (g) that $\Xi|_\cD$ contains a transitive group $\Theta{\,\approx\,}\RR^4$ of 
translations\; (note that each compact subgroup of $\Theta$ is trivial by \cite{cp} 55.28). We have 
$\Lambda{\,\cong\,}\Lambda|_\cD{\,<\,}\Xi|_\cD$ and $\Xi{\hskip1pt:\hskip1pt}\Kappa{\,=\,}7$. 
Consequently $\Xi|_\cD{\,\cong\,}\Lambda{\ltimes}\Theta$, and $\Lambda\Kappa$ is a
compact subgroup of $\Xi$ of maximal dimension. By the Mal'cev-Iwasawa theorem \cite{cp} 93.10 
the group $\Lambda\Kappa$ is connected and coincides with the Lie group $\Lambda\Omega$.
\qed  
\par\medskip 
{\bf 6.0$'$ Lemma.} {\it Suppose that $\Delta$ is a group of automorphisms of an 
$8${\rm -dimensional} plane~$\cP$. If $\Delta$ fixes $2$ distinct points $u,v$ and exactly one 
line $uv$, and if $\dim\Delta{\,\ge\,}8$, then $\Delta$ is a Lie group\/}.
\par\smallskip
{\tt Proof.} If $\Delta/\Nu$ is a Lie group, but the compact central subgroup $\Nu$ is not, 
then  the arguments of 4.0 show that for some $x{\,\notin\,}uv$ the orbit $x^\Delta$ generates 
a proper connected sublane $\cD$. Put $\Delta|_\cD{\,=\,}\Delta/\Kappa$ and apply the stiffness results \cite{sz1} 1.5 to $\Kappa$. If $\cD$ is flat, then $\Delta{\hskip1pt:\hskip1pt}\Kappa{\,\le\,}4$ 
and $\dim\Kappa{\,\le\,}3$ \ (see also \cite{cp} 83.12). Hence $\cD{\;\ldot}\cP$ and 
$\dim\Kappa{\,\le\,}1$. As $\Delta|_\cD$ is a Lie group by \cite{cp} 71.2, we may assume 
$\Nu{\,\le\,}\Kappa$. In the case $\dim\Kappa{\,=\,}0$ it follows fom \cite{cp} 73.22 that $\Delta$ 
does not have two fixed points in $\cD$. Therefore $\dim\Kappa{\,=\,}1$. Choose a line $L$ of
$\cD$ in $\frL_v\sm\{uv\}$ and points $w{\,\in\,}uv\smcap\cD$ and $z{\,\in\,}L\sm\cD$. 
Put $\Lambda{\,=\,}(\Delta_{z,w})^1$. Then $\dim z^{\Delta_L}{\,<\,}4$ by \cite{cp} 53.2, and 
$\dim\Lambda{\,>\,}0$. Let $\Xi{\,=\,}\Kappa^1\Nu$. Then \cite{cp} 93.19 implies 
$\Lambda{\,\le\,}\Cs{}\Xi$. Therefore $z^\Xi{\,\subseteq\,}\cF_\Lambda{\,<\,}\cP$,\:  
$\cF_\Lambda$ is connected, but then $\Xi$ would be a Lie group  by \cite{cp} 32.21 and 71.2. \qed  
\par\medskip
{\bf 6.1 Lie.} {\it If $\Delta$ is semi-simple of dimension $\dim\Delta{\,\ge\,}14$, then $\Delta$ is a 
Lie group\/}.
\par\smallskip
This can be {\tt proved} in the same way as 5.1. The sharper bound is obtained because
$\Gtwo$ is a Lie group, and because  $\cD{\,=\,}\cD^\Delta{\,<\,}\cP$ implies 
$\dim\Delta|_\cD{\,\le\,}10$ instead of $11$, see \cite{cp} 7.3. 
{\tt Remark.} The result is also true if $\Delta$ fixes more than $2$ collinear points but only one line.
\par\medskip
{\bf 6.2 Semi-simple groups.}  {\it If $|\cF_\Delta|{\,=\,}3$ and $\Delta$ is   semi-simple, 
then $\dim\Delta{\,\le\,}20$ or, conceivably, $\Delta$ is a product of $7$ factors each of which is isomorphic to the simply connected covering group $\Omega$ of $\SL2\RR$. If $\Delta$ is almost simple,  then $\dim\Delta{\,\le\,}16$\/}.
\par\smallskip
{\tt Remark.} Results and proof are the same as in the case that $\cF_\Delta$ is a flag, see 5.2.
\par\medskip
{\bf 6.3 Normal torus.} {\it Suppose that $\Delta$ has a one-dimensional compact connected normal subgroup $\Theta$. If $\dim\Delta{\,\ge\,}18$, 
then there exists a $\Delta$-invariant near-field plane $\cH{\,\ldot}\cP$$;$ either $\dim\Delta{\,=\,}18$, or $\cH$ is classical and $\dim\Delta{\,\le\,}20$\/}.  
\par\smallskip
{\tt Proof.} As $\Delta$ is a Lie group, $\Theta{\,\cong\,}\TT$ and the involution $\iota{\,\in\,}\Theta$ is in the center of $\Delta$. Hence $\iota$ is planar.  
Put $\cH{\,=\,}\cF_\iota$ and $\Delta|_\cH{\,=\,}\Delta/\Kappa$. If 
$\Theta|_\cH{\,\ne\,}\1$, Stiffness and  
\cite{sz1} 5.3 imply $\dim\Delta{\,\le\,}9{+}3$. 
Hence $\Theta{\,\le\,}\Kappa$,\: $\Kappa^1{\,\not\hskip1pt\cong\,}\Spin3\RR$,\: 
$\dim\Kappa{\,=\,}1$, and $\Delta{\hskip1pt:\hskip1pt}\Kappa{\,\ge\,}17$. All planes $\cH$ 
with  a group of dimension at least $17$ are described in 
\cite{sz1} 1.10.  Only in the case of the near-field planes 
$\Delta$ fixes two distinct points. Either $\cH$ is  classical, or
$\Delta{\hskip1pt:\hskip1pt}\Kappa{\,=\,}17$ and $\dim\Delta{\,=\,}18$. 
\qed
\par\smallskip
{\tt Remark.} If $\dim\Delta{\,\ge\,}17$, then $\Delta$ induces on $\cH$ a group of dimension
${\ge\,}16$. Such planes are translation planes, they have been determined explicitly
by H\"ahl, cf. \cite{sz1}  3.3.
\par\medskip
Without assumption on $\cF_\Delta$, the following has been proved in \cite{sz10}:
\par\smallskip
{\bf 6.4 Theorem.} {\it If $\Delta$ has a normal vector subgroup and if $\dim\Delta{\,\ge\,}33$, 
then $\Delta$ fixes some element, say a line $W$, and $\Delta$  has a minimal normal subgroup 
$\Theta{\,\cong\,}\RR^t$  consisting  of axial collineations with common axis $W$. Either 
$\Theta{\,\le\,}\Delta_{[a,W]}$ is a group of homologies and $t=1$,  or $\Theta$ is contained 
in the translation group $\Tau = \Delta_{[W,W]}$\/}.
\par\medskip
The next result has been proved in \cite{HS} Theorem~1 under the  assumption  
$\dim\Delta{\,>\,}33$\,: 
\par\smallskip 
{\bf 6.5 Theorem.} {\it If $\dim\Delta{\,\ge\,}34$, then the group $\Tau$ of translations in $\Delta$ 
satisfies $\dim\Tau{\,\ge\,}15$. Either $\Delta$ has a subgroup $\Upsilon{\,\cong\,}\Spin7\RR$ 
and $\dim\Delta{\,\ge\,}36$, or $\Tau$ is transitive, $\dim\Delta{\,=\,}34$, and there exists a maximal semi-simple subgroup $\Upsilon{\,\cong\,}\Spin6\RR$ of $\Delta$\/}. 
\par\smallskip
{\bf 6.6 Normal vector subgroup.} {\it If $\cF_\Delta{\,=\,}\{u,v,W\}$ and if $\dim\Delta{\,=\,}33$, then the translation group $\Tau$ is transitive\/}. 
\par\smallskip
{\tt Proof.} In steps (a--r),  the weaker assumption $\dim\Delta{\,\ge\,}33$ suffices; the exact value 
 of $\dim\Delta$ is needed only in steps (s--v). \\
(a) From 6.\,2 and 3 it follows that 6.4 applies. Suppose that $\Delta_{[a]}{\,\ne\,}\1$ for some 
$a{\,\notin\,}W$. Then $\Delta{\,=\,}\Delta_{[a]}{\ltimes}\Tau$ by \cite{cp} 61.20, and 
$a^\Delta{\,=\,}a^\Tau$ is not contained in a line. Hence $\Tau$ is a vector group $\RR^k$, and 
each  subgroup $\Tau_{\hskip-1pt[z]}$ with $z{\,\in\,}W$ is connected by \cite{cp} 61.9. Let 
$\1{\,\ne\,}\tau{\,\in\,}\Tau_{\hskip-1pt[z]}$ with $z{\,\ne\,}u,v$, and put $\nabla{\,=\,}(\Delta_a)^1$,\: 
$\Lambda{\,=\,}(\nabla_{\hskip-2.5pt a^\tau}\hskip-2pt)^1$. Then 
$$33{\,\le\,}\dim\Delta{\,=\,}\Delta{:}\nabla{\,+\,}\nabla{:\hskip1pt}\Lambda{\,+\,}\dim\Lambda{\,\le\,}
2k{\,+\,}\dim\Lambda{\,\le\,}2k{\,+\,}14\quad{\rm and}\quad k{\,\ge\,}10\,.$$
By the stiffness result 2.6(e), either $\Lambda{\,\cong\,}\Gtwo$ or $\dim\Lambda{\,\le\,}8$. 
In the first case $\Tau{\smcap}\Cs{}\Lambda$ acts effectively on the flat plane $\cF_\Lambda$ 
and $\dim(\Tau{\smcap}\Cs{}\Lambda){\,\le\,}2$. As $\Lambda$ fixes $\tau$ and each non-trivial representation of $\Lambda$  on $\Tau$ has dimension $7$ or $14$, it follows that $k{\,\ge\,}15$.
In the second case, $k{\,\ge\,}13$ and $\dim\Tau_{\hskip-1pt[z]}{\,\ge\,}5$ for each $z{\,\in\,}W$. \\
(b) {\it The elements of a minimal normal subgroup $\Theta$ of $\Delta$,  
$\RR^t{\,\cong\,}\Theta{\,\le\,}\Tau$, have center $u$ or~$v$, say 
$\Theta{\,\le\,}\Tau_{\hskip-1pt[v]}$, or $\Tau$ is transitive and $\cP$ is classical\/}. 
In fact, $t{\,\le\,}8$\: (or   $L{\,\in\,}\frL_v\sm\{W\}$ implies $\1{\,\ne\,}\Theta_L{\,=\,}
\Theta\smcap\Tau_{\hskip-1.5pt L}{\,=\,}\Theta\smcap\Tau_{\hskip-1pt[v]}{\,\triangleleft\,}\Delta$, and $\Theta$ would not be minimal). If some $\tau{\,\in\,}\Theta$ has a center $z{\,\ne\,}u,v$, 
then  $17{\,\le\,}\dim\nabla{\,\le\,}t{\,+\,}\dim\Lambda$,\: $\dim\Lambda{\,>\,}8$, 
$\Lambda{\,\cong\,}\Gtwo$,\: $t{\,>\,}2$,\:  $\Lambda$ acts non-trivially on 
$\Theta$, and $t{\,>\,}7$.  
By minimality, $\Delta$ induces an irreducible group $\tilde\Delta$ on 
$\Theta{\,\cong\,}\RR^8$, and $\tilde\Delta'$ is semi-simple and properly contains $\Gtwo$. 
The list \cite{cp} 95.10 shows that $\tilde\Delta'$ has a subgroup $\Spin7\RR$. This simply connected group can be lifted to a subgroup $\Upsilon$ of $\Delta$, see \cite{cp} 94.27. 
The central involution $\sigma{\,\in\,}\Upsilon$ inverts the elements of $\Theta$ and hence fixes 
$z$. Therefore $\sigma$ is a reflection with axis $W$, the center of $\sigma$ will be denoted by 
$a$. From 2.9 and step (a) it follows that $\sigma^\Delta\sigma{\,=\,}\Tau{\,\cong\,}\RR^k$ with
$k{\,\ge\,}14$ and $a^\Delta{\,=\,}a^\Tau$. Moreover, $\Upsilon$ acts faithfully on 
$\Tau_{\hskip-1pt[u]}$ and $\Tau_{\hskip-1pt[v]}$, but then $\Tau$ is transitive, and 
$\cP{\,\cong\,}\cO$ by \cite{cp} 81.17.  \\
(c) {\it If $\Delta$ has a  subgroup $\Upsilon{\,\cong\,}\Spin7\RR$, then 
$\dim\Tau{\,\ge\,}15$\/}.  By the last part of 2.10, the central involution $\sigma{\,\in\,}\Upsilon$ is a reflection. If $\sigma{\,\in\,}\Delta_{[W]}$, steps (a,b) show transitivity of $\Tau$. If $\sigma$ has center $u$ or $v$, say $\sigma{\,\in\,}\Delta_{[v]}$, the dual of 2.9 implies that $\Upsilon$ acts 
effectively on $\Tau_{\hskip-1pt[v]}$, and then $\Tau_{\hskip-1pt[v]}{\,\cong\,}\RR^8$ by \cite{cp}
95.10. It follows that $\Delta$ induces on $\frL_u$ a trans\-itive group $\Delta/\Delta_{[u]}$. 
Either $\Ypsilon$ is contained in a subgroup of $\Delta$ of type 
${\rm D}_4$,\: $\dim\Delta{\,\ge\,}36$, 
and the claim follows from 6.5, or $\Delta/\Delta_{[u]}$ is contained in 
$e^{\RR}{\cdot\hskip1pt}\Spin7\RR{\hskip1pt\ltimes}\RR^8$ and has dimension 
$\Delta{\hskip1pt:\hskip1pt}\Delta_{[u]}{\,\le\,}30$. In the latter case $\dim\Delta_{[u]}{\,\ge\,}3$. 
According to   \cite{cp} 61.20, the connected component $\Xi$ of $\Delta_{[u]}$ is isomorphic to 
$\Delta_{[u,xv]}{\ltimes}\Tau_{\hskip-1pt[u]}$ and   $\dim\Tau_{\hskip-1pt[u]}{\,>\,}0$. 
Note that the homology group $\Delta_{[u,xv]}$ acts freely on $\Tau_{\hskip-1pt[u]}$.
By 2.15(b), the fixed points  of $\Upsilon$ on the axis of $\sigma$ form a circle~$C$. 
Suppose that  $\Xi{\,\le\,}\Cs{\hskip-2pt}\Upsilon$. Then $\Xi$ acts effectively 
on $C$, and Brouwer's theorem \cite{cp} 96.30 would imply $\dim\Delta_{[u]}{\,\le\,}2$. 
Consequently $\Upsilon$ induces on the Lie algebra $\lie\Xi$ a group $\SO7\RR$ and 
$\dim\Xi{\,\ge\,}7$. Hence $\dim\Tau_{\hskip-1pt[u]}{\,\ge\,}\frac{1}{2}\dim\Xi{\,\ge\,}4$, the action of $\Upsilon$  on $\Tau_{\hskip-1pt[u]}$ is not trivial by 2.15, and then 
$\dim\Tau_{\hskip-1pt[u]}{\,\ge\,}7$. 
In particular: \\
(c$'$) {\it If the central involution of $\Ypsilon$ is a reflection with center $v$, 
then $\Tau_{\hskip-1pt[v]}$ is transitive\/}.  \\ 
{\bf Change of notation.} From now on $\Theta{\,\cong\,}\RR^t$ will denote a minimal 
$\nabla$-invariant subgroup of $\Tau_{\hskip-1pt[v]}$, where $\nabla{\,=\,}(\Delta_a)^1$ for a fixed $a{\,\notin\,}W$; normality of $\Theta$ in $\Delta$ will not be needed. Throughout, let
$w{\,\in\,}S{\,:=\,}W\sm\{u,v\}$ and $\Omega{\,=(\nabla_{\hskip-1.5pt w})^1\,}$.         \\
(d) {\it If $t{\,=\,}1$ and if $\Delta$ has a subgroup $\Gamma{\,\cong\,}\Gtwo$, then $\Tau$ is transitive\/}. \\
Consider the radical $\Rho{\,=\,}\sqrt\Delta$ and note that \cite{HS} Lemma~4 
combined with $\Theta{\,\cong\,}\RR$ yields $\dim\Rho{\,\le\,}19$. Suppose that 
$\Gamma$ is a Levi  complement of $\Rho$ in $\Delta$. Then 
$\Delta{\,=\,}\Rho\Gamma$,\: $\dim\Rho{\,=\,}19$, and $\Rho$ is transitive on the affine point 
set $P\sm W$. Thus, up to conjugacy,   $a^\Gamma{\,=\,}a$ and $\Gamma{\,\le\,}\nabla$. 
Stiffness  and 2.15(a) imply that  $\cF_\Gamma$ is a flat subplane. Consequently, 
$\dim\Cs\nabla\Gamma{\,\le\,}2{\,<\,}\break\dim(\nabla\smcap\Rho){\,=\,}3$. Therefore $\Gamma$ 
acts non-trivially on $\nabla\smcap\Rho$, and $\dim(\nabla\smcap\Rho){\,\ge\,}7$, a contradiction.
Hence $\Gamma$ is properly contained in a maximal semi-simple subgroup $\Psi$ of 
$\Delta$. By Stiffness, $\Cs{}\Gamma$ acts almost effectively on $\cF_\Gamma$, the induced 
group  is solvable by \cite{cp} 33.8. Because of 2.16, it follows that $\Psi$ is an almost simple
orthogonal group containing  $\Ypsilon{\,\cong\,}\Spin7\RR$.  
Step (c) implies that $\dim\Tau{\,\ge\,}15$ and $\dim\Delta{\,\ge\,}36$. Hence 
$\dim\nabla{\,\ge\,}21$ and $\dim\nabla_{\hskip-2pt c,w}{\,\ge\,}12$ for 
$a{\,\ne\,}c{\,\in\,}a^\Theta$. By Stiffness $(\nabla_{\hskip-2pt c,w})^1{\,\cong\,}\Gtwo$, 
and we may assume that $\Gamma{\,\le\,}\nabla$.
If  $\Upsilon{\,\le\,}\nabla$, the central involution $\sigma{\,\in\,}\Upsilon$ is a reflection with axis $av$, since $\Theta{\,\le\,}\Cs{}\Upsilon$, and~(c$'$)  shows that 
$\Tau_{\hskip-1pt[u]}{\,\cong\,}\RR^8$. 
The group $\Tau_{\hskip-1pt[v]}$ is a product of $\Theta{\,=\,}\Theta^\Upsilon$ and a 
$7$-dimensional factor; hence $\Tau$ is transitive in this case.  From $a^\Upsilon{\,\ne\,}a$ 
we will derive a contradiction. First let $\Upsilon_{\hskip-1pt[v]}{\,=\,}\langle \sigma\rangle$ as in 
step (c$'$).   Then $\Tau_{\hskip-1pt[v]}$ is transitive and, up to conjugacy, the axis of $\sigma$ contains the point $a$. 
Now $\Gamma{\,\le\,}\nabla\smcap\Upsilon{\,=\,}\Upsilon_{\hskip-2pt a}$,  
but  $a^\Upsilon{\,\approx\,}\sS_6$ by 2.15, and $\Upsilon_{\hskip-1.5pt a}{\,\cong\,}\Spin6\RR$.
The case that the central involution $\sigma$ of $\Upsilon$ has center $u$
can be dealt with is analogously. \\ 
(e) {\it If $t{\,=\,}1$ and if $\Delta$ has no subgroup $\Gtwo$, then $\Tau$ is transitive\/}. \\ 
Let $a{\,\ne\,}c{\,\in\,}a^\Theta$. Then  $\dim\nabla_{\hskip-2pt c}{\,\ge\hskip1pt}16$, and 2.6(e) implies  $\Lambda{\,=\,}(\nabla_{\hskip-2pt c,w})^1{\,\cong\,}\SU3\CC$. Con\-sequently, 
$a^\Delta$ is open in $P$,\: $\nabla_{\hskip-2pt c}$ 
is transitive on $S$, and a maximal compact subgroup $\Ypsilon$ of $\nabla_{\hskip-2pt c}$ is 
isomorphic to $\SU4\CC{\,\cong\,}\Spin6\RR$\:  (see 2.17(d)\,). The central involution 
$\sigma{\,\in\,}\Ypsilon$ is a reflection with axis $av$ and center $u$,\: 
$\sigma{\,\in\,}\nabla_{\hskip-1.5pt[u]}$. As $(av)^\Delta{\,\ne\,}av$, it follows from the dual of 2.9 
that $\sigma^\Delta\sigma{\,=\,}\Tau_{\hskip-1pt[u]}$ is a vector group. The representation of
$\Ypsilon$ on $\Tau_{\hskip-1pt[u]}$ is faithful. Therefore $\Tau_{\hskip-1pt[u]}{\,\cong\,}\RR^8$, 
and $\Delta$ induces on $\frL_v\sm\{W\}$ a doubly transitive group $\Delta/\Delta_{[v]}$, so that 
$\Xi{\,=\,}\Delta_{av}/\Delta_{[v]}$ is a trans\-itive linear group and $\Xi'$ is semi-simple (cf. \cite{cp} 95.6 and 96.16). As 
$\nabla_{\hskip-2pt c}\smcap\Delta_{[v]}{\,\le\,}\Delta_{[v,av]}{\,=\,}\1$, 
we have $\Upsilon{\,\cong\,}\SU4\CC{\,\le\,}\Xi'$. Moreover 
$\Cs\Xi\Upsilon{\,\le\,}\CC^{\times}$, and $\Xi'$ is almost simple. 
If  $\Upsilon{\,<\,}\Xi'$, then  \cite{cp} 96.10 would show that
$\Gtwo{\,<\,}\Spin7\RR{\,<\,}\Xi$ or $\dim\Xi'{\,\ge\,}28$;  the first possibility is excluded by 
the assumption in step (e),  in the second case $\dim\Delta{\,\ge\,}36$ and theorem 6.5 
would imply $\Xi'{\,\cong\,}\Spin7\RR$. Hence $\Xi{\,\le\,}\CC^{\times}{\cdot\,}\SU4\CC$,\: 
$16{\,\le\,}\dim\Xi{\,\le\,}17$, and $\dim\Delta_{[v]}{\,\ge\,}8$.   
Either $\Delta_{[v]}{\,=\,}\Tau_{\hskip-1pt[v]}$ is transitive, or 
$\Delta_{[v]}$ contains a homology with  axis~$au$. As $\dim(au)^\Delta{\,=\,}8$, the dual of 2.9 
implies again $\Tau_{\hskip-1pt[v]}{\,\cong\,}\RR^8$. \\
(f) {\it If $t{\,=\,}2$ and if $\dim\Lambda{\,<\,}8$ for each stabilizer of a quadrangle, then  
 $\Tau$ is transitive\/}.\break {\tt Proof.}  Let  $a{\,\ne\,}c{\,\in\,}a^\Theta$, and put 
$\Lambda{\,=\,}(\nabla_{\hskip-2pt c,w})^1$. Then $\Delta{\hskip1pt:\hskip1pt}\Lambda{\,\le\,}26$,\: $\dim\Lambda{\,=\,}7$,\: $\dim\nabla{\,=\,}17$,\: $\nabla_{\hskip-2pt c}$ is transitive on $S$ and on $au\sm\{a,u\}$, moreover, $a^\Delta{\,=\,}P\sm W$. Consequently,  $\Delta_{au}$ is doubly transitive on $au\sm\{u\}$. As mentioned in 5.9(e), all such actions are known (Tits 1955, cf. in particular \cite{Vl} Satz 1). As $\nabla_{\hskip-2pt c}$ acts faithfully on $au$ and 
$\dim\nabla_{\hskip-2pt c}{\,=\,}15$, it follows  that $\nabla_{\hskip-2pt c}{\,\cong\,}\SL2\HH$\: (see \cite{cp} 96.10 or \cite{Vl} and note that $\nabla_{\hskip-2pt c}$ is not compact). 
A maximal compact subgroup $\Gamma$ of $\nabla_{\hskip-2pt c}$ is isomorphic to 
$\U2\HH{\,\cong\,}\Spin5\RR$, and the central involution $\sigma$ of $\Gamma$ is a reflection with axis $av$\: (use the last part of 2.10). Therefore 
$\sigma^\Delta\sigma{\,=\,}\Tau_{\hskip-1pt[u]}{\,\cong\,}\RR^8$.
The group $\Delta/\Delta_{[v]}$ induced by $\Delta$ on the pencil $\frL_v\sm\{W\}$ is also doubly trans\-itive, and $\nabla_{\hskip-2pt c}{\,\le\,}\Ypsilon{\,:=\,}\Delta_{av}/\Delta_{[v]}$.  
As $\Theta{\,\le\,}\Xi{\,:=\,}\Delta_{[v]}$, we have $15{\,\le\,}\dim\Ypsilon{\,\le\,}23$. 
By \cite{Vl} either $\dim\Ypsilon'{\,\le\,}18$ or  $\Ypsilon'{\,\cong\,}\Sp4\CC$. The second possibility 
can be excluded however, since  $\SL2\HH$ is not contained in $\Sp4\CC$\: (this follows easily 
from the complete reducibility of the adjoint representation of the smaller group on the Lie algebra 
of the larger one and the fact  that  each representation of $\SL2\HH$ in dimension $<6$ is trivial. 
It is also a consequence of the far more general results of Tits \cite{Ti} p.160/61). Hence 
$\dim\Ypsilon{\,\le\,}19$ and  $\dim\Xi{\,\ge\,}6$. If $\Xi$ contains a non-trivial homology, then 
$\Tau_{\hskip-1pt[v]}$ is transitive by \cite{cp} 61.20. Assume now that $\Xi{\,=\,}\Tau_{\hskip-1pt[v]}$ 
has dimension $\dim\Xi{\,<\,}8$. Then $\Xi{\,\le\,}\Cs{}\nabla_{\hskip-2pt c}{\,\le\,}\Cs{}\Lambda$\: (since $\Theta{\,<\,}\Xi$), and $\Lambda|_{a^\Xi}{\,=\,}\1$, but then $\Lambda$ would be 
trivial. \\  
(g) {\it If $t{\,=\,}2$, then $\Delta$ is transitive on the affine point space $P\sm\{W\}$\/}. \\ If not, then 
$\dim\nabla_{\hskip-2pt c}{\,>\,}15$,\: $\Lambda{\,\cong\,}\SU3\CC$,\: $\dim\nabla{\,=\,}18$, and 
$\nabla_{\hskip-2pt c}$ is transitive on $S$ and on $au\sm\{a,u\}$. Consequently, $\Delta_{au}$ 
is doubly transitive on $au\sm\{u\}$. The exact homotopy sequence shows 
$\nabla\hskip.5pt'_{\hskip-3pt c}{\,\cong\,}\SU4\CC$, see \cite{sz5} Lemma (5) or 2.17(d). The central involution in 
$\nabla\hskip.5pt'_{\hskip-3pt c}$ is a reflection with axis $av$, and $\Tau_{\hskip-1pt[u]}{\,\cong\,}\RR^8$ 
by the dual of 2.9. It follows that $\nabla$ induces on $\Tau_{\hskip-1pt[u]}$ a subgroup of 
$\CC^{\times}{\cdot\hskip1pt}\SU4\CC$. In particular, 
$\nabla{\hskip1pt:\hskip1pt}\nabla_{\hskip-2pt[au]}{\,\le\,}17$  and $\Delta_{[v,au]}{\,\ne\,}\1$. 
As $\Theta{\,<\,}\Tau_{\hskip-1pt[v]}$ and $\Lambda$ acts non-trivially on $\Tau_{\hskip-1pt[v]}$, 
we have even $\dim\Tau_{\hskip-1pt[v]}{\,=\,}8$, and $\Tau$ would be transitive contrary to the assumption. \\
(h) {\it If $t{\,=\,}2$ and $\nabla$ is transitive on $S{\,=\,}uv\sm\{u,v\}$, then $\Tau$ is transitive\/}.\\ 
The {\tt proof} is similar to that in step (g). Because of step (f), we may assume that 
$(\nabla_{\hskip-2pt w})'{\,=\,}\Lambda{\,\cong\,}\SU3\CC$ for $w{\,\in\,}S$. Then  a maximal 
compact subgroup $\Phi$ of $\nabla$ is isomorphic to $\SU4\CC$, see \cite{sz5} Lemma (5) or 2.17(d).
Again the central involution of $\Phi$ is a reflection in 
$\nabla_{\hskip-1pt[u,av]}$, and $\Tau_{\hskip-1pt[u]}{\,\cong\,}\RR^8$ by 2.9. The group 
$\Ypsilon{\,=\,}\Delta_{av}/\Delta_{[v]}$ induced by $\Delta_{av}$ on the pencil $\frL_v\sm\{W\}$
acts faithfully on $\Tau_{\hskip-1pt[u]}$. If $\Tau_{\hskip-1pt[v]}{\,<\,}\Delta_{[v]}$, then 
$\Tau_{\hskip-1pt[v]}{\,\cong\,}\RR^8$ by step (g) and \cite{cp} 61.20. Hence we may assume that 
$\Delta_{[v]}{\,=\,}\Tau_{\hskip-1pt[v]}$. Representation theory as summarized in \cite{cp} 95.10 
shows that $\Ypsilon{\,\le\,}\CC^{\times}{\rtimes\,}\SU4\CC$ or $\Ypsilon'{\,\cong\,}\Spin7\RR$. 
In the first case $\dim\Ypsilon{\,\le\,}17$ and $\dim\Tau_{\hskip-1pt[v]}{\,=\,}8$, in the second 
case  $\dim\Ypsilon{\,\le\,}22$ and $\dim\Tau_{\hskip-1pt[v]}{\,\in\,}\{3,4\}$, but then 
$\Tau_{\hskip-1pt[v]}{\,\le\,}\Cs{}\Phi$ and $\cF_\Lambda{\,\ledot\hskip-3pt}\cP$, which 
contradicts  Stiffness. \\
(i) {\it If $t{\,=\,}2$, then the translation group  $\Tau$ is transitive\/}. \\ Because of steps (f--h), only 
the following situation must be considered: $\dim\nabla{\,=\,}17$, for some $w{\,\in\,}S$ the orbit 
$w^\nabla$  has dimension $\nabla{:}\nabla_{\hskip-2pt w}{\,<\,}8$,\:  
$\Omega{\,=\,}(\nabla_{\hskip-2.5pt w})^1$ satisfies $\dim\Omega{\,=\,}10$, and 
$\Lambda{\,=\,}(\Omega_c)^1{\,\cong\,}\SU3\CC$ for $a{\,\ne\,}c{\,\in\,}a^\Theta$. Put 
$\Rho{\,=\,}\sqrt\Omega$, so that $\Omega{\,=\,}\Rho\Lambda$. As $\Lambda|_\Theta{\,=\,}\1$, 
the radical $\Rho$ is transitive on $\Theta\sm\{\1\}$,\: 
$\Rho/\Rho_{\hskip-2pt c}{\,\approx\,}\CC^{\times}$,\: 
$\dim\Rho_{\hskip-2pt c}{\,=\,}0$, and $\Rho_{\hskip-2pt c}$ is compact by 2.6(\^c). 
Therefore $\Rho_{\hskip-2pt c}$ is finite and $\Rho$ is a finite covering of the cylinder group. 
Hence $\Rho$ contains an involution $\sigma$, and $\sigma$ is a reflection (or else $\Lambda$ 
would act trivially on  $W\smcap\cF_\sigma$), the axis of $\sigma$ is $W$\: (since 
$w^\sigma{\,=\,}w$), and $\sigma{\,\in\,}\Delta_{[a]}$. Now 2.9 and step (g) imply 
$\dim\Tau{\,=\,}16$. \\  
(j) {\it $t{\,\ne\,}3$\/}. Assume the contrary.  Choose $w{\,\in\,}S$ and consider the plane 
$\cB{\,=\,}\langle a^\Theta,u,w\rangle$ of dimension ${\ge\,}8$. Suppose first that $\cB{\;\ldot}\cP$. 
By Stiffness, the connected component $\Omega$ of $\nabla_{\hskip-2pt w}$ induces on $\Theta$ a subgroup  $\overline\Omega{\,\le\,}\GL3\RR$ of dimension at least $6$, and there is a subgroup 
$\Xi{\,\le\,}\Omega$ such that the fixed elements of $\Xi$ in $\Theta$ form a $2$-dimensional subgroup\footnote{\ The group $\overline\Omega$ fixes a proper subspace of $\Theta$ or 
$\overline\Omega$ is irreducible and contains $\SL3\RR$. In any case, the stabilizer $\Chi$ in 
$\overline\Omega$ of some $2$-dimensional subspace $\Eta{\,<\,}\Theta$ has dimension ${>}4$, but 
$\dim\Chi|_\Eta{\,\le\,}4$.}. It follows that 
$\cB\smcap\cF_\Xi{\hskip4pt\ldot}\cB$. For a suitable $c{\,\in\,}a^\Theta\sm\{a\}$ the group 
$\Lambda{\,=\,}(\nabla_{\hskip-2pt c,w})^1$ is compact by \cite{cp} 83.9. Moreover, 
$\dim\Lambda{\,\ge\,}6$, and 2.6(\^e) shows $\Lambda{\,\cong\,}\SO4\RR$. Now 
$\Lambda|_\cB{\,\cong\,}\Lambda|_\Theta{\,\cong\,}\SO3\RR$, but then $c^\Lambda{\,\ne\,}c$. 
Therefore $\cB{\,=\,}\cP$ and $\nabla_{\hskip-2pt w}$ acts faithfully on $\Theta$. Consequently,  
the  connected components of  $\nabla_{\hskip-2pt w}$ and $\GL3\RR$ are isomorphic, and 
$\dim\Lambda{\,=\,}6$.  As $\SO3\RR$ is a maximal subgroup of $\SL3\RR$ (see \cite{cp} 94.34), 
it follows that a maximal compact subgroup of $\Lambda$ is trivial or a torus; hence
$\pi_q\Lambda{\,=\,}0$ for $q{\,>\,}1$. 
Moreover, $\dim\nabla_{\hskip-2pt c}{\,=\,}14$ and $\nabla_{\hskip-2pt c}$ is transitive on $S$. 
Let $\Phi$ be a maximal compact subgroup of $\nabla_{\hskip-2pt c}$. The exact homotopy 
sequence  (cf. \cite{cp} 96.12) implies $\pi_3\Phi{\,\cong\,}\pi_3\Lambda{\,=\,}0$ and 
$\pi_7\Phi{\,\cong\,}\pi_7S{\,\cong\,}\ZZ$. From the first condition it follows by \cite{cp} 94.36 that 
$\Phi$ has no almost simple factor, so that $\Phi$ is a torus group (\cite{cp} 94.31(c)\,), but then 
$\Phi$ cannot satisfy the second condition.   \\
(k) {\tt Remark.} The arguments in step (j) do not use the fact that $\Theta^\nabla{\,=\,}\Theta$. 
They show therefore that there is no  $\Omega$-invariant $3$-dimensional 
subgroup of $\Theta$. Similarly, only the $\Omega$-invariance is used in 
steps  ($\ell$), (m), and (n). \\ 
($\ell$) {\it If $t{\,=\,}4$ and if $\cB{\,=\,}\langle a^\Theta,u,w\rangle{\:\ldot}\cP$, then $\Tau$ is
transitive\/}. \\
Assume first that $\Omega{\,=\,}(\nabla_{\hskip-2pt w})^1$ acts irreducibly on $\Theta$, and  put 
$\Omega^*{\,=\,}\Omega|_\Theta{\,=\,}\Omega/\Kappa$. Stiffness implies $\dim\Kappa{\,=\,}3$ 
or $\dim\Kappa{\,\le\,}1$. Therefore $\dim\Omega^*{\,\ge\,}6$ and even 
$\dim\Omega^*{'}{\,\ge\,}6$ (see \cite{cp} 95.6). It follows that 
$\Omega^*{'}$ is isomorphic to one of the groups $\Sp4\RR$,\: $\SL2\CC$, or 
$\Opr4(\RR,r)$, taken in its standard action on $\Theta$. In each case there is a subset 
$\Xi{\,\subset\,}\Omega$ such that the fixed point set $\Theta\smcap\Cs{}\Xi$ is isomorphic to 
$\RR^2$\: (this is obvious for the complex group; in the other cases $\Xi$ may be chosen as preimage of a $1$-dimensional torus). Now $\cF_\Xi\smcap\cB{\;\ldot}\cB$,\: 
$\Lambda{\,=\,}(\Omega_c)^1{\,\cong\,}\SO4\RR$ by 2.6(c,\^e), and $\dim\Omega^*{\,\le\,}7$. 
Hence $\Kappa^1{\,\cong\,}\Spin3\RR$. As $\Lambda$ has a subgroup $\SO3\RR$, there is a group 
$\Phi{\,\cong\,}\ZZ_2^{\hskip2pt4}$ in $\Omega$, and $\Phi$ contains at least one reflection 
$\sigma$  by 2.10. In fact, 
$\sigma{\,\in\,}\Delta_{[a,W]}$, since $\sigma{\,\in\,}\nabla_{\hskip-2pt w}$.
We have $\dim\Omega{\,\le\,}10$ and $\dim a^\Delta{\,=\,}\Delta{:}\nabla{\,\ge\,}15$. Finally 2.9 
shows that $\sigma^\Delta\sigma{\,=\,}\Tau{\,\cong\,}\RR^{16}$. \\Hence we may suppose that 
there is a proper $\Omega$-invariant subgroup $\Eta{\,<\,}\Theta$, and $\dim\Eta{\,<\,}3$ by 
Remark (k).  Again there exists a subplane $\cF{\,\ldot}\cB$, and for some 
$c{\,\in\,}a^\Eta\sm\{a\}$ Stiffness implies $\Lambda{\,\cong\,}\SU3\CC$, but then 
$\Lambda|_\Theta{\,=\,}\1$, which is impossible. \\
(m)  {\it If $t{\,=\,}4$ and if $\cB{\,=\,}\langle a^\Theta,u,w\rangle{\,=\,}\cP$, then $\Tau$ is transitive\/}. \\
$\cB{\,=\,}\cP$ implies that  $\Omega$ acts effectively on $\Theta$. 
If this action is irreducible, then $\Omega'$ is isomorphic to $\Sp4\RR$ or to $\SL4\RR$. In both 
cases the central involution $\omega{\,\in\,}\Omega'$ is a reflection, since otherwise
$\dim\Omega'|_{\cF_\omega}{\,\le\,}7$ by \cite{cp}  83.17.   Again $\omega{\,\in\,}\Delta_{[a,W]}$ 
and $\Tau{\,=\,}\omega^\Delta\omega$ has even dimension $\Delta{:}\nabla$, see 2.9. 
If $\dim\Omega'{\,=\,}10$, then $\dim\nabla{\,\le\,}19$ and $\dim\Tau{\,\ge\,}14$,\: 
$\dim\Tau_{\hskip-1pt[w]}{\,\ge\,}6$. As $\omega$ inverts each $\tau{\,\in\,}\Tau$, the group 
$\Omega'$ acts faithfully on $\Tau_{\hskip-1pt[w]}$, and then \cite{cp} 95.10 shows 
$\dim\Tau_{\hskip-1pt[w]}{\,=\,}8$ and $\Tau{\,\cong\,}\RR^{16}$. 
Now let $\Omega'{\,\cong\,}\SL4\RR$. Then $\Omega$ has a subgroup $\Phi{\,\cong\,}\SO3\RR$, 
and $\Phi$ fixes some $c{\,\in\,}a^\Theta\sm\{a\}$, so that $\Phi{\,<\,}\Lambda{\,=\,}\Omega_c^{\,1}$. 
Stiffness properties 2.6(\^c,\^e) imply $\Lambda{\,\cong\,}\SO4\RR$, but $\dim\Lambda{\,\ge\,}11$. 
Suppose now that there is an $\Omega$-invariant proper subgroup $\Eta$ of $\Theta$. 
From Remark (k) it follows that $\dim\Eta{\,\le\,}2$. Choose $c{\,\in\,}a^\Eta\sm\{a\}$. Then
$\dim\Lambda{\,\ge\,}9-\dim\Eta$. We will show that $\Lambda$ is compact. 
If $\Eta{\,\cong\,}\RR$ and $\Lambda$ is transitive on $\Theta\sm\Eta$, then the linear group 
$\Lambda$ is also transitive on $\Theta/\Eta{\,\cong\,}\RR^3$; hence $\Lambda$
contains $\SO3\RR$ and $\Lambda{\,\cong\,}\SU3\CC$ by 2.6(\^c,\^e). 
In all other cases,  there are some points  $d{\,\in\,}a^\Theta\sm\{a\}$ and 
$d'{\,\in\,}a^\Theta\sm a^\Eta$ such that  
$\dim\Lambda_d{\,\ge\,}5$,\: $d'{\,\notin\,}\cF_{\Lambda_d}$,  and $\Lambda_{d,d'}{\,\ne\,}\1$. 
Therefore  $\cF_{\Lambda_d}{\,\ldot}\cF_{\Lambda_{d,d'}}{\,\ldot}\cP$. Again 
$\Lambda{\,\cong\,}\SU3\CC$ by Stiffness, but then the action of $\Lambda$ would be trivial.    \\
(n) {\it If $t{\,=\,}5$, then $\Tau$ is transitive\/}. In fact,  the group $\Omega$ acts effectively on 
$\Theta$, since  $\langle a^\Theta,u,w\rangle{\,=\,}\cP$. 
If this action is irreducible, then \cite{cp} 95.10 shows that $\Omega'{\,\cong\,}\Opr5(\RR,r)$; the 
case $r{\,=\,}0$ is excluded by 2.10. Consider a subgroup $\Phi{\,\cong\,}\SO3\RR$ in $\Omega$, 
a point $c{\,\in\,}a^\Theta\sm\{a\}$ such that $c^\Phi{\,=\,}c$,  and $\Lambda{\,=\,}(\Omega_c)^1$. 
Stiffness 2.6(\^c,\^e) implies $\Lambda{\,\cong\,}\SO4\RR$\: (note that each action of $\SU3\CC$ 
on  $\Theta$ is trivial). Consequently $r{\,=\,}1$, and $\Omega$ has a subgroup 
$\Chi{\,\cong\,}\SO4(\RR,1)$. The group  $\Chi$ fixes some point  $d{\,\in\,}a^\Theta\sm\{a\}$.
Now  $\Chi$ is compact by Stiffness 2.6(\^c), but this is not true. Hence $\Theta$ has a minimal 
$\Omega$-invariant subgroup $\Eta{\,\cong\,}\RR^s,\: s{\le}2$ or $s{\,=\,}4$. In the first case, let 
$c{\,\in\,}a^\Eta\sm\{a\}$. Because any action of $\SU3\CC$ on $\Theta$  is trivial, stiffness property 2.6(e) implies $\dim\Omega_c{\,=\,}7$ and $s{\,=\,}2$. It follows that 
$\dim\nabla_{\hskip-2pt c}{\,=\,}15$, and $\dim\Omega{\,=\,}9$. Consider the subplane 
$\cE{\,=\,}\langle a^\Eta,u,w\rangle{\,=\,}\cE^\Omega$. If $\cE{\ledot\hskip-3pt}\cP$ and 
$\Omega|_\cE{\,=\,}\Omega/\Kappa$, then $\dim\Kappa{\,\le\,}3$ and 
$\Omega{\hskip1pt:\hskip1pt}\Kappa{\,\ge\,}6$; on the other hand, 
$\Omega|_\cE{\,\cong\,}\Omega|_\Eta$ is a 
subgroup of the $4$-dimensional group $\GL2\RR$. Therefore $\dim\cE{\,=\,}4$ and 
$\dim\Omega|_\cE{\,\le\,}2$\: (see \cite{cp} 71.7). Hence $\dim\Kappa{\,\ge\,}7$ and 
$\Kappa^1{\,=\,}\Lambda$. Choose any $d{\,\in\,}a^\Theta\sm a^\Eta$. Then 
$\dim\Lambda_d{\,\ge\,}2$,\:  $\langle \cF_\Lambda,d\rangle{\,\ldot}\cP$, and
Stiffness implies that $\Lambda$ is compact, but then the action of $\Lambda$ on $\Theta$ is
completely reducible, and $\Lambda$ acts effectively on a  complement $\RR^3$ of $\Eta$ in 
$\Theta$, which is impossible. Thus $s{\,=\,}4$, and $\Tau$ is transitive by the previous steps. \\
(o) {\it If $t{\,=\,}6$ and $\Omega$ acts irreducibly on $\Theta$, then $\Tau$ is transitive\/}. {\tt Proof:}
$\Omega'$ is semi-simple by \cite{cp} 95.6, and $\dim\Omega'{\,>\,}6$. Clifford's Lemma implies that $\Omega'$ is almost simple or isomorphic to a direct product $\SL2\RR{\times}\SL3\RR$. 
In the second case any involution $\beta$ in the larger factor $\Beta$ is planar, and 
$\Gamma{:=\,}\Cs\Beta\beta{\,\cong\,}\SL2\RR$. The fixed elements of $\Gamma$ on $\Theta$ 
form a $2$-dimensional subgroup. As $\Gamma$ is not compact, $\dim\cF_\Gamma{\,=\,}4$ and
$\cF_\Gamma{\;\ldot}\cF_\beta{\;\ldot}\cP$, but then $\Gamma$ would be compact by 2.6(c).
Now it follows from \cite{cp} 95.10 that $\dim\Omega'{\,=\,}8$ or $\dim\Omega'{\,\ge\,}15$. 
In the latter case $9{\,\le\,}\dim\Omega_c{\,\ne\,}14$ in contradiction to 2.6(e). 
Hence $9{\,\le\,}\dim\Omega{\,\le\,}10$. If $\dim\Omega{\,=\,}10$, then $\dim a^\Delta{\,\ge\,}15$,\:
$\Omega{\,\cong\,}\CC^{\times}{\cdot\hskip1pt}\SU3(\CC,r)$,\:  $\Omega$ contains a central 
involution $\alpha{\,\in\,}\nabla_{\hskip-1pt[a]}$, and 2.9 shows that $\Tau{\,\cong\,}\RR^{16}$. 
If $\dim\Omega{\,=\,}9$ for each choice of $w$, then $\nabla$ is transitive on $S$,\: 
$\dim\nabla{\,=\,}17$,\; $a^\Delta$ is open in $P$, a maximal compact subgroup $\Phi$ of 
$\Omega$ satisfies $\dim\Phi'{\,\in\,}\{3,8\}$,  and 2.17 implies that $\nabla$ must have a  compact subgroup $\U2\HH$ or $\SU4\CC$. These  possibilities will be discussed in the next two  steps. \\
(p) Suppose that $\nabla$ contains $\Ypsilon{\,\cong\,}\U2\HH$, but not $\SU4\CC$. The central involution $\sigma$ of $\Upsilon$ acts trivially on $\Theta$ and  $\Upsilon|_\Theta{\,\cong\,}\SO5\RR$, 
moreover, $\sigma$ is a reflection with axis $av$, and $\Tau_{\hskip-1pt[u]}$ has positive dimension. 
In fact, $\Tau_{\hskip-1pt[u]}{\,\cong\,}\RR^8$ because the representation of $\Upsilon$ on 
$\Tau_{\hskip-1pt[u]}$ is faithful or since $\dim(av)^\Delta{\,=\,}8$. Note that 
$\nabla\smcap\Cs{}\Tau_{\hskip-1pt[u]}{\,\le\,}\Delta_{[v,au]}$. 
If $\Delta_{[v,au]}{\,\ne\,}\1$, then $\Tau_{\hskip-1pt[v]}$ is transitive because 
$\dim(au)^\Delta{\,=\,}8$, and $\Tau{\,\cong\,}\RR^{16}$. Thus we may assume that $\nabla$ acts 
effectively on $\Tau_{\hskip-1pt[u]}$; as $\Ypsilon{<\,}\nabla$, the action is also irreducible. Hence 
the  commutator group $\nabla'$ is semi-simple, and then \cite{cp} 94.33 shows 
$\nabla'{\,\cong\,}\SL2\HH$. Recall that $\Omega'$ is an almost simple group of dimension~$8$, 
and note that $\Omega'$ acts effectively on $\Tau_{\hskip-1pt[u]}{\,\cong\,}\RR^8$.
Let $\Chi{\,=\,}\Tau_{\hskip-1pt[u]}\smcap\Cs{}\Omega'$. Either $\Chi$ is trivial and $\Omega'$ 
acts irreducibly on $\Tau_{\hskip-1pt[u]}$, or $\Chi{\,\cong\,}\RR^2$ by complete reducibiliy and 
\cite{cp} 95.10,\; $\cC{\,=\,}\langle a^\Chi,v,w\rangle$ is a $4$-dimensional subplane, 
$\Omega'|_\cC{\,=\,}\1$, and $\Omega'{\,\cong\,}\SU3\CC$ by Stiffness.
In the first case the centralizer of $\Omega'|_{\Tau_{\hskip-1pt[u]}}$  is  $\RR^{\times}$\: (see 
\cite{cp} 95.10), but the centralizer of the action of $\nabla$ on $\Tau_{\hskip-1pt[u]}$ is 
$\HH^{\times}$,  a contradiction.  In the second case
$\SU3\CC$ would be contained in $\U2\HH$, which is not true. \\
(q) Thus $\nabla'{\,\cong\,}\SU4\CC$ and $\nabla'|_\Theta{\,\cong\,}\SO6\RR$. Again the 
central involution in $\nabla'$ is a reflection in $\nabla_{\hskip-1.5pt[u]}$ and 
$\Tau_{\hskip-1pt[u]}{\,\cong\,}\RR^8$. Let $\Gamma{\,=\,}\Delta_{av}$ and consider the action 
of $\nabla'$ on (the additive group of) the Lie algebra $\frl\hskip1pt\Gamma$. 
As $\Theta^{\hskip-1pt\nabla}{\,=\,}\Theta$ and $\dim\Gamma{\,=\,}25$, the group 
$\Xi{\,=\,}\Gamma\smcap\Cs{}\nabla'$ is $4$-dimensional. Put 
$\Xi\smcap\Cs{}\Tau_{\hskip-1pt[u]}{\,=\,}\Chi$, and note that $\Chi^{\nabla'}{\,=\,}\Chi$ and that
$\Xi/\Chi{\,\le\,}\CC^{\times}$\; (see \cite{cp} 95.10). 
By definition, $\Chi{\,\le\,}\Delta_{[v]}$ and $\dim\Chi{\,\ge\,}2$. 
Either $\Chi{\,\le\,}\Tau_{\hskip-1pt[v]}$, or $\Chi$ contains a homology and \cite{cp} 61.20 
implies $\Tau_{\hskip-1pt[v]}{\,\cong\,}\RR^8$. In the first case, the action of $\nabla'$ on $\Theta$ shows that $\Chi\smcap\Theta{\,=\,}\1$.
Hence $\Chi\Theta{\,=\,}\Tau_{\hskip-1pt[v]}{\,\cong\,}\RR^8$ and $\Tau$ is transitive.\\
(r) If $t{\,=\,}6$ and $\Tau$ is not transitive, then steps (k--o) imply that there exists a minimal 
$\Omega$-invariant subgroup $\Eta{\,<\,}\Theta$ of dimension ${\,\le\,}2$.  As in step (n),
let $a{\,\ne\,}c{\,\in\,}a^\Eta$. Then representation on $\Theta$ shows that 
$\Lambda{\,=\,}(\Omega_c)^1{\,\not\cong\,}\SU3\CC$, and Stiffness 2.6(e) implies 
$\dim\Lambda{\,=\,}7$. Hence $\Eta{\,\cong\,}\RR^2$,\; $\dim\nabla_{\hskip-2pt c}{\,=\,}15$, and 
$\Xi{\,=\,}(\nabla_{\hskip-2pt c})^1$ is transitive and faithful on $S$.
By 2.17(\^c) it follows that $\Xi$ has a compact subgroup $\Psi{\,\cong\,}\U2\HH$, but $\Xi$ is not
compact.  The central involution $\sigma{\,\in\,}\Psi$ is a reflection with axis $av$. The dual of 2.9 shows $\dim\Tau_{\hskip-1pt[u]}{\,>\,}0$, and then $\Tau_{\hskip-1pt[u]}{\,\cong\,}\RR^8$ since 
$\Psi$ acts faithfully on $\Tau_{\hskip-1pt[u]}$. Consequently, $\Xi$ is an irreducible subgroup 
of $\GL8\RR$ and $\Xi'$ is semi-simple. \cite{cp} 94.33 and 95.10 imply
$\Xi{\,\cong\,}\SL2\HH$,  but $c^{\hskip.5pt\Xi}{\,=\,}c$ and each  representation 
of $\SL2\HH$ on  $\Theta{\,\cong\,}\RR^6$ with a fixed point is trivial. This is impossible. \\  
(s) {\it If $\dim\Delta{\,=\,}33$ and  $t{\,=\,}7$, then $\Tau$ is transitive\/}:\: 
The action of $\nabla$ on $\Theta$ is irreducible, because $\Theta$ is a minimal 
$\nabla$-invariant group, and $\dim\nabla{\,<\,}28$, since $\dim\nabla\Theta{\,\le\,}33$.  
By \cite{cp} 95.\:5,\:6, and 10, the group $(\nabla|_\Theta)'$ is isomorphic to $\Gtwo$ or 
locally isomorphic to $\Spin7(\RR,r),\: r{\le}3$, and \cite{cp} 94.27 shows that 
$(\nabla|_\Theta)'$ is covered by a subgroup $\Upsilon{\,\le\,}\nabla$. 
In any case it will turn out that $\nabla$ contains a group $\Gamma{\,\cong\,}\Gtwo$. 
Suppose that $\dim\Upsilon{\,=\,}21$. Then $\dim\Omega{\,\ge\,}13$. Let $\Eta{\,\cong\,}\RR^s$
be a minimal $\Omega$-invariant subgroup  of~$\Theta$. According to {Remark} (k), 
 we have $s{\,\le\,}2$ or $s{\,\ge\,}6$. 
Choose $c{\,\in\,}a^\Eta\sm\{a\}$ and put  $\Lambda{\,=\,}(\Omega_c)^1$. In the first case 
$\dim\Lambda{\,\ge\,}11$, and Stiffness 2.6(d) implies $\Lambda{\,\cong\,}\Gtwo$. 
If $\dim\Eta{\,=\,}6$, then $\Omega$ acts effectively and irreducibly on $\Eta$, and 
$\Omega'$ is almost simple by Clifford's Lemma \cite{cp} 95.5. The list \cite{cp} 95.10 of 
representations shows that $\dim\Omega'{\,\ge\,}15$. Consequently $\dim\Lambda{\,>\,}8$. 
By 2.6(e) again $\Lambda{\,\cong\,}\Gtwo$, and $\nabla$ has indeed a subgroup $\Gamma$
as claimed. The representation of $\Gamma$ on (the additive group of) the Lie algebra 
$\frl\hskip1pt\Delta$ is completely reducible, and $\Delta{\,:\,}\Gamma\Theta{\,=\,}12$. 
The only representation of $\Gtwo$ in dimension ${<\,}14$ is the natural one as  
$\Aut\OO$ on the space of pure octonions. Therefore $\Rho{\,=\,}(\Cs{}\Gamma)^1$ 
satisfies $\dim\Rho{\,=\,}5$, but $\Rho$ acts faithfully on the flat plane $\cF_\Gamma$,  
fixes $u$ and $v$, and hence has dimension ${\le\,}4$, a contradiction. \\
(t) {\bf Corollary.} {\it If $\dim\Delta{\,=\,}33$ and $z{\,\in\,}\{u,v\}$, then 
$\dim\Tau_{\hskip-1pt[z]}{\,=\,}0$ or  $\Tau_{\hskip-1pt[z]}{\,\cong\,}\RR^8$\/}. \\
(u) {\it If $\dim\Delta{\,=\,}33$ and if $\dim\Delta_{[u]}{\,>\,}0$ or $\nabla_{\hskip-1pt[u]}{\,\ne\,}\1$,  
then $\Tau$ is transitive\/}. \\
This is an immediate consequence of \cite{cp} 61.20, the assumption that $W$ is the only fixed line,  and Corollary (t). \\ 
Only the case that $\nabla$ acts faithfully and irreducibly on $\Theta{\,\cong\,}\RR^8$ has still to be considered. Put $\Gamma{\,=\,}\Delta_{au}$ and note that $\nabla{\,\le\,}\Gamma$ 
and $\dim\Gamma{\,=\,}25$. \\
(v)  {\it $\Gamma{:\,}\Gamma_{\hskip-1pt[u]}{\,\le\,}22$,\: 
$\dim\Gamma_{\hskip-1pt[u]}{\,>\,}0$, and $\Tau$ is transitive\/}. \\
The group $\Upsilon{=\,}\Gamma|_\Theta{\,=\,}\Gamma/\Gamma_{\hskip-1pt[u]}$ is an 
irreducible subgroup of $\GL8\RR$. Suppose that 
$\Gamma{:\,}\Gamma_{\hskip-1pt[u]}{\,>\,}20$. Then $\dim\Upsilon'{\,>\,}18$, and Clifford's Lemma \cite{cp} 95.6 implies that $\Upsilon'$ is almost simple and irreducible (cf. \cite{HS} proof, step 18)   for details). The list \cite{cp} 95.10 of 
representations shows $\Upsilon'{\hskip1pt\cong\,}\Spin7(\RR,r)$ with $r{\,=\,}0,3$; in particular, 
$\dim\Upsilon'{\,=\,}21$ and $\dim\Upsilon{\,\le\,}22$. \qed
\par\medskip
The planes  of 6.5 with $\dim\Delta{\,>\,}34$ have been described  in \cite{HS}:
\par\medskip
{\bf 6.7 Cartesian planes.} {\it A plane $\cP$ satisfies the conditions of {\rm 6.5} with 
$\dim\Delta{\,\ge\,}35$ if, and only if, $\cP$ can be coordinatized by a topological Cartesian field 
$(\OO,+,\diamond)$ defined as follows\/}: {\it let $(\RR,+,\ast,1)$ be an arbitrary real Cartesian 
field  with unit element $1$ such that identically $(-r){\ast}s{\,=\,}-(r{\ast}s)$, and let 
$\rho{\,:\,}[0,\infty){\,\approx\,}[0,\infty)$ be any homeomorphism with $\rho(1){\,=\,}1$. 
Write each octonion in the form $x{\,=\,}x_0{\,+\,}\frx$ as in $2.5$, and put 
$$s\diamond x{\,=\,}|s|^{-1}s\hskip1.5pt(|s|{\ast\hskip1pt}x_0+\rho(|s|){\hskip1pt\cdot\hskip1pt}\frx) \ 
{\rm for} \ s\ne 0, \ 0\diamond x{\,=\,}0\,.$$\/}\par\smallskip
Recall that $\Sigma{\,=\,}\Aut\cP$.\quad (a) {\it If $\dim\Sigma = 39$, then $\cP$ is a translation plane\/}. \par
(b) {\it The plane $\cP$ is a translation plane if\/, and only if\/, it can be coordinatized  by a quasi-field $\OO_{\textstyle\diamond}$  where $*$ is the ordinary multiplication of the 
reals. In this case $\dim\Sigma=39$ if\/, and only if\/, $\rho$ is a multiplicative homomorphism$;$ otherwise $\dim\Sigma=38$\/}.\par
(c) {\it If $\cP$ is not a translation plane, then
$\dim\Sigma{\,=\,}38$ if, and only if, $\cP$ can be coordinatized by a Cartesian field 
$\OO_{\textstyle\diamond}$  where 
$r{\ast}s{\,=\,}rs\hskip6pt(s{\,\ge\,}0)$ and $r{*}s{\,=\,}|r|^{\gamma\,} rs\hskip6pt(s{\,<\,}0)$ for some 
$\gamma{\,>\,}0$  and $\rho:[0,\infty) \to [0,\infty)$ is a multiplicative homomorphism\/}. \par
For the cases $\dim\Sigma{\,=\,}37$ the reader is referred to the last part of \cite{HS}.
\par\medskip

Proofs will not be repeated here. There seems to be little chance to improve these results.
\par\bigskip 
\newpage
{\Bf 7. Collinear fixed points}
\par\medskip
Assume in this section that $\Delta$ fixes more than $2$ points but only one line. Some results will be stronger than in the previous section. Write   $\cF_\Delta{\,=\,}\{u,v,w,...,W\}$.  
\par\medskip
{\bf 7.0 Lie.} {\it If $\Delta$ fixes at least $3$ distinct points $u,v,w$ and exactly one line $W$, 
and if $\dim\Delta{\,\ge\,}15$, then $\Delta$ is a Lie group\/}. 
\par\smallskip
(a) The {\tt proof} is similar to that of 6.0, and the same notation will be used. 
Assume that $\Nu$ is not a Lie group. 
If $\cD{\,=\,}\langle x^\Delta\rangle$ is flat, then 2.6(\^a) implies $\dim\Delta{\,\le\,}3{+}11$, 
otherwise $\Kappa{\,\cong\,}\Gtwo$ would be a maximal compact subgroup. 
If $\dim\cD{\,=\,}4$, then $\dim\Delta{\,\le\,}6{+}8$, again by Stiffness.
Hence $\cD{\;\ldot}\cP$, the  kernel $\Kappa$ is compact, and $\dim\Kappa{\,\le\,}7$. \\
(b) Choose a line $L{\,\in\,}\cD, \ v{\,\in\,}L$,   a point  $z{\hskip1pt\in\hskip1pt}L\sm\cD$,  and let 
$\Lambda{\,=\,}(\Delta_z)^1$. According to Lemma 7.0$'$ 
below, $\Delta/\Kappa$ is a Lie group, and so is $\Delta/(\Kappa\smcap\Nu)$. Therefore
we  may assume that $\Nu{\,\le\,}\Kappa$. We have  $\Kappa_z{\,=\,}\1$\; because 
$\langle \cD,z\rangle{\,=\,}\cP$. Hence $\Lambda\smcap\Nu{\,=\,}\1$ and $\Lambda$ 
is a Lie group. Moreover, $\Delta_L{\hskip1pt:\hskip1pt}\Lambda{\,<\,}8$, or $\Kappa$ would be 
a Lie group by \cite{cp} 53.2. Consequently $\dim\Lambda{\,\ge\,}4$. \\
(c) If  $\dim\Kappa{\,=\,}0$, then $\Delta{\hskip1pt:\hskip1pt}\Kappa{\,=\,}15$ 
and $\cD$ is the classical quaternion plane 
(see \cite{sz1} 5.\,1,3, and~5). In this case,  
$\Delta$ fixes a connected subset of $uv\smcap\cD$ pointwise, and $\cF_\Lambda$ is a 
connected subplane  of dimension at most $4$. By \cite{cp} 55.4 the plane 
$\langle z^\Nu,u,w\rangle{\,\le\,}\cF_\Lambda$ 
is connected, and then $\Nu$ is a Lie group by \cite{cp} 32.21  and 71.2. \\
(d)  As in the proof of 6.0, 
the connected component  $\Kappa^1$ is a product of a commutative connected group $\Alpha$ with at most one almost simple factor $\Omega{\,\cong\,}\Spin3\RR$, and $\Alpha$ is in the centralizer 
of $\Lambda$. Either $\Kappa^1{\,\cong\,}\Omega$ or $z{\,\ne\,}z^\Alpha{\,\subseteq\,}\cF_\Lambda$
and $\cF_\Lambda$ is connected. In the second case  $\Nu$ would again be a Lie group. \\
(e) If $\Kappa^1{\,\cong\,}\Omega$, then 
$\Lambda/\Cs\Lambda\Omega{\,\le\,}\Aut\Omega{\,\cong\,}\SO3\RR$ and 
$\dim\Cs\Lambda\Omega{\,\ge\,}1$; in other words, $\Lambda$ has a connected subgroup $\Rho$ which centralizes $\Omega$ and 
$\Kappa{\,\approx\,}z^\Kappa{\,\subseteq\,}\cF_\Rho{\,\ge\,}\cF_\Lambda$. Therefore 
$\cF_\Rho{\;\ldot}\cP$ and $\Rho$ is compact. In particuar, there is an involution 
$\rho{\,\in\,}\Rho$ such that $\cF_{\hskip-2pt\rho}{\,=\,}\cF_\Rho$. As $z^\rho{\,=\,}z$, the 
induced map $\overline\rho{\,=\,}\rho|_\cD$ is either a Baer involution or a reflection. \\
(f) If $\overline\rho$ is planar, i.e., if $\cD\smcap\cF_{\hskip-2pt\rho}{\;\ldot}\cD,\cF_{\hskip-2pt\rho}$, then the lines of $\cD$ are $4$-spheres (cf. \cite{cp} 53.10 or \cite{sz3} 3.7) and so are the lines of $\cF_{\hskip-2pt\rho}$. Richarson's theorem $(\dagger)$ 
implies that  the action of $\Omega$ on $L\smcap\cF_{\hskip-2pt\rho}$ is equivalent to the 
standard action  of $\SU2\CC$ on $\sS_4$. Hence $\Omega$ cannot fix the $2$-sphere 
$L\smcap\cF_{\overline\rho}$. This contradiction shows that $\overline\rho$ is a reflection with 
axis $W\smcap\cD$ and some center $a$. \\
(g) In the latter case, the lines of $\cD$ are again homeomorphic to $\sS_4$; this follows from the action of $\Delta$ on $\cD$: if $L\smcap\cD$ is not a manifold, then 
$\Delta{\hskip1pt:\hskip1pt}\Delta_L{\,<\,}4$ by \cite{cp} 53.2 \ (cf. also \cite{HK} 5.5), and 
$\Delta_L{\hskip1pt:\hskip1pt}\Kappa{\,\ge\,}9$.  On the other hand,  
$\Delta_a{\hskip1pt:\hskip1pt}\Kappa{\,<\,}7$, or $L\smcap\cD$ would be a manifold by
 \cite{sz15}~$(\ast\ast)$. Therefore $\dim a^{\Delta_L}{\,=\,}3$, 
and $\overline\rho^{\Delta_L}\overline\rho$ is a $3$-dimensional translation group $\Theta$ 
of $\cD$, but $\dim\Theta$ is even by 2.9. \\
(h)  The lines of $\cK{\,=\,}\cF_{\hskip-2pt\rho}$ are also $4$-spheres, in fact,  
$W\smcap\cD{\,=\,}W\smcap\cK$ and $S{\,=\,}L\smcap\cK{\,\approx\,}\sS_4$. This is a consequence of \cite{cp} 92.16 and 51.20, see the  proof of 4.0. \\
(i) Now Richarson's theorem $(\dagger)$ can be applied to the action of $\Omega$ on $S$. 
It follows that $S\smcap\cD{\,=\,}\{a,v\}$ and that $\Omega$ acts freely on $S\sm\{a,v\}$. The orbit space $S/\Omega$ is an interval and the compact group $\Nu$ in $\Cs{}\Omega$ induces the 
identity on this interval. Hence $z^{\Omega\Nu}{\,=\,}z^\Omega$, but $(\Omega\Nu)_z{\,=\,}\1$. 
This contradiction completes the proof. \qed
\par\medskip
{\bf 7.0$'$ Lemma.} {\it Suppose that $\Delta$ is a group of automorphisms of an 
$8${\rm -dimensional} plane~$\cP$. If $\Delta$ fixes at least $3$ distinct points $u,v,w$ and exactly one line $W$, and if $\dim\Delta{\,\ge\,}8$, then $\Delta$ is a Lie group\/}.
\par\smallskip
The {\tt proof} is similar to the previous one but less involved. If $\Delta/\Nu$ is a Lie group, 
but the compact central subgroup $\Nu$ is not, then  the arguments of 4.0 show that 
for some $x{\,\notin\,}W$ the orbit $x^\Delta$ generates a proper connected subplane. 
Put  $\cD{\,=\,}\langle x^\Delta,u,v,w\rangle$ and 
$\Delta|_\cD{\,=\,}\Delta/\Kappa$, and apply the stiffness results 
\cite{sz1} 1.5 to $\Kappa$. 
If $\cD$ is flat, then $\Delta{\hskip1pt:\hskip1pt}\Kappa{\,\le\,}3$ and $\dim\Kappa{\,\le\,}3$ \ (see 
also \cite{cp} 83.12); if $\cD{\;\ldot}\cP$, then $\Delta{\hskip1pt:\hskip1pt}\Kappa{\,\le\,}6$ and 
$\dim\Kappa{\,\le\,}1$ \ (\cite{cp} 83.11). \qed
\par\medskip
 {\bf 7.1.} {\it If $\Delta$ is semi-simple of dimension $\dim\Delta{\,\ge\,}14$, then $\Delta$ is a 
Lie group\/}, cf. 6.1.
\par\medskip
{\bf 7.2 Semi-simple groups.} {\it If $\Delta$ is semi-simple with more than $2$ collinear fixed 
points, then $\dim\Delta{\,\le\,}18$\/}.
\par\smallskip
{\tt Proof.} (a) Suppose that $\dim\Delta{\,>\,}14$. Then each involution in $\Delta$ is planar, 
there is no $\Delta$-invariant proper subplane, and the center $\Zeta$ of $\Delta$ does not 
contain  an involution,  cf.~5.2. Moreover, $\Zeta$  acts freely on the affine point set, because 
$x^\zeta{\,=\,}x$ and $\zeta{\,\in\,}\Zeta\sm\{\1\}$ imply 
$\langle x^\Delta\rangle{\,\le\,}\cF_\zeta$. \\
(b) Let $\Gamma$ be an almost simple factor of $\Delta$, and put 
$\Psi{\,=\,}(\Cs\Delta\Gamma)^1$. In steps (c--g), assume that $\Gamma$ is a proper factor
of minimal dimension. \\
(c) {\it If $\Gamma$ contains an involution $\iota$, then $\dim\Delta{\,\le\,}18$\/}: consider 
$\Psi|_{\cF_\iota}{\,=\,}\Psi/\Kappa$ and $\Lambda{\,=\,}\Kappa^1$. According to 
\cite{sz1} 5.2, 
 we have  $\Psi{\hskip1pt:\hskip1pt}\Kappa{\,\le\,}9$, or $\Psi/\Kappa$ is locally isomorphic 
to  $\Opr5(\RR,2)$ and $\Psi$ has a subgroup $\SO3\RR$ by step (a). The factor $\Lambda$ 
of $\Delta$ is trivial or isomorphic to $\Spin3\RR$\, (see 2.6(\^b)\hskip1pt), but the latter is 
excluded  by  (a). If $\dim\Psi{\,=\,}9$, then $\dim\Gamma{\,=\,}3$ and $\dim\Delta{\,=\,}12$; 
if $\dim\Psi{\,<\,}9$, then $\dim\Delta{\,\le\,}16$. In the case 
$\dim\Gamma{\,=\,}\dim\Psi{\,=\,}10$ both factors would contain $\SO3\RR$, but this  contradicts 
 Lemma 2.11. Hence  $\dim\Gamma{\,\le\,}8$ and $\dim\Delta{\,\le\,}18$. \\
 (d) Suppose that $\{\1\}$ is the only compact subgroup of $\Delta$. Then each factor of 
 $\Delta$ is isomorphic to the simply connected covering group $\Omega$ of $\SL2\RR$. 
 If such a factor $\Gamma$ is straight, then Baer's theorem \cite{ba} implies
 $\Gamma$ consists of translations, because $\Gamma$ is not planar and does not consist of 
 homologies by \cite{cp} 61.2,  cf. also 5.2(g). All translations in 
 $\Delta$ have the same center (or the translation group would be commutative). Hence there 
 are at most two straight factors. \\
(e) For each $x{\,\notin\,}W$ and each factor $\Gamma$ of $\Delta$ the orbit 
$x^{\Gamma\smcap\Zeta}$ is non-trivial by step (a), and $x^\Gamma{\,\ne\,}x$. Consequently
$\cE{\,=\,}\cE(x,\Gamma){\,=\,}\langle x^{\Gamma\hskip-2pt},u,v,w\rangle$ is a connected subplane, and $\cE$ is not flat because $\Gamma$  acts non-trivially on $\cE$\, (see 
\cite{cp} 3.8).  If $\cE{\,\ledot\hskip-3pt}\cP$, then  
$\Psi_{\hskip-1pt x}|_\cE{\,=\,}\1$,\, $\Psi_{\hskip-1pt x}$ is 
compact by  2.6(b), and then $\Psi_{\hskip-1pt x}{\,=\,}\1$,\, $\dim\Psi{\,\le\,}15$, and 
$\dim\Delta{\,\le\,}18$.  We may assume, therefore, that each plane $\cE(x,\Gamma)$ is $4$-dimensional. \\
(f) Choose a factor $\Gamma$ which is not straight. If $L$ is a line in $\cE$ through 
$z{\,\in\,}\{u,v,w\}$, then $\dim\Gamma_{\hskip-2pt L}{\,>\,}0$ by the assumption in step~(e).
From $L^\Psi{\,=\,}\frL_z\sm W$ it follows that $\Gamma_{\hskip-2pt L}$ and hence $\Gamma$ consists of translations with center $z$,  and $\Gamma$ would be straight against our 
assumption. Consequently there are lines $K{\,=\,}xu$ and $L{\,=\,}xv$ such that 
$\dim K^\Psi\!,\:\dim L^\Psi{\,<\,}8$ and $\Psi{\hskip1pt:\hskip1pt}\Psi_{\hskip-1pt x}{\,\le\,}14$. \\
(g) Consider a factor $\hat\Gamma{\,\ne\,}\Gamma$, the product
$\Chi{\,=\,}\Psi\smcap\hat\Psi$ of the remaining factors,  and the corresponding planes 
$\cE,\,\hat\cE{\,=\,}\cE(x,\hat\Gamma)$, where $x$ is chosen as in step  (f) such that 
$\Chi{\hskip1pt:\hskip1pt}\Chi_x{\,\le\,}14$. By 2.14, we have 
$\cE^{\hskip1pt\hat\Gamma}{\,\ne\,}\cE$,\; $x^{\hat\Gamma}{\,\not\subseteq\,}\cE$,\; 
$\hat\cE{\,\ne\,}\cE$, and $\cF{\,=\,}\langle \cE,\hat\cE\rangle{\,\ledot\!}\cP$. 
Obviously, $\Chi_x|_\cF{\,=\,}\1$,\: $\Chi_x$ is compact by 2.6(b), and then 
$\Chi_x{\,=\,}\1$. Now $\dim\Chi{\,\le\,}12$ because $\Chi$ is a product of $3$-dimensional
factors, and $\dim\Delta{\,\le\,}18$. \\
(h) Corollary. {\it If $\dim\Delta{\,>\,}18$, and if $\Delta$ is not almost simple, then all factors of minimal dimension are isomorphic to the simply connected covering group $\Omega$ of $\SL2\RR$ and there is at least some other factor\/}. \\
(i) Assume from now on that $\dim\Delta{\,>\,}18$, and choose  a factor 
$\Gamma{\,\not\hskip1pt\cong\,}\Omega$ of $\Delta$. Then $\dim\Gamma{\,\ge\,}6$ and $\Gamma$ 
contains a planar involution $\iota$. The arguments in step (c) show that either 
$\dim\hskip-.5pt\Psi{\,\le\,}9$, or $\Psi{\,\circeq\,}\Opr5(\RR,2)$ has a subgroup $\SO3\RR$. 
In both cases $\dim\Gamma{\,\ge\,}10$. Note that this is true for each admissible choice of 
$\Gamma$. Consequently, $\Delta$ has no factors of dimension $6$ or~$8$. If 
$\dim\Gamma{\,=\,}10$, then $\Gamma$ contains $\SO3\RR$ because of step (a), and 
Lemma 2.11  implies $\dim\Psi{\,\le\,}9$. There is another involution 
$\iota'{\,\in\,}\Gamma$ which commutes with $\iota$, and $\iota'$ fixes a subplane 
$\cC{\,\ldot}\cF_\iota$, or it  induces a reflection on $\cF_\iota$ with some center 
$a{\,\notin\,}W$. In the first case, the stiffness result 
\cite{sz1} 1.5(4) shows that $\Psi$ acts 
almost effectively on $\cC$, and 
$\dim\Psi{\,\le\,}3$ by 2.14; in the second case, $a^\Psi{\,=\,}a$,\: $\Psi|_{a^\Gamma}{\,=\,}\1$,\: 
$\langle a^\Gamma\rangle{\,=\,}\langle a^\Delta\rangle{\,=\,}\cP$, and $\Psi{\,=\,}\1$. Hence we have even $\dim\Gamma{\,\ge\,}16$. By Priwitzer's results \cite{pw1} (see 5.2(j)\,), 
we may assume that $\dim\Gamma{\,\le\,}30$. This proves: \\
(j) {\it If $\dim\Delta{\,>\,}18$, then $\Delta$ is almost simple or a product of an almost simple 
group with $\Omega$\/}. \\
(k) The almost simple groups $\Gamma$ with $\dim\Gamma{\,\ge\,}20$ have been dealt with 
in 5.2 steps (j--l) under the assumption that $\Gamma$ fixes a flag (and possibly further 
elements).  In each case a contradiction has been obtained. Thus, only the possibility 
$\Gamma{\,\cong\,}{\rm (P)}\SL3\CC$ and $\Delta{\,=\,}\Gamma\Omega$ remains. 
An involution $\iota{\,\in\,}\Gamma$ is the center of a subgroup $\Ypsilon{\,\cong\,}\SL2\CC$ of 
$\Gamma$, and $\Ypsilon$ acts almost effectively on $\cF_\iota$. A maximal compact subgroup 
$\Phi{\,\cong\,}\SU2\CC$ of $\Ypsilon$ induces on $\cF_\iota$ a group $\SO3\RR$, and the 
involutions of $\Phi|_{\cF_\iota}$ are planar. Hence $\cF_\iota$ has lines homeomorphic to $\sS_4$, 
and Richardson's theorem$(\dagger)$ implies that $\cE{\,=\,}\cF_{\iota,\Phi}$ is flat. 
By Stiffness, $\Omega$ acts almost effectively on $\cE$, and $\Omega|_\cE$   is almost simple; 
on the other hand, $\Omega|_\cE$ is solvable by \cite{cp} 33.8. This contradiction completes 
the proof.  \qed
 \par\medskip
 {\bf 7.3 Normal torus.} {\it If $\Delta$ has a one-dimensional compact connected normal subgroup 
$\Theta$, and if $\dim\Delta{\,\ge\,}15$, then there exists a $\Delta$-invariant classical quaternion subplane $\cH$, and $\dim\Delta{\,\le\,}16$\/}.
\par\smallskip
 {\tt Proof.} $\Delta$ is a Lie group by 7.0, hence $\Theta{\,\cong\,}\TT$ is a torus.
 As in 6.3, the involution $\iota{\,\in\,}\Theta$ is in the center of $\Delta$, and 
 $\cH{\,:=\,}\cF_\iota{\,=\,}\cH^\Delta$ is a Baer subplane. Put again 
 $\Delta|_\cH{\,=\,}\Delta/\Kappa$.  If $\Theta|_\cH{\,\ne\,}\1$, then it 
 follows from 
\cite{sz1} 7.3 and Stiffness that $\dim\Delta{\,\le\,}7{+}3$. 
 Therefore $\Theta{\,\le\,}\Kappa$,\: $\Kappa^1{\,\not\cong\,}\Spin3\RR$,\: $\dim\Kappa{\,=\,}1$, 
 and $14{\,\le\,}\Delta{\hskip1pt:\hskip1pt}\Kappa{\,\le\,}15$. Now the claim is a consequence of \cite{sz1} 7.3 or  \cite{sz2} Th.~5. \qed
 \par\medskip
{\bf 7.4 Theorem.} {\it If $\dim \Delta \ge 32$ and $\Delta$ has\/} ({\it at least\/}) {\it $3$ fixed 
points, then $\Delta$ contains a trans\-itive translation group $\Tau$. Either $\dim \Delta = 32$ 
and a maximal semi-simple subgroup $\Psi$ of $\Delta$ is isomorphic to  $\SU 4 \CC$, 
or $\dim \Delta \ge 37$ and $\cP \cong \cO$\/}. 
\par\smallskip
This has been {\tt proved} in \cite{sz14}; the result seems to be best possible. 
Translation planes with a group locally isomorphic to $\SU4\CC$ have been studied in a long paper 
by  H\"ahl \cite{Ha4} without any assumption on the existence of a fixed point on the translation axis:  
\par\smallskip
{\bf 7.5 Theorem.} {\it Let $\Sigma$ be the automorphism group of a translation plane with 
translation axis $W$. If $\Sigma$ has a subgroup $\Ypsilon{\,\circeq\,}\SU4\CC$ with a unique 
fixed point $o{\,\notin\,}W$, then $\Ypsilon$ is normal in the connected component of $\Sigma_o$ 
or $\cP$ is the classical octonion plane\/}. 
\par\smallskip
{\tt Remark.} Specific non-classical examples have been given in 4.12.
\par\bigskip
 
{\Bf 8. Fixed double flag}
\par\medskip
Throughout this section, $\Delta$ fixes exactly $2$ points and $2$ lines, i.e., 
$\cF_\Delta{\,=\,}\langle u,v,av\rangle$ is a {\it double flag\/}. In the classical plane $\cO$ a maximal 
semi-simple subgroup $\Phi$ of $\Sigma_{u,v,av}$ is isomorphic to $\Spin8\RR$ and fixes even a 
triangle.
\par\medskip
{\bf 8.0  Semi-simple Lie groups.} {\it If $\Delta$ is semi-simple, $\cF_\Delta{\,=\,}\langle u,v,av\rangle$ is a double flag, and $\dim\Delta{\,\ge\,}26$, then $\Delta$ is a Lie group\/}.
\par\smallskip  
{\tt Proof.} By \cite{psz} or 2.3 it suffices to consider the case $\dim\Delta{\,=\,}26$. Write 
$\Delta{\,=\,}\Gamma\Psi$, where $\Gamma$ is an almost simple factor of {\it maximal\/} dimension. Note that $\dim\Gamma{\,\ge\,}8$ since $\dim\Delta{\,\equiv\,}2\bmod3$. \par
(a) $\Delta$ does not contain a reflection $\sigma$ with center $v$ or axis $uv$\:  (or else 
$\sigma^\Delta\sigma$ would be a commutative normal subgroup, see 2.9). Therefore  each involution in $\Delta$ is either planar or a reflection in $\Delta_{[u,av]}$, and one of two commuting involutions is planar, see 2.10. Moreover, any non-central involution $\iota$ of an almost simple factor $\Chi$   is planar: if $\iota$ is a reflection, then $\langle \iota^\Chi\rangle{\,=\,}\Chi{\,\le\,}\Delta_{[u,av]}$, and 
$\Chi$ is compact or two-ended by \cite{cp} 61.2, hence a compact Lie group of rank $1$. 
If $\Chi{\,\cong\,}\SO3\RR$, the involutions in $\Chi$ are planar; if $\Chi$ is simply connected,
the only involution in $\Chi$ is central. \par
(b) {\it If $\Ypsilon$ is a semi-simple group in the centralizer of a planar involution $\beta$, 
then~$\dim\Ypsilon{\le13}$\/}:\: Consider the induced group 
$\Ypsilon|_{\cF_\beta}{\,=\,}\Ypsilon/\Kappa$. From 
\cite{sz1} 6.1 it follows that 
$\Ypsilon{\hskip1pt:\hskip1pt}\Kappa{\,\le\,}10$. On the other hand, $\Kappa^1{\,=\,}\Lambda{\,\triangleleft\,}\Ypsilon$ 
is compact and semi-simple, and then $\Lambda$ is a Lie group (cf. \cite{cp} 93.11), and Stiffness implies  $\Lambda{\,\cong\,}\Spin3\RR$ or $\Lambda{\,=\,}\1$. \\
(c) Suppose that $\dim\Gamma{\,=\,}20$. Then $\Gamma{\,\cong\,}\Sp4\CC$\; (because 
$\PSp4\CC$ contains $\SO5\CC{\,>\,}\SO5\RR$) and 
$\beta{\,=\,}{\rm diag}(1,1,-1,-1)$ is planar. We have $\Cs{}\beta{\,\ge\,}(\SL2\CC)^2{\cdot}\Psi$, and 
$\dim\Cs{}\beta$ would be too large. \\
(d) If $\dim\Gamma{\,=\,}16$, then $\Gamma{\,\cong\,}{\rm(P\hskip-1pt)}\SL3\CC$ and 
${\rm diag}(1,-1,-1)$ corresponds to a planar involution $\beta$ containing 
$\SL2\CC{\hskip1pt\cdot\Psi}$  in its centralizer, and $\dim\Cs{}\beta$ is again too large. \\
(e) Each group $\Gamma$ of type ${\rm A}_3$ contains a planar involution $\beta$ centralizing 
a $6$-dimensional semi-simple subgroup of $\Gamma$. This is easy in the cases $\SU4\CC$,\,
$\SL2\HH$,\, $\SL4\RR$ and its universal covering group; it is less obvious for the groups 
related to $\SU4(\CC,r)$, $r{\,=\,}1,2$, which are not necessarily Lie groups, 
see 8.2(j) Case 4) below. Once again this contradicts step (b). \\
(f) If $\Gamma$ is the compact group $\Gtwo$, then $\Psi$ would be solvable by Ê\cite{cp} 33.8. 
In the case $\Gamma{\,\circeq\,}\Gtwo(2)$ a maximal compact subgroup $\Phi$ of $\Gamma$ is 
isomorphic to $\SO4\RR$ or to $(\Spin3\RR)^2$, and $\Phi\Psi$ would centralize a planar involution  in contradiction to step (b). \\ 
(g) {\tt Corollary.} {\it All factors of $\Delta$ have dimension at most $10$\/}. \\
(h) {\it If $\dim\Gamma{\,=\,}10$, then $
\Delta$ is a Lie group\/}. Suppose that 
$x^\Gamma{\,=\,}x$ for some point $x{\,\notin\,}av\smcup uv$. Then $x^\Delta{\,=\,}x^\Psi$,\: 
$\cD{\,:=\,}\langle x^\Psi\rangle{\,=\,}\cD^\Delta$, and $\Gamma|_\cD{\,=\,}\1$. 
Now $\cD$ is flat by Stiffness, but then $\dim\Delta$ would be too small. 
Dually each fixed line of $\Gamma$ contains either $u$ or $v$. 
As $\dim\Delta_{[u,av]}{\,\le\,}8$, we have $\Gamma{\,\ne\,}\Gamma_{\hskip-2pt[u]}$.
If $\Gamma{\,=\,}\Gamma_{\hskip-2pt[v]}$, then the dual of \cite{cp} 61.20(a) shows that 
$\Gamma{\,\cong\,}\Gamma_{\hskip-2pt[au]}{\ltimes}\Gamma_{\hskip-2pt[uv]}$ is not 
almost simple. Hence there is a point $p{\,\notin\,}av\smcup uv$ such that  
$\cE{\,=\,}\langle p^\Gamma\rangle$ is a subplane, and  $\cE{\ledot\hskip-3pt}\cP$, or else 
$\dim\Gamma|_\cE{\,\le\,}6$.
As $\Psi_{\hskip-1pt p}|_\cE{\,=\,}\1$, Stiffness shows that $\Psi_{\hskip-1pt p}$ is compact and 
$\dim\Psi_{\hskip-1pt p}{\,\le\,}7$.
Therefore $\dim p^\Psi{>\,}8$, \: $\langle p^\Psi\rangle{\,=\,}\cP$, and 
$\Gamma_{\hskip-2pt p}{\,=\,}\1$. It follows that $\cE{\,=\,}\cP$,\: $\Psi_{\hskip-1pt p}{\,=\,}\1$, and
$\dim p^\Psi{\,=\,}16$. According to \cite{cp} 96.11(a) the orbit $p^\Psi$ is open in $P$, and then
$\Psi$ is a Lie group by \cite{cp} 53.2. The center of $\Gamma$ contains arbitrarily small compact 
subgroups $\Nu$ such that $\Gamma/\Nu$ is a Lie group (cf. \cite{cp} 93.8.) 
For each $\nu{\,\in\,}\Nu$
there exists a unique $\hat\nu{\,\in\,}\Psi$ with $p^\nu{\,=\,}p^{\hat\nu}$. Because 
$\Nu{\,\le\,}\Cs{}\Delta$, the map $(\nu{\,\mapsto\,}\hat\nu){\,:\,}\Nu{\,\to\,}\Psi$ is a continuous homomorphism, $\Nu$ is a Lie group, and so are $\Gamma$~and~$\Delta$. \\
(i) {\it If $\dim\Gamma{\,=\,}8$, then $\Delta$ is a Lie group\/}. 
The first arguments are similar to those in step~(h), but not quite the same. 
If $x^\Gamma{=\,}x{\,\notin\,}av\smcup uv$ and 
$\cD{\,=\,}\langle x^\Psi\rangle$, then $\Gamma|_\cD{\,=\,}\1$ and
$\dim\cD{\,\le\,}4$. Let $y{\,\in\,}\cD$ such that $\Lambda{\,=\,}\Delta_{x,y}$ fixes a quadrangle. 
It follows that $\Delta{\,:\,}\Lambda{\,\le\,}8$ and $\dim\Lambda{\,\le\,}14$. Thus $\dim\Delta$ 
would be too small. Dually $L^\Gamma{\,=\,}L{\,\in\,}\frL$ implies $L{\,\in\,}\frL_u\smcup\frL_v$. 
If $\Gamma{\,=\,}\Gamma_{\hskip-1.5pt[u]}$, then $\Gamma$ is compact by \cite{cp} 61.2 
and transitive on $uv\sm\{u,v\}$, which is impossible. If 
$\Gamma{\,=\,}\Gamma_{\hskip-2pt[v]}$, then $\Gamma{\,\ne\,}\Gamma_{\hskip-2pt[v,uv]}$ \  
(note that $\Gamma$ contains an involution and apply \cite{cp} 55.28). Therefore
$\Gamma{\,\cong\,}\Gamma_{\hskip-2pt[au]}{\ltimes}\Gamma_{\hskip-2pt[uv]}$ by
\cite{cp} 61.20, and $\Gamma$ would not be almost simple. 
Again there is a point $p$ such that  $\cE{\,=\,}\langle p^\Gamma\rangle{\,\ledot\!}\cP$.
We have  $\Psi_{\hskip-1pt p}|_\cE{\,=\,}\1$ and $\dim\Psi_{\hskip-1pt p}{\,\ge\,}2$ because $\dim\Psi{\,=\,}18$. 
Hence $\cE$ is a Baer subplane, $\dim\Psi_{\hskip-1pt p}{\,\le\,}7$ by Stiffness, $\dim p^\Psi{\,\ge\,}11$, 
and $\langle p^\Psi\rangle{\,=\,}\cP$. The center of $\Delta$ contains a compact subgroup 
$\Xi{\,<\,}\Psi$ such that $\Psi/\Xi$ is a Lie group. Obviously
$\Psi_{\hskip-1pt p}|_{p^{\Gamma\Xi}}{\,=\,}\1$,\;   $\cE{\,<\,}\cP$, and 
$\langle p^{\Gamma\Xi}\rangle{\,=\,}\cE$. As $\langle p^\Psi\rangle{\,=\,}\cP$, it follows that
$(\Gamma\Xi)_p{\,=\,}\1$ and $\dim p^{\Gamma\Xi}{\,=\,}8$. Now  $p^{\Gamma\Xi}$ is 
an open point set in $\cE$, and $\Gamma\Xi$ is a Lie group by \cite{cp} 53.2. 
Consequently $\Gamma$ and $\Psi$ are Lie groups. \qed
\par\medskip
{\bf 8.1 Lemma.} {\it If $\Delta$ is almost simple, and if $\cF_\Delta$ is a double flag, then each non-central involution in $\Delta$ is planar\/}, see step (a) in the previous {\tt proof}.
\par\medskip
{\bf 8.2 Semi-simple groups.} {\it If $\cF_\Delta$ is a double flag and $\Delta$ is a 
semi-simple Lie group, then $\dim\Delta{\,\le\,}24$\/}.
\par\smallskip
Because of the possible existence of a reflection in $\Delta_{[u,av]}$, the {\tt proof} differs in 
several steps from the proof of 7.2; accordingly, the result is much weaker. By Priwitzer's 
Theorem 4.2, we may assume $\dim\Delta{\,\le\,}28$.  \\
(a) {\it If there exists a $\Delta$-invariant proper subplane $\cF$, then 
$\dim\Delta{\,\le\,}18$\/}:  write $\Delta|_\cF{\,=\,}\Delta/\Kappa$, and put $\Lambda{\,=\,}\Kappa^1$ 
and $\dim\cF{\,=\,}d$. Note that $\cF_\Delta$ need not be contained in $\cF$.  If $d{\,\le\,}2$, then $\Delta{\,=\,}\Kappa$ or $\Delta/\Kappa$ is solvable 
by \cite{cp} 33.8, hence $\Delta{\,=\,}\Lambda$ and $\dim\Lambda{\,<\,}14$ by Stiffness. 
For $d{\,=\,}4$, there are $3$ subcases: if $\cF{\,=\,}\cF_\Kappa$, then $u,v,av{\,\in\,}\cF$,\; 
$\Delta{\hskip1pt:\hskip1pt}\Kappa{\,\le\,}3$ by 2.14, and $\dim\Kappa{\,\le\,}8$.
If $\cF{\,\ldot}\cF_\Kappa{\,\ldot}\cP$, then $v{\,\in\,}\cF$ or $v$ is incident with a unique line 
$L{\,=\,}L^\Delta{\,\in\,}\cF$, $\Delta{\hskip1pt:\hskip1pt}\Kappa{\,\le\,}8$ by \cite{cp} 71.8 and 
$\dim\Kappa{\,\le\,}7$. If $\cF_\Kappa{\,=\,}\cP$, then $\Kappa{\,=\,}\1$ and 
$\Delta$ is a subgroup of the $16$-dimensional group $\PSL3\CC$.
Finally, let $d{\,=\,}8$. Either $u,v,av{\,\in\,}\cF$,\; $\Delta{\hskip1pt:\hskip1pt}\Kappa{\,\le\,}10$ by \cite{sz1} 6.1, and $\Lambda$ is compact and semi-simple, hence a Lie group contained in 
$\Spin3\RR$, or $\cF_\Kappa{\,=\,}\cP$,\, $\Kappa{\,=\,}\1$, and $\Delta$ acts faithfully on $\cF$. 
Each exterior  fixed element is incident with a unique interior element, and this is 
$\Delta$-invariant. From \cite{sz1} 7.3 it follows that  $\dim\Delta{\,\le\,}13$, or 
$\cF\smcap\cF_\Delta$ is a non-incident point-line pair and $\cF$ is the classical quaternion 
plane. In the latter case $\dim\Delta{\,\le\,}18$. (This is best possible: 
the stabilizer $\Chi$ of an exterior flag in a proper Hughes plane fixes a non-incident interior point-line pair and hence a double flag in $\cP$; the group $\Chi$ has a  semi-simple subgroup $\Delta$ of  dimension $18$.) This proof  shows also the following: \\
(b) Corollary. {\it If $u,v,av{\,\in\,}\cF$, in particular, if
$\cF_\Xi{\,=\,}\cF_{\hskip-.7pt\Xi}^{\hskip2pt\Delta}{\,<\,}\cP$ for some set 
$\Xi{\,\subseteq\,}\Delta$, then $\dim\Delta{\,\le\,}13$. Hence each central involution in $\Delta$ 
is a reflection with center $u$ and axis $av$ whenever $\dim\Delta{\,>\,}13$\/} (see 2.\,9 and 10). \\
(c) Let $\Gamma$ be a proper almost simple factor of $\Delta$ and put 
$\Psi{\,=\,}(\Cs\Delta\Gamma)^1$. Assume first that $\Gamma$ has $minimal$ dimension. If 
$\Gamma$ contains a planar involution $\iota$, we write $\Psi|_{\cF_\iota}{\,=\,}\Psi/\Kappa$.   
According to \cite{sz1} 6.1, we have $\Psi{\hskip1pt:\hskip1pt}\Kappa{\,\le\,}10$. Either $\Kappa$ 
is discrete or $\Kappa^1{\,\cong\,}\Spin3\RR$. In the latter case, $\dim\Gamma{\,=\,}3$ by 
minimality, $\dim\Psi{\,\le\,}13$, and $\dim\Delta{\,\le\,}16$. \\
(d) Suppose from now on that $\dim\Delta{\,\ge\,}20$. If $\Delta$ has two distinct factors 
$\Gamma_{\!\nu}{\,\cong\,}\Spin3\RR$, then step (b) implies that the central involution 
$\iota_\nu{\,\in\,}\Gamma_{\!\nu}$ is a reflection in $\Delta_{[u,av]}$. As 
$\iota_1\iota_2{\,=\,}\iota_2\iota_1$, it follows from 2.10 that $\iota_1{\,=\,}\iota_2{\,=\,}\iota$, and 
$\Gamma_{\!1}\Gamma_{\!2}{\,\cong\,}\SO4\RR$ has a subgroup $\Phi{\,\cong\,}\SO3\RR$. 
Any involution $\beta{\,\in\,}\Phi$ is distinct from $\iota$, hence planar, and the arguments of step (c) show that $\dim(\Psi_1\smcap\Psi_2){\,\le\,}13$,\: $\dim\Delta{\,<\,}20$ contrary to the assumption.  More generally, the same reasoning shows that at most one $3$-dimensional factor has a non-trivial torus subgroup. (Note that a subgroup 
$\TT^2{\,<\,}\Gamma_{\!1}\Gamma_{\!2}$ contains a planar involution.)\\
(e) {\it If $\Gamma$\hskip1pt is straight, then $\Gamma{\,\cong\,}\Spin3\RR$ is a compact group of homologies in $\Delta_{[u,av]}$, or $\Gamma$ is con\-tained in the translation group 
$\Tau{\,=\,}\Delta_{[v,uv]}$ and $\Gamma$ has no compact subgroup other~than~$\{\1\}$. 
In particular, $\dim\Gamma{\,=\,}3$. At most one straight factor of $\Delta$ is compact, at most 
two consist of translations\/}. By Baer's theorem, either $\cF_\Gamma{\;\ldot}\cP$ is  a $\Delta$-invariant subplane, which contradicts steps (a,b), or center and axis of $\Gamma$ are fixed 
elements of $\Delta$. If $\Gamma{\,\le\,}\Delta_{[u]}$, then $\Gamma$ is compact by the second part of  2.12, and $6{\,\ne\,}\dim\Gamma{\,<\,}8$; if  $\Gamma{\,\le\,}\Delta_{[v]}$ and hence  
$\Gamma{\,\le\,}\Tau$, use \cite{cp} 55.28. \\
(f) {\it If $\Delta$ is a product of $3$-dimensional factors, then $\dim\Delta{\,\le\,}21$\/}: By the assumption in step~(d), there are at least $7$ such factors; at most $3$ of them can be straight, 
see step~(e). Hence there is a point $x$ and a factor $\Gamma$ such that 
$\cE{\,=\,}\langle x^\Gamma,u,v,av\rangle$ is a subplane, $\cE$~is not flat by \cite{cp}~3.8, and 
$\Psi_{\hskip-1pt x}|_\cE{\,=\,}\1$. If  $\cE{\,\ledot\!}\cP$, then  Stiffness implies 
$\dim\Psi{\,\le\,}16{+}3$.  As 
$\dim\Psi{\,\equiv\,}0{\,\bmod\,}3$, it follows that $\dim\Psi{\,\le\,}18$.   We may assume, therefore, 
that $\dim\cE{\,=\,}4$. If possible, choose $\Gamma$ as the factor of torus rank $1$. Otherwise let 
$\hat\Gamma$ denote the straight compact factor if such a factor exists. In the second case, 
$\cG{\,=\,}\langle x^{\Gamma\hskip1pt\hat\Gamma},u,v,av\rangle{\ledot\hskip-3pt}\cP$ because 
$\hat\Gamma$ consists of homologies; in the first case, $\cE^{\hat\Gamma}{\,\ne\,}\cE$ for some factor $\hat\Gamma$ by step (b), so that again $\cG{\,\ledot\hskip-3pt}\cP$.
Let $\Chi{\,=\,}\Psi\smcap\hat\Psi$ 
be the product of the remaining factors. Then $\Chi_x|_\cG{\,=\,}\1$,\: $\Chi_x$ is compact, and 
then $\Chi_x{\,=\,}\1$ by the choice of $\Gamma$ and $\hat\Gamma$. Consequently, 
$\dim\Chi{\,<\,}16$, and (f) is proved. \\
(g) {\it If some almost simple factor $\Gamma\!$ of $\Delta$ has dimension 
$\dim\Gamma{\,\in\,}\{6,8\}$, then  $\dim\Delta{\,\le\,}26$\/}; {\it moreover, $\dim\Gamma{\,=\,}6$ or 
$\dim\Psi{\,<\,}19$\/}. 
By step (e), $\Gamma$ is not straight. Hence there is a point $x$ such that $\cE{\,=\,}\langle x^\Gamma,u,v,av\rangle$ is a subplane. 
2.14 implies  $\cE{\ledot\hskip-3pt}\cP$ (note that $\Gamma$ acts non-trivially on~$\cE$). Again 
$\Psi_{\hskip-1pt x}|_\cE{\,=\,}\1$,\: $\dim\Psi_{\hskip-1pt x}{\,\le\,}3$, and $\dim\Psi{\,\le\,}16{+}3$. In the case of equality, 
$\cE$ is a Baer subplane, $\dim\Psi_{\hskip-1pt x}{\,=\,}3$,\: $\dim x^\Psi{\,=\,}16$ and $x^\Psi$ is open in $P$, see \cite{cp} 96.11. 
By construction, $x^\psi{\,\in\,}\cE{\,\Leftrightarrow\,}\cE^\psi{\,=\,}\cE$ for $\psi{\,\in\,}\Psi$. Put 
$\Xi{\,=\,}\Psi_\cE$. Then $x^\Xi$ is an open point set in $\cE$, and $\Xi/\Psi_{\hskip-1pt x}{\,\approx\,}x^\Xi$\: 
(cf. \cite{cp} 96.9). Consequently, $\dim\Xi{\,=\,}8{+}3$. If $\dim\Gamma{\,>\,}6$, then
 $\dim\Gamma\Xi|_\cE{\,\ge\,}8{+}11{-}3$, but the stabilizer of a double flag in $\cE$ has 
dimension  ${\le\,}15$ by \cite{cp} 83.26. Hence $\dim\Gamma{\,=\,}6$ or $\dim\Psi{\,\le\,}18$. 
More can be shown, using the structure of $\Delta$: \\
(h) {\it If $\Delta$ has an almost simple factor $\Gamma$ of dimension $6$, then 
$\dim\Delta{\,\le\,}24$\/}. \\ 
Suppose that $\dim\Delta{\,>\,}24$. Then $\dim\Psi{\,=\,}19$. This case has been discussed in the previous step, and the same notation will be used. Step (b) shows that $\Gamma$ 
contains a unique involution, and this is a reflection. Consequently, $\Gamma{\,\cong\,}\SL2\CC$,\: 
$\Gamma$ is not a group of homologies by \cite{cp} 61.2, and $\Gamma|_{av}{\,\cong\,}\SO3\CC$. 
Moreover, $\cE{\,\ldot}\cP$,\; $x^\Xi$ is open in the point space of $\cE$, and $\Xi$ has an 
open orbit on $av\smcap\cE$. Therefore the lines of $\cE$ are manifolds  homeomorphic to 
$\sS_4$, see  \cite{cp} 53.2. According to Richardson's theorem ($\dagger$), 
a maximal compact subgroup $\Phi$ of $\Gamma$ fixes a circle $C{\,\subset\,}av{\smcap}\cE$. Consider the $11$-dimensional group 
$\Chi{\,=\,}\Phi\Xi|_\cE$, and  note that $C^\Chi{\,=\,}C$. Let $a,b{\,\in\,}C\sm\{v\}$. Then $\dim\Chi_{a,b}{\,\ge\,}9$. On the other hand, $\dim\Chi_{a,b}{\,\le\,}7$ by the stiffness result \cite{sz15} (*) or \cite{sz1} 1.5(6). \\ 
(i) {\it If some almost simple factor $\Gamma\!$ of $\Delta$ has dimension 
$\dim\Gamma{\,=\,}8$, then $\dim\Delta{\,\le\,}24$\/}. \\
Let $\dim\Delta{\,>\,}24$. Step (g) implies  $\dim\Psi{\,\in\,}\{17,18\}$. 
If $\Psi$ is an almost direct product $\Alpha\Beta$ with $\dim\Beta{\,>\,}10$ and 
$\rk\Beta{\,\ge\,}2$, then $\Beta$ contains a planar  involution $\beta$, and 
$\dim\Alpha\Gamma|_{\cF_\beta}{\,\le\,}10$ by \cite{sz1} 6.1. The kernel of the action of 
$\Alpha\Gamma$ on $\cF_\beta$ is a compact normal  subgroup of $\Alpha\Gamma$ of positive dimension ${\le\,}3$, hence it coincides with $\Alpha$,\: $\cF_\Alpha{\,=\,}\cF_\beta$, and 
$\Alpha{\,\cong\,}\Spin3\RR$ by Stiffness. Now $\cF_{\hskip-2pt\Alpha}^{\hskip1pt\Beta}{\,=\,}
\cF_\Alpha$ is $\Delta$-invariant in contradiction to  step (b). \\ 
Next, let $\Delta$ be~a product of two $8$-dimensional almost 
simple factors $\Gamma\hskip-1pt,\, \hat\Gamma$ and one $10$-dimensional factor~$\Omega$. 
Each of these factors has a $3$-dimensional compact subgroup. If some factor contains a planar involution $\iota$\: (in particular, if the factor has rank $2$, or if it has a subgroup $\SO3\RR$), 
then the product of the other two factors of $\Delta$ induces on $\cF_\iota$ a semi-simple group of dimension at least $16$, but this contradicts 
\cite{sz1} 6.1. Hence each of the $3$ factors of 
$\Delta$ has a maximal compact subgroup isomorphic to $\Spin3\RR$. The central involutions of 
these subgroups coincide with a unique reflection $\sigma{\,\in\,}\Delta_{[u,av]}$. There is some point $x$ such that   $\cE{\,=\,}\langle x^{\hat\Gamma},u,v,av\rangle{\:\ledot\!}\cP$. 
Let $\Lambda$ be the connected component of   $(\Gamma\Omega)_x$, and note  that 
$\dim\Lambda{\,\ge\,}2$ and that 
$\Lambda|_\cE{\,=\,}\1$. By Stiffness, $\Lambda{\,\cong\,}\Spin3\RR$, and $\Lambda$ is contained 
in a maximal compact subgroup $\Phi{\,\cong\,}\SO4\RR$ of $\Gamma\Omega$. 
It follows that $\Lambda$ is one of the two normal factors of $\Phi$. Consequently, $\Lambda$ is contained in $\Gamma$ or in $\Omega$, but $\Lambda$ is planar, a contradiction. \\
The case of two $8$-dimensional factors $\Gamma\hskip-1pt,\:\hat\Gamma$ and factors 
$\Alpha_\nu,\: \nu{=}0,1,2$ of dimension $3$ can be dealt with similarly: by step (d), at most one of the $\Alpha_\nu$, say $\Alpha_2$, can have positive rank. Because of 
\cite{sz1} 6.1, the fixed point set of 
$\Alpha_\nu$ is at  most $4$-dimensional. By Baer's theorem, there is some point $x$ such that 
 $\cE{\,=\,}\langle x^{\Alpha_1\Alpha_2},u,v,av\rangle$ is a subplane and that $\Alpha_1\Alpha_2$ 
 induces a $6$-dimensional group on $\cE$. From 2.14 it follows that $\dim\cE{\,>\,}4$, and 
 $\cE{\ledot\hskip-3pt}\cP$. Put $\Xi{\,=\,}\Gamma\hat\Gamma\Alpha_0$. Then $\dim\Xi{\,=\,}19$,\: 
 $\dim\Xi_x{\,=\,}3$,\: $\Xi_x|_\cE{\,=\,}\1$, and $\Lambda{\,=\,}(\Xi_x)^1{\,\cong\,}\Spin3\RR$ 
is contained in a maximal compact subgroup $\Phi$ of $\Xi$. Again $\Phi{\,\cong\,}\SO4\RR$. For the same reason as before, this is impossible.  
The only remaining case that $\Psi$ is a product of six  $3$-dimensional factors can be excluded by the arguments  of step  (f).   \\
(j) {\it If $\Delta$ has a $10$-dimensional factor $\Gamma$, then  $\dim\Delta{\,\le\,}24$\/}. \\
We may assume that $25{\,\le\,}\dim\Delta{\,\le\,}28$ and that there is no factor of dimension $6$ 
or $8$. If $\Gamma$ contains two commuting involutions, then there is a planar involution 
$\beta{\,\in\,}\Gamma$, $\dim\Psi|_{\cF_\beta}{\,\le\,}10$ by 
\cite{sz1} 6.1, and $\dim\Delta{\,\le\,}23$. 
Therefore $\Gamma$ is isomorphic to the simply connected covering group  of  $\Opr5(\RR,2)$, and the unique central involution $\sigma{\,\in \,}\Gamma$ is a reflection  in $\Delta_{[u,av]}$. \\ 
Case 1) There are exactly two almost simple factors $\Alpha,\Ypsilon{\,\ne\,}\Gamma$, where 
$\dim\Ypsilon{\,\in\,}\{14,15\}$ and $\rk\Ypsilon{\,=\,}2$. Then $\Ypsilon$ contains a planar involution 
$\iota$, and $\dim\Gamma\Alpha|_{\cF_\iota}{\,\le\,}10$ by 
\cite{sz1} 6.1. Hence 
$\Alpha|_{\cF_\iota}{\,=\,}\1$ and $\Alpha{\,\cong\,}\Spin3\RR$ by Stiffness. The involution 
$\alpha{\,\in\,}\Alpha$ is planar, and $\Gamma\Ypsilon|_{\cF_\alpha}$ would be too large. \\
Case 2) $\Delta$ has at least two $3$-dimensional factors, $\Alpha$ and $\Beta$. At most one of 
 them consists of homologies with axis $av$\: (see steps (d,e)\,). Let $\Alpha|_{av}{\,\ne\,}\1$. 
 Then there is a point $x$ such that $\langle x^\Alpha,u,v,av\rangle$ is a subplane, and so is  
 $\cE{\,=\,}\langle x^{\Alpha\Beta},u,v,av\rangle$. If $\dim\cE{\,<\,}8$, then 2.14 shows that
$\Alpha\Beta|_\cE$ is almost simple. By assumption $\Alpha|_\cE{\,\ne\,}\1$. Hence
 $\Beta|_\cE{\,=\,}\1$,\: $\cE{\,\le\,}\cF_\Beta$,\:  
 $\cF_{\hskip-1pt\Beta}^{\hskip1pt\Delta}{\,=\,}\cF_\Beta{\,<\,}\cP$,  and  $\dim\Delta{\,\le\,}13$ by  
Corollary (b). Therefore $\cE{\ledot\hskip-4pt}\cP$. Put $\Chi{\,=\,}({\Alpha\Beta})^1$.
  Then $\Chi_x|_\cE{\,=\,}\1$,\: $\dim\Chi{\,\le\,}16{+}3$, and 
$\dim\Delta{\,\le\,}25$. This takes care of all possibilities except the following ones: $\Psi$ is an almost simple group of dimension $16$, or $\dim\Psi{\,=\,}15$ and $\Psi$ is almost simple or a product of
$3$-dimensional factors. \\
Case 3) $\Psi$ is isomorphic to (P\hskip-1pt)$\hskip-1pt\SL3\CC$. There are $3$ pairwise commuting involutions conjugate to $\beta{\,=\,}{\rm diag}(-1,-1,1)$; they are necessarily planar, and 
$(\Cs\Psi\beta)'{\,\cong\,}\SL2\CC$. Hence $(\Cs\Delta\beta)'$ induces on $\cF_\beta$ a $16$-dimensional semi-simple group in contradiction to 
\cite{sz1}  6.1. \\
Case 4) Similarly, the groups $\Psi$ of type ${\rm A}_3$ can be dealt with by exhibiting a suitable 
semi-simple group $\Ypsilon$ in the centralizer of a planar involution $\beta$. \\
\begin{center} 
\begin{tabular}{|c|c|c|c|}
\hline
$\Psi$ & $\beta$ & $\Ypsilon$& $\Psi/\langle -1\rangle$ \\
\hline
$\SU4\CC$ & $\:{\rm diag}(-1,-1,1,1)\:$ & $(\SU2\CC)^2$ & $\SO6\RR$ \\
$\SL2\HH$ & ${\rm diag}(-1,1)$ & $(\SU2\CC)^2$ & $\Opr6(\RR,1){>}\SO5\RR$ \\
\hline
$\SU4(\CC,1)$ & ${\rm diag}(-1,-1,1,1)$ & $\:\SU2\CC{\times}\SL2\RR\:$& $\SaU3$ \\
$\SU4(\CC,2)$ & ${\rm diag}(-1,-1,1,1)$ & $(\SU2\CC)^2$ & $\Opr6(\RR,2)$  \\ 
\hline
$\SL4\RR$ & ${\rm diag}(-1,-1,1,1)$ & $(\SL2\RR)^2$ & $\Opr6(\RR,3){>}(\SO3\RR)^2$  \\
$\tilde{\rm SL}_4\RR$ & $\beta{\,\in\,}\Spin4\RR$ & $\Spin4\RR$ & \\
\hline
\end{tabular}\end{center}
Here $\tilde\Psi$ denotes the simply connected covering group of $\Psi$.
Recall from 2.\,10, and 11  that $\Delta$ has no subgroup $\SO5\RR$ or $(\SO3\RR)^2$. Hence $\Psi$ 
cannot be a proper factor group of $\SU4\CC,\: \SL4\RR$, or $\SL2\HH$. (Note that 
$\{\hskip1pt({a \atop  } { \atop \raise2pt\hbox{${\scriptstyle b}$}})\mid a,b{\,\in\,}\HH'\,\} {\,\le\,}\SL2\HH$.)
For  $\Psi{\,\cong\,}\SU4(\CC,1)$, the involution $\beta$ is contained in $\SU3\CC$ and corresponds 
to a planar involution in each factor group of $\Psi$ as well as in each covering group; the centralizer 
of such an involution contains a semi-simple group locally isomorphic to  
$\SU2\CC{\times}\SU2(\CC,1){\,\cong\,}\SU2\CC{\times}\SL2\RR$.
Similarly, if $\Psi{\,\cong\,}\SU4(\CC,2)$, then $\beta{\,\in\,}(\SU2\CC)^2$ commutes with the 
reflection $-\1$; hence $\beta$ corresponds to a planar involution in each factor group and in each covering group of $\Psi$, centralizing a $6$-dimensional compact semi-simple group. \\
Case 5) All factors of $\Delta$ other than $\Gamma$ are $3$-dimensional. Recall that $\Gamma$ 
is the simply connected covering group of $\Opr5(\RR,2)$. This case
 is more complicated. Choose $\Alpha$ and $\Beta$ among the $3$-dimensional factors and a line $ux$ such that $(ux)^\Alpha{\,\ne\,}ux$ and 
$\cE{\,=\,}\langle x^{\Alpha\Beta},u,v,av\rangle{\ledot\hskip-4pt}\cP$ as in Case 2).  
Let again $\Chi{\,=\,}(\Cs{}\Alpha\Beta)^1$.
If $\dim\Delta{\,\ge\,}25$, 
then $\dim\Chi{\,=\,}19$,\: $\Chi_x|_\cE{\,=\,}\1$,\: $\dim\Chi_x{\,\le\,}3$,\: $\dim x^\Chi{\,=\,}16$, and 
$(vx)^\Chi$ is open in the pencil $\frL_v$  by \cite{cp} 96.11. 
As $x{\,\notin\,}av$ is an arbitrary point of $ux\sm\{u\}$, it follows that $\Delta$ is transitive on 
$\frS{\,=\,}\frL_v\sm\{av,uv\}$. The exact homotopy sequence (\cite{cp} 96.12) will be applied to 
this action. Let $L$  be some line in $\frS$, and consider the part 
$$\dots\to\pi_7\Delta\to\pi_7\frS\to\pi_6\Delta_L\to\dots$$ of the homotopy sequence.  
Note that $\frS$ is homotopy equivalent ($\simeq$) to $\sS_7$ by 2.1. Choose maximal compact 
subgroups $\Phi$ in $\Delta$ and $\Kappa$ in $\Delta_L$ such that $\Kappa{\,\le\,}\Phi$. 
Then $\Delta{\,\simeq\,}\Phi$ and
$\Delta_L{\,\simeq\,}\Kappa$ by the Mal'cev-Iwasawa theorem \cite{cp} 93.10, and we get 
 an exact sequence $$\dots\to\pi_7\Phi\to\pi_7\sS_7\to\pi_6\Kappa\to\dots.$$ 
 It is well-known that 
$\pi_7\sS_7{\,\cong\,}\ZZ$\: (see, e.g., \cite{sp} 7.5.6 or \cite{br} II.16.4). A maximal compact 
subgroup of $\Gamma$ is isomorphic to $\Spin3\RR$, hence homeomorphic to $\sS_3$. 
At most one of the $3$-dimensional factors has a non-trivial maximal compact subgroup, it is  
homeomorphic to $\sS_3$ or to $\sS_1$. If $k$ is odd, then all homotopy groups 
$\pi_q\sS_k$ with $q{\,\ne\,}k$ are finite (\cite{sp} 9.7.7) and $\pi_q\sS_1{\,=\,}0$  for $q{\,\ne\,}1$ 
(\cite{sp} 7.2.12); in fact, $\pi_7\sS_3{\,\cong\,}\ZZ_2$ and $\pi_6\sS_3{\,\cong\,}\ZZ_{12}$.
Therefore both $\pi_7\Phi$ and $\pi_6\Kappa$ are finite. Exactness would force $\pi_7\sS_7$ to be finite, a contradiction. This completes the proof of step (j) and shows that $\Delta$ is almost 
simple or has some almost simple factor of dimension ${\hskip-2pt\ge\,}14$. \\
(k) {\it If each almost simple factor of $\Delta$ has dimension $3$ or $14$, then 
$\dim\Delta{\,\le\,}20$\/}. \\
By step (f), we may assume that $\dim\Gamma{=\,}14$, and then $\rk\Gamma{=\,}2$. Lemma 8.1
implies that $\Gamma$ contains a planar involution $\iota$. For the compact group $\Gamma$ of 
type ${\rm G}_2$ it has been stated in 2.5 that $\Cs\Gamma\iota$ is semi-simple and 
$\dim\Cs\Gamma\iota{\,=\,}6$. Now let $\Gamma$ 
be a non-compact group of type ${\rm G}_2$. If $\Gamma$ is strictly simple, then $\Gamma$ has 
a subgroup $\Phi{\,\cong\,}\SO4\RR$; by Lemma 8.1  we may choose $\iota$ as the central involution of $\Phi$; 
the double covering group of $\Gamma$ contains $\tilde\Phi{\,\cong\,}(\Spin3\RR)^2$, and the central involution of one of the two factors of $\tilde\Phi$ can be taken as $\iota$, so that in each case 
$\dim\Cs\Gamma\iota{\,\ge\,}6$. Now the following lemma shows  $\dim\Psi{\,\le\,}6$ and the assertion is proved. \\
($\ell$) {\tt Lemma.} {\it If $\iota$ is a \emph{planar} involution in $\Gamma$ and if $\Chi$ is a 
maximal semi-simple subgroup of $\Cs\Gamma\iota$, then $\dim\Chi\Psi|_{\cF_\iota}{\,\le\,}10\,$\/}
(by  \cite{sz1} 6.1) {\it and $\dim\Chi{\,+\,}\dim\Psi{\,\le\,}10{+}3$\/}. \\
(m) {\it If each almost simple factor of $\Delta$ has dimension $15,16$, or $3$, then 
$\dim\Delta{\,\le\,}22$\/}. \\
In fact, if $\dim\Gamma{\,\in\,}\{15,16\}$, then there is a planar involution $\iota{\,\in\,}\Gamma$  
such  that $\Cs\Gamma\iota$ contains a $6$-dimensional semi-simple group, see cases 3) and 4) of 
step (j). Lemma ($\ell$) implies $\dim\Delta{\,\le\,}16{\,+\,}\dim\Psi{\,\le\,}22$. \\
(n) {\it If $\Delta$ has an  almost simple factor $\Gamma\!$ of  dimension $20$, then 
$\Delta{\,\cong\,}\Sp4\CC$\/}. \\
By the last part of 2.10, the simple group $\SO5\CC$ cannot act on a plane. Its covering group 
$\Spin5\CC{\,\cong\,}\Sp4\CC$ contains a planar involution $\iota{\,=\,}{\rm diag}(-1,-1,1,1)$ 
centralized by the $12$-dimensional semi-simple group $(\Sp2\CC)^2{\,\cong\,}(\SL2\CC)^2$.
Hence  $\Delta{\,=\,}\Gamma$ by Lemma ($\ell$). (In fact, the group $(\SL2\CC)^2$ cannot act on 
$\cF_\iota$ and case (n) is impossible). \\
(o) {\it Suppose that $\Delta$ has an almost simple factor $\Gamma$ of type ${\rm B}_3$. 
Then $\dim\Delta{\,\le\,}24$\/}. \\
By 2.\,10 and 11, none of the simple groups $\Opr7(\RR,r)$ can act on $\cP$. The central involution 
of the compact group $\Gamma{\,\cong\,}\Spin7\RR$ is a reflection with axis $av$, and the action 
of $\Gamma$ on $av\sm\{v\}$ is equivalent to the linear action of $\SO7\RR$ on $\RR^8$, see 2.15. 
Hence $\Gamma$ fixes a circle on $av$. For $z{\,\in\,}uv\sm\{u,v\}$, Stiffness implies that
 $\Gamma_{\hskip-2pt z}{\,=\,}\Lambda{\,\cong\,}\Gtwo$, that $\cF_\Lambda$ is a flat subplane, and that $\Psi{\,=\,}\Cs{}\Gamma$ acts almost effectively on $\cF_\Lambda$. Now $\Psi|_{\cF_\Lambda}$
is solvable by \cite{cp} 33.8.  Hence~$\Delta{\,=\,}\Gamma$. \\
 Next, let $\Gamma{\,\cong\,}\Spin7(\RR,1)$. As $\Gamma$ contains a $3$-torus, there are (at least) 
 $6$ non-central planar involutions in $\Gamma$. By the covering homomorphism, they are mapped 
 to involutions in $\Opr7(\RR,1)$, up to conjugacy to diagonal matrices with entries $\pm1$ and 
determinant $1$. In the centralizer of each of these matrices there is a semi-simple group of 
dimension ${\!\ge\,}9$, which is  covered by a semi-simple group in the centralizer  of  a planar involution in $\Gamma$.  Consequently, Lemma ($\ell$) shows $9{\,+\,}\dim\Psi{\,\le\,}13$,\: 
$\dim\Psi{\,\le\,}3$, and $\dim\Delta{\,\in\,}\{21,24\}$. The same arguments apply in the cases 
$r{\,>\,}1$ even if there are more possibilities than just a double covering of $\Opr7(\RR,r)$. \\
(p) {\it If $\Delta$ has an almost simple factor $\Gamma$ of type ${\rm C}_3$, then 
$\dim\Delta{\,=\,}21$ and $\Delta{\,=\,}\Gamma$\/}. \\
The planar involution ${\rm diag}(-1,-1,1)$ in $\U3(\HH,r)$ has a semi-simple centralizer
$\U2\HH{\times}\HH'$ of dimension $13$. From Lemma ($\ell$) it follows that $\Delta{\,=\,}\Gamma$, and the claim holds 
for the $4$ groups (P\hskip-.5pt)$\hskip-2pt\U3(\HH,r)$ with $r{\,=\,}0,1$. The symplectic group 
$\Gamma{\,=\,}\Sp6\RR$ has a maximal compact subgroup $\U3\CC$. The involution 
$\iota{\,=\,}{\rm diag}(-1,-1,1){\,\in\,}\SU3\CC$ is planar, its centralizer contains the semi-simple group 
$\Ypsilon{\,=\,}\Sp4\RR{\times}\Sp2\RR$ of dimension $13$. The element $\iota$ is mapped to 
a non-central involution in $\PSp6\RR$, and it belongs to each of the infinitely many covering groups 
of $\Gamma$ because $\SU3\CC$ is simply connected. In any case, the centralizer of $\iota$ 
is locally isomorphic to $\Ypsilon$. Again $\dim\Ypsilon{\,+\,}\dim\Psi{\,\le\,}13$ and 
$\Delta{\,\circeq\,}\Gamma$. \\ 
(q) Each group $\Gamma{\,\circeq\,}\SU5(\CC,r)$ contains a non-central and hence  planar involution corresponding to  $\iota{\,=\,}{\rm diag(-1,-1,-1,-1,1)}$, and $\iota$ is centralized by the $15$-dimensional group  $\Ypsilon{\,=\,}\SU4(\CC,r{-}1)$, but planarity implies  $\dim\Ypsilon{\,\le\,}10$. 
Analogously, the simply connected covering group of $\SL5\RR$ can be excluded. \\
(r) Only one possibility remains: $\Gamma{\,=\,}\Delta$ is an almost simple group of dimension
$28$, in fact, an orthogonal group of type D$\hskip-1.5pt_4$, see \cite{cp} 94.33 and 2.16 above. 
From 2.10 and 2.11 it follows that $\Delta$ is a proper covering of a  group $\POpr8(\RR,r)$
different from $\Opr8(\RR,r)$.  By step (b) there is at most one central involution. As $\Spin8\RR$ 
has center $(\ZZ_2)^2$, the group  $\Delta$ is not compact and $r{\,\ne\,}0$. If $r{\,=\,}4$, 
a maximal compact subgroup of the simply connected covering group  $\tilde\Delta$ is isomorphic 
to  $(\Spin4\RR)^2{\,\cong\,}(\Spin3\RR)^4$. Hence each non-trivial element in the center of 
$\Delta$  is an involution, and $\Delta$ would contain $(\SO3\RR)^2$ contrary to 2.11. \\
(s) The cases $0{\,<\,}r{\,<\,}3$ can be excluded by similar arguments as in step (p):\:   
$\Delta$ contains the simply connected double covering $\Phi{\,=\,}\Spin6\RR$ of $\SO6\RR$.  
Choose a $1$-torus $\Theta{\,<\,}\Phi$ which is disjoint from the kernel of the covering map 
$\kappa{\,:\,}\Delta{\,\to\,}\POpr8(\RR,r)$. Up to conjugation, the involution $\iota{\,\in\,}\Theta$ 
is mapped onto ${\rm diag}(-1,-1,1,...,1){\,\in\,}\Opr8(\RR,r)$. Hence $\Cs\Delta\iota$ contains 
a $15$-dimensional semi-simple group $\Ypsilon{\,\circeq\,}\Opr6(\RR,r)$, but   
$\dim\Ypsilon{\,\le\,}10$ by \cite{sz1} 6.1.  \\
(t) Finally, let $\Delta$ be a double covering of $\Opr8(\RR,3)$, and consider the covering map of a $2$-torus $\Theta$ in the subgroup $\Spin5\RR$ of $\Delta$ into $\SO5\RR$. Up to conjugation, 
some non-central (and hence planar) involution $\iota{\,\in\,}\Delta$ is contained in $\Theta$ and  is mapped onto 
${\rm diag}(-1,-1,1,1,1){\,\in\,}\SO5\RR$, and $\Cs\Delta\iota$ is locally isomorphic to 
$\Opr6(\RR,3)$, again a contradiction. \qed
\par\medskip
{\bf 8.3 Normal torus.} {\it Suppose that $\cF_\Delta$ is a double flag  and that $\Delta$ has 
a one-dimensional compact connected normal  subgroup $\Theta{\,\triangleleft\,}\Delta$. Then 
$\dim\Delta{\,\le\,}30$\/}.
\par\smallskip
{\tt Proof.}  (a) We may assume that $\dim\Delta{\,\ge\,}27$. Then $\Delta$ is a Lie group by 2.3, 
and $\Theta$ is a central torus. If $\Lambda{\,<\,}\Delta$ and if $\cF_\Lambda$ is a subplane, then $\Theta$ acts non-trivially on $\cF_\Lambda$\: (or else 
$\cF_\Lambda{\,\le\,}\cF_\Theta{\,=\,}\cF_{\hskip-1.5pt\Theta}^{\hskip1pt\Delta}{\,<\,}\cP$ and 
$\dim\Delta{\,\le\,}18$ by Stiffness). Hence $\cF_\Lambda$ is not flat, $\dim\cF_\Lambda{\,\ge\,}4$,\:
$\dim\Lambda{\,\le\,}8$,\: $\Delta{\hskip1pt:\hskip1pt}\Lambda{\,\le\,}8{+}16$, and $\dim\Delta{\,\le\,}32$. In fact, $\Theta{\,\le\,}\Delta_{[av]}$: 
let $a,c{\,\in\,}K{\,=\,}av\sm\{v\}$ and $z{\,\in\,}S{\,=\,}uv\sm\{u,v\}$, and consider  
$\Gamma{\,=\,}\Delta_{a,c}$ and $\Lambda{\,=\,}\Gamma_{\hskip-2pt z}$. 
If $a{\,\ne\,}c{\,\in\,}a^\Theta$, then $\Gamma {\,= \,}\Delta_a$ and 
$\Delta{\hskip1pt:\hskip1pt}\Lambda{\,\le\,}16$ contrary to the assumption. \\
(b) Suppose that  $\dim\Delta{\,=\,}32$.  Then $\Delta_z$ is doubly transitive on $K$,\: $\Gamma$
is transitive on $S$,\: $\dim\Gamma{\,=\,}16$,  and $\Lambda{\,\cong\,}\SU3\RR$. If $\Phi$ is a maximal compact subgroup of $\Gamma$, then $\Phi'{\,\cong\,}\SU4\CC$ by 2.17(d).
Now $\Theta\Phi'{\,\le\,}\Phi$ and $\Phi{\hskip1pt:\hskip1pt}\Lambda{\,=\,}8$, but then $\Phi$ would be transitive on $S$, which  is impossible.  Therefore $\dim\Gamma{\,<\,}16$ and 
$\dim\Delta{\,<\,}32$. \\
(c) Let $\dim\Delta{\,=\,}31$. Then $\dim\Gamma{=\,}15$ by step (b), and $\Delta$ is 
doubly transitive on~$K$\: (notation as in (a)\,). Put $\nabla{\,=\,}\Delta_a$ and 
$\hat\nabla{\,=\,}\nabla|_K{\,\cong\,}\nabla/\Nu$. The kernel $\Nu{\,=\,}\Delta_{[av]}$ consists of 
homologies with center $u$ and $\Theta{\,\trianglelefteq\,}\Nu$ by step (a).
Let $\Kappa$ be a maximal compact connected subgroup of $\Nu$. Then 
$\Nu{\hskip1pt:\hskip1pt}\Kappa{\,\le\,}1$ by \cite{cp} 61.2, and $\Kappa$ has torus rank 
$\rk \Kappa{\,=\,}1$\; (use \cite{cp} 55.35). Therefore $\Kappa{\,=\,}\Theta$,\: $\dim\Nu{\,\le\,}2$, and 
$21{\,\le\,}\nabla{\hskip1pt:\hskip1pt}\Nu{\,<\,}23$. The structure of doubly transitive groups
has been determined by Tits, cf. \cite{cp} 96.16. As $K$ is not compact, $\hat\nabla$ is a transitive  subgroup of $\GL8\RR$; these groups are described explicitly in \cite{Vl}, 
see also \cite{cp} 96.\,19--22. In particular,  $\hat\nabla{\,=\,}\Eta\hat\Ypsilon$ is a product of an  almost simple Lie group $\hat\Ypsilon$ and a subgroup $\Eta$ of its centralizer; moreover 
$\Eta{\,\le\,}\HH^{\times}$,  a maximal compact subgroup of $\hat\Ypsilon$ is transitive on the 
$7$-sphere formed by the rays in $\RR^8$, and  $16{\,<\,}\dim\hat\Ypsilon{\,<\,}23$.  In this 
dimension range, only three almost simple  groups have an $8$-dimensional irreducible representation,  viz. the groups  $\Sp4\CC$,\:   $\Spin7\RR$, and  $\Spin7(\RR,3)$,  
see \cite{cp} 96.10. The last one does not act transitively on $\sS_7$, the other two are simply connected. By \cite{cp} 94.27 it follows that $\nabla$ has a subgroup $\Ypsilon$ which is 
isomorphic to one of $\Sp4\CC$ or $\Spin7\RR$. First, let $\Ypsilon$ be compact. Then 
$\nabla{\,\cong\,}e^{\RR}{\cdot}\Spin7\RR{\cdot\hskip.5pt}\Theta$ and 
$\nabla_{\hskip-2pt c}{\,=\,}\Lambda{\times}\Theta$ with $\Lambda{\,\cong\,}\Gtwo$\: (see 2.15). 
Now $\cF_\Lambda$ is a flat subplane and $\Theta$ acts as a group of homologies on 
$\cF_\Lambda$. This  contradicts \cite{cp} 32.17 and shows that $\Ypsilon{\,\cong\,}\Sp4\CC$.
The fixed elements of the involution $\iota{\,=\,}(-1,-1,1,1)$ on $av$ form a $4$-sphere. Hence 
$\cF_\iota{\;\ldot}\cP$. Choose $c$ and $z$ as before, but in $\cF_\iota$. Then 
$\cF_\Lambda{\;\ldot}\cF_\iota$, and Stiffness 2.6(c,\^e) implies $\Lambda{\,\cong\,}\SU3\CC$.
This group, however, is not contained in the  maximal compact subgroup $\U2\HH$ of $\Ypsilon$. 
Alternatively, $\Cs\Ypsilon\iota$ has a subgroup $(\SL2\CC)^2{\,\cong\,}(\Sp2\CC)^2$ contrary to 
\cite{sz1} 6.1. \qed
\par\medskip
{\bold Remark.} {\it Under the assumptions of {\rm8.3}, suppose that $\dim\Delta{\,\ge\,}20$. 
Then $\Delta$ is a Lie group or the center $\Zeta$ of $\Delta$ consists of homologies in 
$\Delta_{[u,av]}$\/}.
\par\smallskip
{\tt Proof.} (a) Assume that $a^\Zeta{\,\ne\,}a$ and that $\Delta$ is not a Lie group.  Let $\Nu$ be a compact $0$-dimensional normal subgroup such that $\Delta/\Nu$ is a Lie group\; (cf. \cite{cp} 93.8). Put $\nabla{\,=\,}\Delta_a$ and note that $\Theta{\,\le\,}\Zeta$ and that 
$\nabla|_{a^\Zeta}{\,=\,}\1$. 
There is some point $p{\,\in\,}au\sm\{a,u\}$ such that  $p^\Theta{\,\ne\,}p$\; (or else 
$\Theta{\,\le\,}\Delta_{[v,au]}$, which contradicts $a^\Zeta{\,\ne\,}a$). The group 
$\Lambda{\,=\,}(\Delta_p)^1$ fixes $a^\Zeta$,  hence a quadrangle,  
$\Delta{\hskip1pt:\hskip1pt}\Lambda{\,\le\,}16$,
and $\cE{\,:=\,}\cF_\Lambda{\,=\,}\cE^\Theta$ is a connected proper subplane.
From  $\Theta|_\cE{\,\ne\,}\1$ it follows that $\dim\cE{\,\ge\,}4$\; (use \cite{cp} 32.21{b} and 17). \\
(b) If $\dim\cE{\,=\,}4$, then $\Nu|_\cE{\,=\,}\Nu/\Xi$ is a Lie group by  \cite{cp} 71.2, and 
$\Xi{\,\ne\,}\1$. Hence $\cE{\,\le\,}\cF_\Xi{\,<\,}\cP$, and $\cF_\Xi$ is $\Delta$-invariant. 
We have $\Delta{\hskip1pt:\hskip1pt}\nabla{\,\le\,}8$,\;  
$\nabla{\hskip1pt:\hskip1pt}\Lambda{\,=\,}\dim p^\nabla$, and $p^\nabla{\,\subseteq\,}\cF_\Xi$.
Consequently $\Delta{\hskip1pt:\hskip1pt}\Lambda{\,\le\,}12$ and $\dim\Lambda{\,\le\,}8$.
If $\cF_\Xi{\,=\,}\cE$, then $\dim p^\nabla{\,\le\,}2$ and $\dim\Delta{\,\le\,}18$; if 
$\cF_\Xi{\;\ldot}\cP$ and  $\dim\Lambda{\,=\,}8$, then $\Lambda{\,\cong\,}\SU3\CC$
 and the lines of $\cF_\Xi$ are $4$-spheres, but then $\Lambda$ cannot act on 
$\cF_\Xi$, see $(\dagger)$. Therefore $\dim\Delta{\,<\,}20$.  
If $\cE{\;\ldot}\cP$ and $\Lambda\smcap\Nu{\,=\,}\1$, then $\Lambda$ is a Lie group and 
$\dim\Lambda{\,\le\,}3$ by Stiffness, hence $\dim\Delta{\,\le\,}16{+}3{\,<\,}20$. If 
$\1{\,\ne\,}\zeta{\,\in\,}\Lambda\smcap\Nu$, then $\cE{\,=\,}\cF_\zeta{\,=\,}\cE^\Delta$ and 
$\Delta{\hskip1pt:\hskip1pt}\Lambda{\,\le\,}8$. In this case $\dim\Delta{\,\le\,}15$. \qed
\par\medskip
The following  theorems have been proved in \cite{hs}:
\par\smallskip
{\bf 8.4.} {\it If $\dim\Delta{\,\ge\,}33$ and if $\cF_\Delta{\,=\langle u,v,av\rangle\,}$ is a double flag, 
then the translation group $\Tau{\hskip1pt=\,}\Delta_{[v,uv]}$ is transitive and 
$\Phi{\,=\,}(\Delta_a)'{\,\cong\,}\Spin8\RR$. In particular, $\dim\Delta{\,\ge\,}36$. 
If $\dim\Delta{\,\ge\,}38$, then the plane is classical\/}.
\par\medskip
{\bf 8.5 Distorted octonions.}  {\it Let $(\RR,+,\ast,1)$ be a topological  Cartesian field with unit
element $1$ such that  $\:(-r){\,\ast\,}s{\,=\,} {\,-\,} (r{\,\ast\,}s){\,=\,}r{\,\ast\,}(-s)$ holds identically. 
Define a  new multiplication  on the octonion algebra  $(\OO,+, \: )$ by
$ a {\,\circ\,} x {\,=\,} {|a|{\ast}|x|\,|ax|^{-1}}\,  a\hskip1pt x $ for $a,x{\,\ne\,} 0$  and
$0{\,\circ\,}x {\,=\,} a{\,\circ\,}0 {\,=\,} 0$.
Then the {\emph distorted octonions} $(\OO,+,\circ)$ form a topological Cartesian field\/}.
\par\medskip
{\bf 8.6 Theorem.} {\it A plane $\cP$ can be coordinatized by distorted octonions if, and only if, 
$\cP$ has the properties of Theorem\/} 8.4. 
\par\bigskip 
{\Bf 9. Fixed triangle}
\par\medskip 
Let $\cF_\Delta{\,=\,}\langle a,u,v\rangle$ be a triangle. If $\cP$ is a translation plane with 
translation axis $uv$, then $\Delta$ has a compact subgroup $\Phi$ of codimension 
$\Delta{\,:\,}\Phi{\,\le\,}2$, see H\"ahl \cite{Ha5}, \cite{cp} 81.8. There seems to be no way to extend this basic result to general $16$-dimensional planes; cf. 9.1, however. 
\par\medskip
{\bf 9.0 Lie.}	{\it If  $\cF_\Delta{\,=\,}\langle a,u,v\rangle$ is a triangle, and if
$\dim\Delta{\,\ge\,}26$, then $\Delta$ is a Lie group\/}.
\par\smallskip
By \cite{psz}, only the case $\dim\Delta{\,=\,}26$ requires a {\tt proof.} Suppose that 
$\Delta$ is not a Lie group, and consider the action of 
$\Delta$ on the sides $S_\nu$ of the triangle (without the vertices). 
If $\Delta$ is transitive on $S_\nu$, then $\Delta|_{S_\nu}{\,=\,}\Delta/\Kappa_\nu$ is a Lie group, 
see \cite{cp} 53.2. Suppose that this happens for two distinct indices $\mu,\nu$. It follows that  
$\Kappa_\nu$ is a Lie group, since $\Kappa_\mu\smcap\Kappa_\nu{\,=\,}\1$, and then $\Delta$ 
itself is a Lie group by \cite{cp} 94.3(d). For the same reason, all orbits on two sides of the triangle 
have dimension~${<\,}8$, and $\dim x^\Delta{\,\le\,}14$ for each point 
$x{\,\notin\,}au{\smcup}av{\smcup}uv$. Put $\Lambda{\,=\,}(\Delta_x)^1$. Then 
$\Delta{\hskip1pt:\hskip1pt}\Lambda{\hskip1pt\le\hskip1pt}14$ and $\dim\Lambda{\,>\,}11$. Stiffness implies 
$\Lambda{\,\cong\,}\Gtwo$, and representation of $\Lambda$ on the Lie algebra 
$\frl\hskip1pt\Delta$ shows $\dim\Cs{}\Lambda{\,\in\,}\{5,12\}$. Again by Stiffness, $\Cs{}\Lambda$ 
acts almost effectively on the flat plane $\cF_\Lambda$. This contradicts the assumption on 
$\cF_\Delta$. \qed
\par\medskip
{\bf Remark.} {\it If  $\cF_\Delta{\,=\,}\langle a,u,v\rangle$ is a triangle, if $\dim\Delta{\,\ge\,}19$, 
and if $x{\,\notin\,}au{\smcup}av{\smcup}uv$, then $\Delta_x$ is a Lie group\/}.
\par\smallskip
{\tt Proof.} We may assume that $\Delta$ is not a Lie group. Again there exist arbitrarily small compact central  subgroups $\Nu{\,\le\,}\Delta$ of dimension $0$ such that $\Delta/\Nu$ is a Lie group, see \cite{cp} 93.8. We will show that $\Nu$ acts freely on 
$P\sm (au{\smcup}av{\smcup}uv)$. Suppose that $a,u,v,x$ is a non-degenerate quadrangle, 
and that $x^\zeta{\,=\,}x$ for some $\zeta{\,\in\,}\Nu\sm\{\1\}$. Then 
$\zeta|_{\langle x^\Delta\rangle}{\,=\,}\1$ and 
$\cD{\,=\,}\langle x^\Delta, a,u,v\rangle{\,=\,}\cD^\Delta$ is a proper subplane; in fact,
$\cD{\;\ldot}\cP$\; (or else $\dim x^\Delta{\,\le\,}4$, and $\Lambda{\,=\,}(\Delta_x)^1$ would have 
dimension $\dim\Lambda{\,>\,}14$, which contradicts Stiffness). 
Put $\Delta|_\cD{\,=\,}\Delta/\Kappa$.  By 
\cite{sz1} 1.7 we have 
$\Delta{\hskip1pt:\hskip1pt}\Kappa{\,\le\,}11$ and $\dim\Kappa{\,\ge\,}8$, but  
the stiffness property 2.6(b) implies $\dim\Kappa{\,\le\,}7$. This contradiction shows that 
$\Delta_x{\smcap}\Nu{\,=\,}\1$. \qed

\par\medskip
{\bold 9.1 Semi-simple groups.} {\it  If $\Delta$ is semi-simple, if $\cF_\Delta$ is a triangle $a,u,v$, and if $\dim\Delta{\,\ge\,}28$,  then $\Delta{\,\cong\,}\Spin8\RR$\/}.
\par\smallskip
{\tt Proof.} (a) Choose a point $e$ such that $a,e,u,v$ is a non-degenerate quadrangle, and consider the connected component $\Lambda$ of $\Delta_e$. Then 
$\Delta{\hskip1pt:\hskip1pt}\Lambda{\,\le\,}16$,\: $\dim\Lambda{\,>\,}11$, and 
$\Lambda{\,\cong\,}\Gtwo$ by Stiffness. In particular, $\dim\Delta{\,\le\,}30$. \\
(b) {\it $\Delta$ is  almost simple\/}. If not, then some proper  factor $\Gamma$ of $\Delta$ 
contains an 
isomorphic copy  $\hat\Lambda$ of $\Lambda$. Repeated application of 2.15(a) and Stiffness show that  $\cF_{\hat\Lambda}$ is a flat subplane and that $\Psi{\,=\,}(\Cs\Delta\Gamma)^1$ acts 
effectively on $\cF_{\hat\Lambda}$.  By \cite{cp} 33.8 the group $\Psi$ is solvable, in fact 
$\dim\Psi{\,\le\,}2$\: (see \cite{cp} 32.10, cf. also 33.10);  hence $\Psi{\,=\,}\1$  and $\Delta{\,=\,}\Gamma$, a contradiction. \\
(c) If $\dim\Delta{\,=\,}30$, then a maximal compact subgroup of $\Delta$ is locally isomorphic to 
$\SU4\CC$ and does not contain $\Gtwo$. Hence $\dim\Delta{\,=\,}28$, and $\Delta$ is a group of 
type ${\rm D}_4$ by 2.16. From 2.10 and $\Lambda{\,<\,}\Delta$ it follows that 
$\Delta{\,\cong\,}\Spin8(\RR,r)$ with $r{\,\le\,}1$. Suppose that $r{\,=\,}1$. Then $\Delta$ has 
a maximal compact subgroup $\Phi{\,\cong\,}\Spin7\RR$ such that $\Lambda$ is contained in $\Phi$. 
Note that $\Lambda$ is even maximal in $\Phi$\; (if not, then the action of $\Lambda$ on 
$\frl\Phi$ shows that $\Chi{\,=\,}\Cs\Phi\Lambda$ is $7$-dimensional and $\Phi{\,=\,}\Lambda\Chi$, but then $\Lambda$ would be normal in $\Phi$). 
The central involution $\sigma{\,\in\,}\Phi$ is a reflection, say with axis $uv$. Therefore $\Delta$ acts effectively on $av$.
Let $c{\,=\,}eu\smcap av$ and $\Gamma{\,=\,}(\Delta_c)^1$. Each orbit of $\Phi$ on $av\sm\{a,v\}$ 
is a $7$-sphere (cf. 2.15), and $\Phi$ acts linearly on $av\sm\{v\}{\,\approx\,}\RR^8$, see \cite{cp} 96.36. The stabilizer $\Phi_c{\,=\,}\Lambda{\,\cong\,}\Gtwo$ is a maximal compact subgroup of 
$\Gamma$\; (use the Mal$'$cev-Iwasawa theorm).
We have $\Delta{\hskip1pt:\hskip1pt}\Gamma{\,\ge\,}\dim c^\Phi$ and 
$\Gamma{\hskip1pt:\hskip1pt}\Lambda{\,\in\,}\{6,7\}$. Stiffness implies 
$\dim(\Cs\Gamma\Lambda){\,\le\,}1$. As the action of $\Lambda$ on 
(the vector space underlying) the Lie algebra $\frl\hskip1pt\Gamma$ is completely reducible, it 
follows from \cite{cp} 95.10 that $\dim\Gamma{\,=\,}21$. Hence $\dim c^\Delta{\,=\,}7$, and then 
$c^\Delta{\,=\,}c^\Phi$. In particular, $\Delta_c{\,=\,}\Gamma$ is connected. It follows that the map  
$\kappa{\,:\,}\delta{\,\mapsto\,}c^\delta$ on $\Delta$ induces a homeomorphism 
$\Delta/\Gamma{\,\approx\,}c^\Phi{\,\approx\,}\sS_7$, see \cite{cp} 96.9(a). Moreover, 
$\Lambda{\,<\,}\Gamma$ implies that $\Gamma$ is neither almost simple nor semi-simple. Consequently, $\sqrt\Gamma{\,=\,}\Xi{\,\cong\,}\RR^7$ and $\Gamma{\,=\,}\Xi{\rtimes}\Lambda$ is given by the action of $\Aut\OO$ on the pure octonions.  The same arguments may be applied 
to $b{\,=\,}ev\smcap au$ instead of $c$.  In particular, $b^\Delta{\,=\,}b^\Phi{\,=\,}S{\,\approx\,}\sS_7$.
The vector group $\Xi$ acts freely on $S$ because each stabilizer $\Xi_s$ wwith $s{\,\in\,}S$ 
fixes a quadrangle and is compact according to step (a). By the open mapping theorem \cite{cp} 
51.19 or 96.11, each orbit $s^\Xi$ is open in $S$. Hence $b^\Xi{\,=\,}S$ and 
$\xi{\,\mapsto\,}b^\xi{\,:\,}\Xi{\,\to\,}S$ is a homeomorphism, but obviously this is not true. \qed

\par\smallskip
{\bold Remark.} {\it If $\Delta$ is semi-simple and $\cF_\Delta{\,=\,}\langle a,u,v\rangle$ 
is a triangle, then $\dim\Delta{\,\ne\,}25$\/}. 
\par\smallskip  
{\tt Proof.} Assume that $\dim\Delta{\,=\,}25$. This number being odd, $\Delta$ has 
either a $15$-dimensional factor or at least one factor of dimension $3$. \\
(a) In the first case, $\Delta{\,=\,}\Gamma\Psi$ has exactly two factors. Let $\dim\Gamma{\,=\,}15$.    
If $\Gamma$ is transitive on $au\sm\{a,u\}{\,\ni\,}x$, then $\Psi_{\hskip-1pt x}$  is a group of homologies with axis $au$, and so is $\Psi$ because $\Psi$ is almost simple,  but then 
$\dim\Psi{\,\le\,}8$, a contradiction. Analogously, $\Gamma$ has an orbit 
$y^\Gamma{\subset\,}av\sm\{a,v\}$ of dimension ${<\,}8$. Hence there is a point 
$p{\,\notin\,}au\smcup av\smcup uv$ such that $\Gamma_{\hskip-2pt p}{\,\ne\,}\1$, and 
$\cE{\,=\,}\cF_{\Gamma_{\hskip-1pt p}}$ is a proper $\Psi$-invariant  subplane. Consequently, 
$\dim p^\Psi{\,\le\,}8$ and $\dim\Psi_{\hskip-1pt p}{\,\ge\,}2$. It follows that 
$\cH{\,=\,}\cF_{\Psi_{\hskip-1pt p}}{<\,}\cP$, and $\Gamma|_\cH$ is too large. \\
(b) Suppose that $\Delta$ has a $10$-dimensional almost simple factor $\Psi$. Then its 
complement $\Gamma{\,=\,}(\Cs{}\Psi)^1$ has dimension $15$. In step (a), the fact that 
$\Gamma$ is almost simple has been used only at the very end. Therefore there is again a subplane 
$\cH{\,=\,}\cH^\Gamma{<\,}\cP$. Put $\Gamma|_\cH{\,=\,}\Gamma/\Kappa$. 
If $\cH{\;\ldot}\cP$, then $\Kappa$  is compact and semi-simple, hence a Lie group, and 
$\Kappa{\,\cong\,}\Spin3\RR$ by Stiffness. On the other hand, 
$\Gamma{\hskip1pt:\hskip1pt}\Kappa{\,\le\,}9$, see \cite{sz9} 6.1 or 
\cite{sz1} 7.3. Hence 
$\dim\cH{\,\le\,}4$, and then $\Gamma{\hskip1pt:\hskip1pt}\Kappa{\,\le\,}3$ by 2.14. Now Stiffness implies   $\Kappa{\,\cong\,}\Gtwo$, which is impossible. \\
(c) If $\Delta$ has an almost simple factor $\Gamma$ of dimension $16$, then a maximal 
compact subgroup of $\Gamma$ is locally isomorphic to $\SU3\CC$, and 2.17 implies that 
$\Gamma$ cannot be transitive on one of the sides of the fixed triangle (without its vertices).
Hence there is again a point $p$ such that $\Gamma_{\hskip-2pt p}{\,\ge\,}2$ and
$\cE{\,=\,}\cF_{\Gamma_{\hskip-1pt p}}$ is a proper $\Psi$-invariant  subplane. Now 
$\dim p^\Psi{\,\le\,}8$, and $ \cH{\,=\,}\cF_{\Psi_{\hskip-1pt p}}{<\,}\cP$, but then $\dim\Gamma|_\cH{\,\le\,}11$ 
by \cite{sz9} 6.1 or \cite{sz1} 1.7. \\
(d) The case that  $\Delta$ has two $8$-dimensional factors is similar but somewhat more complicated.  Put $D{\,=\,}P\sm(au\smcup av\smcup uv)$, let $\dim\Gamma{\,=\,}8$, and denote the product of the other factors by $\Psi$. Then $\dim\Psi{\,=\,}17$,\, $\Psi_{\hskip-1pt q}{\,\ne\,}\1$ for $q{\,\in\,}D$, and  $\cD{\,=\,}\cF_{\Psi_{\hskip-.5pt q}}{\,=\,}\cD^{\hskip.5pt\Gamma}$ is a proper subplane. 
Either $\Gamma$ is sharply transitive on $\cD\smcap D$, or $\Gamma_{\hskip-2pt p}{\,\ne\,}\1$ 
for some $p$ and $\cE{\,=\,}\cF_{\Gamma_{\hskip-1pt p}}{\,=\,}\cE^\Psi{\,<\,}\cP$. In the first case, 
let $\Phi$ be a maximal compact subgroup of $\Gamma$. Note that $\dim\Phi{\,<\,}8$ and that 
$\cD\smcap D$ is homotopy equivalent to $\sS_3^{\;2}$. It follows that $\Phi'{\,\circeq\,}\Spin3\RR$, and  the exact homotopy sequence \cite{cp} 96.12 yields 
$0\to\pi_3\Phi'\to\pi_3\sS_3^{\;2}{\hskip1pt=\hskip1pt}\ZZ^2\to0$, a contradiction.
Hence $\cE{\,<\,}\cP$. For each possible dimension 
of $\cE$ it turns out that $\Psi$ is too large: write $\Psi|_\cE{\,=\,}\Psi/\Kappa$. If 
$\dim\cE{\,\le\,}4$, then $\Psi{\,=\,}\Kappa$\; (see \cite{cp} 33.8 and 2.14 above), but 
$\dim\Kappa{\,\le\,}14$ by Stiffness; if $\cE{\;\ldot}\cP$, then $\dim\Kappa{\,\le\,}3$ and 
$\Psi{\hskip1pt:\hskip1pt}\Kappa{\,\ge\,}14$, which contradicts 
\cite{sz1} 1.7.   \\
(e) Only one possibility remains: $\Delta$ has a factor $\Gamma$ of type ${\rm G}_2$,
a factor $\Psi$ of dimension $8$, and a $3$-dimensional factor $\Omega$. If $\Gamma$ is
compact, then $\cF_\Gamma$ is flat and $\Psi{\,\cong\,}\Psi|_{\cF_\Gamma}$ would be trivial. Hence $\Gamma{\,\circeq\,}\Gtwo(2)$. In a similar way as before, this leads to a contradiction: 
$(\Gamma\Omega)_p{\,\ne\,}\1$, and there is a subplane 
$\cE{\,=\,}\cE^\Psi{\,<\,}\cP$. If $\Psi$ is sharply transitive on $\cE\sm(au\smcup av\smcup uv)$, then a maximal compact subgroup $\Phi$ of $\Psi$ is homeomorphic to $(\HH')^2$ by the Mal'cev-Iwasawa theorem, but $\dim\Phi{\,\in\,}\{3,4,8\}$. Hence there is some point $q{\,\in\,}\cE$ such that 
$\Psi_{\hskip-1pt q}{\,\ne\,}\1$, and $\Gamma$ acts almost faithfully on 
$\cH{\,=\,}\cF_{\Psi_{\hskip-1pt q}}$, but $\dim\Gamma|_\cH{\,\le\,}11$. \qed
\par\medskip
{\bf 9.2  Lie.} {\it Assume that $\Delta$ fixes a triangle $a,u,v$\/}. \\
(A) {\it If $\Delta$ is semi-simple and if $\dim\Delta{\,\ge\,}22$, then $\Delta$ is a Lie group\/}. \\
(B) {\it If $\Delta$ has a $1$-dimensional compact normal subgroup $\Theta$  and if 
 $\dim\Delta{\,\ge\,}18$, then $\Delta$ is a Lie group and $\dim\Delta{\,\le\,}23$\/}. \\
(C) {\it If $\Delta$ has a minimal normal vector subgroup $\Theta{\,\cong\,}\RR^t$  and if 
$\dim\Delta{\,\ge\,}22$, then $\Delta$ is a Lie group\/}. 
\par\smallskip
Part (A) {\tt Proof.} Because of the previous Remark and 9.0 we may assume that 
$\dim\Delta{\,\le\,}24$. Suppose that $\Delta$ is not a Lie group, und use the same notation as in 9.0. Let $\Gamma$ be a proper factor of $\Delta$ of minimal dimension and denote the product of the 
other factors by $\Psi$.  Recall from the proof of 9.0 and the Remark added to 9.0 that 
$\Delta{\hskip1pt:\hskip1pt}\Delta_p{\,\le\,}14$ and that $\Delta_p\smcap\Nu{\,=\,}\1$ for any point
$p{\,\notin\,}au\smcup av\smcup uv$.  \\
(a) If $p^\Gamma{\,=\,}p$, then $\Gamma|_{p^\Delta}{\,=\,}\1$. The action of $\Nu$ on 
$\cF_\Lambda$ implies $\Lambda{\,\not\cong\,}\Gtwo$. Hence $\dim\Delta_p{\,\le\,}11$ and 
$\dim p^\Delta{\,\ge\,}11{\,>\,}8$, but then $\Gamma{\,=\,}\1$. Therefore $p^\Gamma{\,\ne\,}p$. \\
(b) {\it  $\dim\Gamma{\,\ge\,}6$\/}: we have 
$\cE{\,=\,}\langle p^{\Gamma\Nu},a,u,v\rangle\ledot\!\cP$ by  \cite{cp} 32.21 and 71.2. 
Obviously, $\Psi_p|_\cE{\,=\,}\1$ and $\Psi_p$ is a Lie group. Hence $\dim\Psi_p{\,\le\,}3$,\;  
$\dim\Psi{\,\le\,}17$, and $\dim\Gamma{\,\ge\,}5$.   \\
(c) Suppose in steps (c--f) that $\dim\Gamma{\,\ge\,}23$. Then
 $\dim\Delta{\,=\,}24$, or $\dim\Gamma{\,=\,}8$,\; $\dim\Psi{\,=\,}15$, and 
$\Psi$ is almost simple. In both  cases, $\dim\Gamma{\,>\,}6$. \\
(d)  If  $\dim\Gamma{\,=\,}8$, then again $\cE{\,=\,}\langle p^{\Gamma\Nu},a,u,v\rangle\ledot\!\cP$ 
and $\Psi_p|_\cE{\,=\,}\1$. Stiffness 2.6(\^b) shows that $\dim\Psi_p{\,\in\,}\{1,3\}$,\, 
$\cE{\,\ldot\,}\cP$, and $\langle p^\Psi\rangle{\,=\,}{\bold\cP}$. Consequently, 
$(\Gamma\Nu)_p{\,=\,}\1$,\,  $p^\Gamma$ is open in the point space of $\cE$, and $\Nu$ would be a Lie group by \cite{cp} 53.2 contrary to the assumption. \\
(e) {\it $\Delta$  is almost simple\/}. The only other possibility is $\dim\Gamma{\,=\,}10$. In this case $\Psi$ is a Lie group of type ${\rm G_2}$, and $\Psi$ is simple or a twofold covering of 
a simple group. Each non-central involution of $\Psi$ acts non-trivially on all three sides of the fixed triangle. Hence $\Psi$ contains a planar involution $\beta$, and 
$\Gamma|_{\cF_\beta}{\,\ne\,}\1$, but then $\dim\Gamma{\,\le\,}9$ by \cite{sz9} 6.1. \\ 
(f) Because of step (e),  the group $\Delta$ maps onto $\PSU5(\CC,r)$ with $r{\,=\,}1$ 
or $r{\,=\,}2$, and $\Delta$ has a subgroup $\Phi{\times}\Chi$, where $\Phi{\,\cong\,}\SU3\CC$ and 
$\Chi{\,\cong\,}\SU2(\CC,2{-}r)$. If $\iota$ denotes the involution ${\rm diag}(1,1,1,-1,-1)$, then
$\Phi{\times}\Chi{\,\le\,}\Cs\Delta\iota$ and $\dim\Cs\Delta\iota{\,\le\,}9{+}4$. Consequently, $\iota$ 
is not in the center of $\Delta$ and $\iota$ acts non-trivially on the sides of the fixed triangle. 
Therefore $\iota$ is planar. As $\dim\Phi{\,=\,}8$, Stiffness implies that $\Phi|_{\cF_\iota}{\,\ne\,}\1$. 
An involution $\kappa{\,\in\,}\Phi$ is not in the center of $\Phi$ and induces a Baer involution on 
$\cF_\iota$ for analogous reasons as those which proved the planarity of $\iota$. Therefore 
$\cC{\,=\,}\cF_{\iota,\kappa}$ is a $4$-dimensional subplane. Note that $\Nu$ acts freely on $\cC$. 
By \cite{cp} 53.2, the group $\Nu$ is a Lie group and so is $\Delta$. \\
(g) Now let $\dim\Delta{\,=\,}22$. Then $\Delta$ has no almost simple factor of dimension $16$, 
or else $\dim\Gamma{\,=\,}6$, both $\Gamma$ and $\Psi$ are Lie groups, and so is 
$\Delta{\,=\,}\Gamma\Psi$. (Note that the maximal compact subgroups of $\Gamma$ and $\Psi$ 
are almost simple, and use the approximation theorem \cite{cp} 93.8 together with 93.11.)
If $\Psi$ is almost simple of type ${\rm G}_2$, then 
$\Psi{\,\circeq\,}{\rm G}_2(2)$. As in step (e), a non-central involution $\beta{\,\in\,}\Psi$ is planar. 
There is a compact $6$-dimensional group $\Phi{\,\le\,}\Cs\Psi\beta$, and 
$\dim\Gamma\Phi|_{\cF_\beta}{\,\ge\,}8{+}3$. This contradicts \cite{sz1} 7.3. \\
(h) Thus $\dim\Gamma{\,=\,}6$ and $\Psi$ is a product of two almost simple factors $\Chi$ and 
$\Ypsilon$, where $\dim\Ypsilon{\,\i\,}\{8,10\}$. Again a non-central involutiom $\beta{\,\in\,}\Ypsilon$ 
is planar; moreover,  $\dim\Chi{\,\ge\,}6$, and $\dim\Gamma\Chi\Nu|_{\cF_\beta}{\,\ge\,}12$, but then 
$\Nu$ is a Lie group by \cite{pw3}. \\
Part  (B) {\tt Proof.}  (a) Suppose that $\Delta$ is not a Lie group. As in the proof of 9.0,
all orbits on two sides of the fixed triangle 
have dimension~${<\,}8$, and $\dim p^\Delta{\,\le\,}14$ for each point 
$p{\,\notin\,}au{\smcup}av{\smcup}uv$. Put $\Lambda{\,=\,}(\Delta_p)^1$.
By the approximation theorem 
\cite{cp} 93.8, there is a compact $0$-dimensional central subgroup $\Nu{\,\triangleleft\,}\Delta$
such that $\Delta/\Nu$ is a Lie group. Note that $\Theta{\,\le\,}\Cs{}\Delta$.
If $p^\Theta{\,=\,}p$, then $\cF_\Theta{\,=\,}\cF_{\hskip-1pt\Theta}^{\:\Delta}{\,<\,}\cP$ and 
$\Delta{\hskip1pt:\hskip1pt}\Delta_p{\,=\,}\dim p^\Delta{\,\le\,}\dim\cF_\Theta{\,=\,}d{\,\in\,}\{0,2,4,8\}$. Either $d{\,\le\,}2$,\,  $\Theta{\,\triangleleft\,}\Lambda{\,=\,}(\Delta_p)^1{\,\not\cong\,}\Gtwo$,\, 
$\dim\Lambda{\,\le\,}11$,  and $\dim\Delta{\,\le\,}13$, or $d{\,\ge\,}4$,\; $\dim\Lambda{\,\le\,}8$, 
and $\dim\Delta{\,\le\,}16$. If $p^\Theta{\,\ne\,}p$, then $\Theta$ acts non-trivially on the subplane 
$\cE{\,=\,}\langle p^{\Theta\Nu},a,u,v \rangle$, $\cE$ is not flat by \cite{cp} 32.17, and 
$\dim\cE{\,=\,}4$ or $\cE{\,\ledot\!}\cP$. In the first case, $\Nu|_\cE{\,=\,}\Nu/\Kappa$ is a Lie group by \cite{cp} 71.2, and $\Kappa{\,\ne\,}\1$. 
Hence $\cF_\Kappa{\,=\,}\cF_\Kappa^{\hskip2pt\Delta}{\,<\,}\cP$,\,
$p^\Delta{\,\subseteq\,}\cF_\Kappa$,\, $\Delta{\hskip1pt:\hskip1pt}\Delta_p{\,\le\,}8$, 
$\Delta_p|_\cE{\,=\,}\1$,\, $\dim\Delta_p{\,\le\,}8$, and $\dim\Delta{\,\le\,}16$. 
In the second case, we have $\Nu_p|_{\langle p^\Delta\rangle}{\,=\,}\1$. Therefore 
$\Delta{\hskip1pt:\hskip1pt}\Delta_p{\,=\,}\dim p^\Delta{\,\le\,}8$ or $\Delta_p\smcap\Nu{\,=\,}\1$
and $\Delta_p$ is a Lie group. Stiffness implies $\dim\Delta_p{\,\le\,}7$ or $\dim\Delta_p{\,\le\,}3$, 
respectively, and $\dim\Delta{\,\le\,}17$. \\
(b) If $\Delta$ is transitive on the complement of the triangle, then 
$\Delta'{\,\cong\,}\Spin8\RR$ by Corollary 2.18, $\dim\Delta{\,\ge\,}28$, $\Delta$ is a Lie group, 
the torus rank $\rk\Delta$ is $4$, and $\Delta$ cannot have a normal torus subgroup. Hence there is some point $p$ such that  $p,a,u,v$ form  a quadrangle and
$\Delta{\hskip1pt:\hskip1pt}\Delta_p{\,<\,}16$. Put $\Lambda{\,=\,}(\Delta_p)^1$ and 
$\cE{\,=\,}\cF_\Lambda$. If $\Theta|_\cE{\,=\,}\1$, then $\Theta{\,\le\,}\Lambda$, 
$\cE{\,\le\,}\cF{\,=\,}\cF_\Theta{\,=\,}\cF^\Delta{\,<\,}\cP$, 
$\Delta{\hskip1pt:\hskip1pt}\Delta_p{\,=\,}\dim p^\Delta{\,\le\,}8$, $\Lambda{\,\not\cong\,}\Gtwo$, 
$\dim\Lambda{\,\le\,}11$, and $\dim\Delta{\,\le\,}19$. If $\Theta|_\cE{\,\ne\,}\1$, however, then 
$\cE{\ledot\!}\cP$ or $\Theta|_\cE$ is a Lie group by \cite{cp} 32.21 and 71.2, and
$\Theta|_\cE{\,\cong\,}\TT$.  From \cite{cp} 32.17 
it follows that $\cE$ is not flat, and $\dim\cE$ is at least $4$. Stiffness implies $\dim\Lambda{\,\le\,}8$ and $\dim\Delta{\,<\,}24$.   \qed 
\par\smallskip
Part  (C) {\tt Proof.} (a) As in part (A) and with the same notation, $\dim p^\Delta{\,\le\,}14$ for each point $p{\,\notin\,}au{\smcup}av{\smcup}uv$. Select $p$ in such a way that $p^\Theta{\,\ne\,}p$
and choose a one-parameter subgroup $\Pi{\,\le\,}\Theta$ and an element $\rho{\,\in\,}\Pi$ with
$p^\rho{\,\ne\,}p$. The action of $\Delta$ on $\Theta$ being linear, we have
$\Theta_p{\,\le\,}\Gamma{\,=\,}\Delta_{p,\rho}{\,\le\,}\Cs{}\Pi\Nu$ and 
$\dim\Gamma{\,\ge\,}\dim\Delta{\hskip1pt-\hskip1pt}14{\hskip1pt-\hskip1pt}t$. As $p^\Pi{\,\ne\,}p$, 
the plane $\cE{\,=\,}\langle p^{\Pi\Nu},a,u,v\rangle$ is a connected subplane. 
By the Remark added to 9.0, the central subgroup $\Nu$ satisfies $\Delta_p\smcap\Nu{\,=\,}\1$. 
Therefore $\Delta_p$ and $\Gamma$ are Lie groups and $\Nu$  acts effectively on $\cE$.
Because of \cite{cp} 32.21 and 71.2 it follows that $\cE{\hskip1pt\ledot\hskip-3pt}\cP$. Note that 
$\Gamma|_\cE{\,=\,}\1$. Stiffness shows that $\Gamma$ is compact. In particular, $\Theta_p$ is compact and hence $\Theta_p{\,=\,}\1$. Moreover, $\Gamma^1{\,\le\,}\Spin3\RR$ and 
$\dim\Gamma{\,\le\,}3$. This proves the claim for $t{\,<\,}5$.  If $t{\,=\,}5$, then 
$\dim\rho^{\Delta_p}{\,=\,}\Delta_p{:\hskip1pt}\Gamma{\,=\,}5$ and $\dim\Delta_p{\,=\,}8$. 
As $\Theta_p{\,=\,}\1$, there is no restriction for the choice of $\rho$. Hence 
$\Delta_p$ is transitive on $\Theta\sm\1$, but a transitive linear group on $\RR^5$ 
contains $\SO5\RR$ and has dimension at least $11$, see \cite{Vl} Satz~1 
or use \cite{cp} 96.19--22. \\
(b) Let $\Lambda{\,=\,}(\Delta_p)^1$ and consider a  minimal $\Lambda$-invariant  subgroup
$\Eta{\,\cong\,}\RR^s$ of $\Theta$. According to step (a), we may assume that $s{\,>\,}5$.  
By definition, $\Lambda$ induces an irreducible re\-presentation 
$\overline\Lambda$ on $\Eta$; the kernel of this representation is contained in~$\Gamma$. Consequently  $\dim\Lambda{\,\ge\,}\dim\Delta{\hskip1pt-\hskip1pt}14{\,\ge\,}8$ and 
$\dim\overline\Lambda{\,\ge\,}5$.  
A Levi complement $\overline\Upsilon$ in $\overline\Lambda$ is the image of a maximal semi-simple  subgroup $\Upsilon$ of $\Lambda$, see \cite{cp} 94.27. \\
(c) In the case $s{\,=\,}6$ and $\dim\Lambda{\,>\,}8$ it follows that $\dim\rho^\Lambda{\,=\,}6$ for each choice of $\rho$ in~$\Eta$. Hence $\overline\Lambda$ is a transitive linear group of $\RR^6$, and $\Lambda$ is transitive on the $5$-sphere consisting of the rays in~$\RR^6$.  
From \cite{cp} 96.\hskip1pt19,21,22 and 94.27 we conclude that $\Upsilon{\,\cong\,}\SU3\CC$\; 
(since $\dim\Ypsilon{\,\le\,}14$).
In particular $\rk\Lambda{\,\ge\,}2$, and $\Lambda$ contains two commuting planar involutions
$\alpha,\beta$ such that $\dim\cF_{\alpha,\beta}{\,=\,}4$. As $\Nu$ acts freely on
$\cF_{\alpha,\beta}$,\; $\Nu$ would be a Lie group by \cite{cp} 71.2, a contradiction to our assumption. More generally, the last arguments prove the following \\
(d) {\tt Lemma.} {\it If there exists a pair of commuting involutions in $\Delta_p$, then $\Delta$ 
is a Lie group\/}. \\
(e) The possibility  $s{\,=\,}6$ and $\dim\Lambda{\,=\,}8$ will be postponed to step (h). 
If $s{\,=\,}7$, then Clifford's Lemma \cite{cp} 95.5 implies that $\overline\Upsilon$ is 
almost simple of dimension $\dim\overline\Upsilon{\,>3}$, and representation theory \cite{cp} 95.10 shows that $\overline\Upsilon$ 
is a group of type ${\rm G_2}$ and torus rank 2. Hence $\Delta$ is a Lie group by Lemma (d). \\
(f) Similarly, the other prime numbers can dealt with: if $s{\,=\,}11{\;\rm or\;}13$, then 
$\overline\Upsilon{\,\circeq\,}\SO s\RR$\; (again by \cite{cp} 95.10), and $\overline\Ypsilon$ is by far too large. \\
(g)  For compound $s$ the group $\overline\Upsilon$ may be properly semi-simple. If $s$ is even, complex representations are possible, and we can only infer $\dim\overline\Upsilon{\,\ge\,}3$. 
The kernel $\Kappa{\,=\,}\Cs\Upsilon\Eta$ of the action of $\Upsilon$ on $\Eta$ is contained 
in $\Gamma$; hence $\Kappa$ is compact and $\dim\Kappa{\,\le\,}3$. In fact, 
$\Kappa^1{\,\cong\,}\Spin3\RR$ or $\Kappa$ is finite (note that $\Kappa$ is semi-simple). As proper covering groups of $\SL2\RR$  do not have a faithful linear representation, each factor of 
$\overline\Upsilon$ has positive torus rank. The same holds for the factors of $\Ypsilon$.
Because of Lemma (d) we may assume that $\overline\Ypsilon$ is almost simple and that $\Kappa$ 
is finite. If $\rk\Ypsilon{\,>\,}1$ or if $\Ypsilon$ has a subgroup $\SO3\RR$,  Lemma (d) applies also. Thus only the cases 
$\overline\Ypsilon{\,\circeq\,}\SL2\CC,\ \SL3\RR$, and $\Sp4\RR$ need further consideration.  \\
(h) In the first case, $\Ypsilon{\,\cong\,}\SL2\CC$ contains a central involution $\iota$, which is planar because it fixes $p$. Let $\Ypsilon^*{\,=\,}\Ypsilon|_{\cF_\iota}$. By Stiffness, $\Ypsilon^*{\,\ne\,}\1$, 
and then $\dim\Ypsilon^*{\,=\,}6$; on the other hand, the stiffness property 
\cite{sz1} 1.5(1) implies 
$\dim\Ypsilon^*{\,\le\,}4$. For similar reasons, $\dim\Ypsilon{\,\ne\,}8$: the central involution $\kappa$ 
of the simply connected two-fold covering $\Ypsilon$ of $\SL3\RR$ is again planar, Stiffness shows 
$\Ypsilon|_{\cF_\kappa}{\,\cong\,}\SL3\RR$, but $\dim\Ypsilon|_{\cF_\kappa}{\,\le\,}4$ by   
\cite{sz1} 1.5(1). 
Finally, let $\dim\Ypsilon{\,=\,}10$. Each faithful linear representation $\overline\Ypsilon$ of 
$\Ypsilon$ is isomorphic to ${\rm(P)}\Sp4\RR$\; (see\cite{cp} 95,10), and $\U2\CC$ is a maximal compact subgroup of  $\Sp4\RR$. By assumption, $\Ypsilon$ is a finite covering of $\overline\Ypsilon$. Hence 
$\rk\Ypsilon{\,=\,}2$, and Lemma (d) applies. \qed
\par\medskip
{\bold Corollary.} {\it If $\cF_\Delta$ is a triangle and if $\dim\Delta{\,\ge\,}22$, then $\Delta$ is a Lie group\/}.
\par\smallskip
{\tt Proof.} There is a compact $0$-dimensional central subgroup $\Nu$ such that $\Delta/\Nu$ 
is a Lie group, and $\Delta/\Nu$ is semi-simple or has a central torus $\Tau$ or a minimal normal vector group $\Xi$, see 2.3. In the first case, each commutative connected normal subgroup 
$\Alpha{\,\triangleleft\,}\Delta$ maps onto the identity of $\Delta/\Nu$. Hence $\Alpha{\,\le\,}\Nu$, 
and $\Alpha$ is trivial, in other words, $\Delta$ is semi-simple and (A) applies. In the other two cases, 
the connected component $\Theta$ of the pre-image of $\Tau$ or $\Xi$ satisfies the assumptions in 
(B) or (C) respectively. \qed
\par\medskip 
{\bf 9.3.} {\it  If $\dim\Delta{\,=\,}30$ and if $\Delta$ fixes a triangle, then 
$\Delta{\,=\,}\Rho{\times}\Phi$, where $\Phi{\,\cong\,}\Spin8\RR$ and $\Rho$ is a product of two 
one-parameter groups of homologies with distinct centers\/}.
\par\smallskip
{\tt Proof.} Let $c{\,\in\,}av\sm\{a,v\}$ and $z{\,\in\,}S{\,=\,}uv\sm\{u,v\}$, and put 
$\Gamma{=\,}\Delta_c$ 
and $\Lambda{\,=\,}\Gamma_{\hskip-1.5pt z}$. Then $\dim z^\Gamma{\,=\,}8$ and $\Gamma$ is 
transitive on $S$.  As $\Lambda{\,\cong\,}\Gtwo$,  it 
follows by 2.17(e) that $\Gamma{\,\cong\,}e^{\RR}{\times}\Spin7\RR$. Similarly, 
2.17(f)  shows that a maximal compact subgroup $\Phi$ of $\Delta$ is isomorphic to $\Spin8\RR$. 
The $2$-dimensional radical $\Rho{\,=\,}\sqrt\Delta$ acts effectively on the flat plane 
$\cF_\Lambda$, and $\Rho{\,\cong\,}\RR^2$ by \cite{cp} 33.10. Therefore 
$\Rho_{\hskip-2pt z}{\,\le\,}\Cs{}\Gamma$,\: $\Rho_{\hskip-2pt z}|_{z^\Gamma}{\,=\,}\1$, and 
$\Rho_{\hskip-2pt z}{\,\le\,}\Delta_{[a,uv]}$.  Analogously, $\Rho_{\hskip-2pt c}{\,\le\,}\Delta_{[u,av]}$. \qed
 \par\bigskip
\break 
{\Bf 10. Summary}
\par\medskip
The entries in the table below have the following meaning:\quad  
\par\smallskip
\begin{tabular}{l l}
$\dim\Delta\ge b \Rightarrow \cP$ is known, &
$\dim\Delta\ge b' \Rightarrow \cP$ is a translation plane, \\
$\dim\Delta\ge b''\Rightarrow\cP$ is a  Cartesian plane,\hspace{30pt} &
$\dim\Delta\ge b^*\Rightarrow \cP$ is a Hughes plane, \\
$\dim\Delta\ge c \Rightarrow \cP$ is classical (Moufang), &
$\dim\Delta\le d$, \   $\dim\Delta\ge g \Rightarrow \Delta$ known.
\end{tabular}
\par\smallskip 
\begin{center}
\begin{tabular}{|c||c|c|c|c||l|} 
\hline
$\cF_\Delta$ & $\Delta$ $s$-$s$  & $\TT\triangleleft \Delta$ & $\RR^t\triangleleft\Delta$   & 
$\Delta$ arbitrary & References   \\  \hline 
$\emptyset$ & $b^*{=}25$ & $b^*{=}23$ & $d{\,\le}23$ & $b{\,=}25$\quad$c{\,=}37$ 
&  3.1,2,3,\;\cite{cp}\,86.35  \\
$\{W\}$ & $d{\,=\,}21\hskip5pt$ & $4.4\hskip3pt^{1)}$ & $b'{\,=}35$ & $b{\,=}35$\quad$c{\,=}36$  
& 4.3,6,8,11,12   \\ 
flag & $d{\,=}21\,^{2)}$ & $5.3\hskip3pt^{1)}$ & ? & $b{\,=}40$\quad$c{\,=}41$ 
& 5.2,\,\cite{cp}\,87.7  \\  \hline
$\{o,W\}$ & $g{\,=}29 \ $ & $d{\,=}30$ & $d{\,=\,}32\,^{7)}$ &\hfill $c{\,=}38$ 
& 4.2,\,5.\,6/12 \\  
$\langle u,v\rangle$ & $d{\,=}21\,^{2)}$ & $d{\,=}20 \ $ & $b{\,=}35\hskip3pt^{3)} $ & $b{\,=}35$\quad
$b'{=}39$ & 6.2,3,\,5/7 \\
$\langle u,v,w\rangle$ & $d{\,=}18 \ $ & $d{\,=}16 \ $ &$g{\,=}32\quad c{\,=}33$ & $b'{=}32$\quad$c{\,=}33 $ & 7.2,3,4 \\  \hline
$\langle u,v,ov\rangle$ & $d{\,=}24\,^{4)}$ & $d{\,=}30 \ $ & $b''{=}33 \ $ & $b{\,=}33$\quad 
$c{\,=}38$  &  8.2,3,\hskip1pt4/6 \\  
$\langle o,u,v\rangle$ & $g{\,=}28 \ $ & $d{\,=}23 \ $ & ? & $d{\,=}30$\quad$g{\,=}30$ 
&  9.\,1,\hskip1pt2,\hskip1pt3 \\
arbitrary & $b^*{=}29\;^{5)}$ & $b^*{=}31$ & $b{\,=}35\,^{6)}$ & $b{\,=}40$\quad$c{\,=}41$ & 
\cite{sz5},\,\cite{cp}\,87.7 \\  \hline
\end{tabular}
\end{center}
\par 
\begin{tabular}{l l l l}
\quad $^{1)}$ & 
$\dim\Delta{\,\ge\,}18\Rightarrow\exists_{\cB{\,\cong\,}\cH}\,\cB^\Delta{\,=\,}\cB{\,\ldot}\cP$ &
\quad $^{2)}$ & $\dim\Delta{\,\le\,}16$ {\it if $\Delta$ is almost simple\/} \cr 
\quad $^{3)}$ & $\dim\Delta{\,\in\,}\{33,34\}\Rightarrow\cP\ is\ transl.\ plane$  &  
\quad $^{4)}$ & {\it if $\Delta$ is a Lie group\/} \cr 
\quad $^{5)}$ & or $\dim\Delta{\,\ge\,}36$ &
\quad $^{6)}$ & {\it or $\cF_\Delta$ is a flag\/} \cr
\quad $^{7)}$ & {\it if $t{\,\ne\,}1$\/} &   \cr
\end{tabular}
\par\medskip
{\bf Lie groups.} {\it If $\Theta$ denotes a compact conected $1$-dimensional subgroup and
if $\dim\Delta{\,\ge\,}k$, then $\Delta$ is a Lie group\/}.
\begin{center}
\begin{tabular}{|c||c|c|c|c||l|} 
\hline
$\cF_\Delta$ & $\Delta$ $s$-$s$  & $\Theta\triangleleft \Delta$ & $\RR^t\triangleleft\Delta$   & 
$\Delta$ arbitrary & References   \\  \hline 
$\emptyset$ & $21$ &  && 23 &  3.1, 3.0  \\
$\{W\}$ & $14$ &  && 23 & 4.3, 4.0 \\
flag & $15$ &  && 22 & 5.1, 5.0 \\  \hline
$\{o,W\}$ &26& && 27 &   \cite{psz}\hskip1pt({\bold  a}), \cite{psz}  \\  
$\langle u,v\rangle$ & $14$ & && 18 & 6.1, 6.0 \\
$\langle u,v,w\rangle$ & $14$ & && 15 & 7.1, 7.0 \\  \hline
$\langle u,v,ov\rangle$ & $26$ & $20\,^{1)}$ && 27 & 8.0, \cite{psz}  \\  
$\langle o,u,v\rangle$ & 22 & $18\,^{\ \ }$  & 22 & 22 & 9.\hskip1pt2 \\
arbitrary & 26 &  && 27&  \cite{psz} \\  \hline
\end{tabular}
\end{center}
\par 
\begin{tabular}{l l l l}
\quad $^{1)}$ & or $\Cs{}\Delta|_{av}{\,=\,}\1$\; (8.3 Remark) \cr
\end{tabular}
\par

\qquad Helmut R. Salzmann \par
\qquad Mathematisches Institut \par
\qquad Auf der Morgenstelle 10 \par
\qquad D-72076 T\"ubingen \par
\qquad helmut.salzmann@uni-tuebingen.de 

\end{document}